\documentclass[11pt]{article}

\usepackage[numbers]{natbib}
\usepackage{enumitem}
\usepackage{amsmath}
\usepackage{amsfonts}
\usepackage{amssymb}
\usepackage{epsfig}
\usepackage{dsfont}
\usepackage{mathrsfs}
\usepackage{bbm}
\usepackage{bm}
\usepackage{multirow}
\usepackage{makecell}
\usepackage{threeparttable}
\usepackage{authblk}
\usepackage{mathabx} 
\usepackage{amsxtra}
\usepackage{graphicx}
\usepackage{animate}
\usepackage{tikz}
\usepackage{algorithmic}
\usepackage{float}  
\usepackage{lipsum}
\usepackage{caption}
\usepackage{soul}
\usepackage{url}
\usepackage{algorithm}
\usepackage[title]{appendix}
\usepackage{boldline}
\usepackage{diagbox}
\usepackage{xcolor}
\usepackage{pifont}
\usepackage{harmony}
\usepackage{scalerel}
\usepackage[colorlinks=true]{hyperref}
\hypersetup{urlcolor=blue, citecolor=blue,linkcolor=blue}
\usepackage[capitalise,nameinlink]{cleveref}

\RequirePackage{xcolor}[2.11]
\colorlet{inlinkcolor}{green!60!black}
\colorlet{exlinkcolor}{red!50!black}
\colorlet{reviewcolor}{black!50}
\crefname{equation}{}{}
\crefname{figure}{Figure}{Figures}
\crefname{assumption}{Assumption}{Assumptions}
\crefname{subappendix}{Appendix}{Appendices}

\allowdisplaybreaks[4]
\usepackage{subfigure}
\counterwithin{algorithm}{section}
\newtheorem{theorem}{Theorem}[section]

\newtheorem{proposition}{Proposition}[section]
\newtheorem{example}{Example}[section]
\newtheorem{remark}{Remark}[section]
\newtheorem{definition}{Definition}[section]
\newtheorem{assumption}{Assumption}[section]

\newtheorem{criterion}{Criterion}[section]
\numberwithin{equation}{section}
\numberwithin{definition}{section}

\newcommand{\dd}{\mathsf {d\kern -0.07em l}} 
\newcommand{\bgeqn}{\begin{eqnarray}}
\newcommand{\edeqn}{\end{eqnarray}}
\newcommand{\bgeq}{\begin{eqnarray*}}
\newcommand{\edeq}{\end{eqnarray*}}
\newcommand{\bec}{\begin{center}}
\newcommand{\enc}{\end{center}}
\newcommand{\R}{{\rm I\!R}}

\newcommand{\inmat}[1]{\mbox{\rm {#1}}}

\renewcommand{\Box}{\framebox{\rule{0.3em}{0.0em}}}

\newcommand{\be}{\begin{equation}}
\newcommand{\ee}{\end{equation}}

\newcommand{\vt}{{\vartheta}}

\newcommand{\bdy}{\bm{y}}

\newcommand{\proof}{\noindent{\bf Proof. }}

\def\bbe{{\mathbb{E}}}

\title
{Modified Polyhedral Method to Elicit Shape-Free Utility and Reduction of Conservatism in Robust Optimization
\footnote{This project is supported by a CUHK startup grant and RGC grant
14204624.}
}

\author[1]{Sainan Zhang}
 \author[2]{Shaoyan Guo}
 \author[3]{Melvyn Sim}
 \author[1]{Huifu Xu}
 \affil[1]{Department of Systems Engineering and Engineering Management,
              The Chinese University of Hong Kong, Shatin, N.T., Hong Kong.}
\affil[2]{School of Mathematical Sciences, Dalian University of Technology, Dalian, China.}
\affil[3]{Department of Analytics and Operations, NUS Business School, National University of Singapore, Singapore.}              
\affil[]{\textsf{snzhang.m@gmail.com,
              syguo@dlut.edu.cn,
              dscsimm@nus.edu.sg,
              hfxu@se.cuhk.edu.hk}}

\begin{document}
  \maketitle
\vspace{-3em}
\begin{abstract} 
In this paper, we propose a modified polyhedral method to elicit a decision maker’s (DM’s) nonlinear univariate utility function without explicit knowledge of its shape, Lipschitz modulus, or inflection point. Inspired by \cite{THS04} for linear multivariate utility elicitation, we extend the polyhedral method to the nonlinear setting by addressing two major challenges. First, we approximate the nonlinear utility function using a continuous piecewise linear function (PLF) represented by its increment vector, such that utility elicitation can be viewed as progressively reducing a polyhedral feasible region of increment vectors. Second, we iteratively shrink this polyhedron through hyperplane cuts generated from pairwise lottery comparisons. Since PLF approximation errors may cause the resulting hyperplanes to deviate from ideal cutting directions, we dynamically enrich the PLF by incorporating newly generated lottery outcomes as additional breakpoints. The elicited responses are then used to construct an ambiguity set of utility functions for decision making under worst-case expected utility. To mitigate the potential conservatism of the resulting maximin model, we use the analytic center of the ambiguity set as a nominal estimate of the DM’s utility increment vector and derive, through dual reformulation, a flexible weighting scheme for controlling conservatism. The resulting model accommodates asymmetric bounds on utility increments around the analytic center and admits a single mixed-integer linear programming reformulation for piecewise linear utilities. Numerical results demonstrate the effectiveness of the proposed elicitation method and conservatism-reduction approach.
Finally, we discuss several potential extensions, including applications to multivariate nonlinear utility elicitation and dynamic pricing problems.
\end{abstract}

\textbf{Keywords}: 
Polyhedral method,
adaptive preference elicitation,
shape-free utility function,
conservatism reduction,
piecewise linear function approximation.

\vspace{-0.15cm}

\section{Introduction}

\vspace{-0.15cm}

The theory of expected utility, established by \cite{VNM47}, 
provides a fundamental framework for decision making under uncertainty.
However, in practice, 
the utility function that characterizes a decision maker's (DM's) preferences over uncertain outcomes is typically unknown.
Consequently, various methods for eliciting preference 
have been proposed 
to estimate the DM's utility function from observed or elicited preference information. When such information is limited, the estimated utility function may be subject to substantial uncertainty, potentially leading to poor decisions.
To account for this uncertainty,
Armbruster and Delage \cite{AmD15} propose a preference robust optimization (PRO) model that determines the optimal decision based on the worst-case utility function derived from an ambiguity set, 
a collection of utility functions constructed from the available preference information.
A primary concern with PRO models is
that the optimal decision may become excessively conservative
 when the ambiguity set is large.
This raises two related questions: 
how to elicit a DM's true utility efficiently to mitigate ambiguity while simultaneously 
modifying existing PRO frameworks to reduce conservatism.
This paper aims to address these challenges.

\vspace{-0.1cm}
\subsection{Literature Review}
\vspace{-0.1cm}

The literature related to our work can be divided into two streams: preference elicitation and PRO.

\vspace{-0.3cm}
\subsubsection{Preference Elicitation} 
\vspace{-0.1cm}

The literature can broadly be categorized into two main groups based on the structure of utility functions and their respective applications.
In marketing,
consumer preferences for a multi-attribute product or service
are often characterized by  an observable and deterministic  multi-variate linear utility (a scalar product of the vector of product features and the vector of consumer's evaluations of the features which is also known as preferences)
plus a random preference shock
\citep{Fis99}.
In practice,
information about
consumers' preferences is often incomplete,
and it is necessary to elicit preferences.
The well-known
{\em conjoint analysis} 
is a survey-based statistical 
market research technique 
that helps to determine the linear utility function,
e.g.,
by providing preference ratings for product attributes.
The most widely adopted method
is the
{\it choice-based conjoint analysis},
where there is a choice between two or more products.
Toubia et al.~\cite{THS04}  first 
proposed 
a polyhedral method to elicit
a deterministic vector of partworths for 
a linear multi-attribute utility function.
 Toubia et al.~\cite{THG07}
proposed a probabilistic polyhedral method which incorporates response errors
in the polyhedral method.
Bertsimas and O'Hair~\cite{BeO13} considered the case where inconsistent responses occur 
in the elicitation process using robust optimization. 
Saur\'e and Vielma~\cite{SaV19} considered the case where 
the vector of partworths
is random and proposed an approximate adaptive 
ellipsoidal method which incorporates a Bayesian learning scheme.

Another stream of the preference elicitation literature, particularly in behavioral economics and decision analytics, focuses on nonlinear univariate utility functions.
The well-known
Von Neumann–Morgenstein’s expected utility theory (EUT, see \cite{VNM47}), 
Yaari's dual theory of choice \citep{Yaa87_b},
rank-dependent expected utility theory (RDEU, see \cite{Qui82,Wak94}),
and cumulative prospect theory (CPT, see \cite{TvK92})
characterize preference in this manner.
In typical practice,
the specific nonlinear univariate utility function 
that captures a DM's preference is  
unknown.
Various methods have been proposed to determine the utility functions in EUT \citep{BDM64} and the utility and distortion functions in RDEU and CPT for parametric forms \citep{TvK92,CaH94,HeO94,TvF95} and parametric-free forms \citep{WaD96,Abe00,ABP07}.

Many utility elicitation methods under EUT assume that prior shape information, such as monotonicity, convexity/concavity, or the locations of inflection points, is 
available. Such information 
restricts the set of candidate utility functions 
substantially
by ruling out functions that violate the prescribed shape.
However, in practice, reliable shape information is often not
available. Even when the utility function exhibits structured nonlinearities, such as the concave--convex--concave pattern reported in \cite{FrS48}, the locations of the inflection points are typically unknown. This motivates 
us to develop an approach using a polyhedral elicitation
that progressively reduces the set of candidate utility functions
from a geometric perspective, without requiring assumptions about
the utility shape or the Lipschitz
modulus, or~the~locations~of~inflection~points.

\vspace{-0.1cm}
\subsubsection{PRO}
\vspace{-0.1cm}

Recent advances
in PRO can be classified into three categories: utility-based preference robust optimization (UPRO), which addresses ambiguity in the true utility function; preference robust risk optimization (PRRO), which considers ambiguity in the true risk measure; and preference robust choice function-based optimization models (PRCF), which encompass choice functions that may represent utility and/or risk functionals.
The literature on UPRO can be further classified based on the characteristics of the true utility function, including whether it is univariate or multivariate, deterministic or stochastic,
linear or nonlinear, parametric or non-parametric, single-stage or multistage, as well as whether elicitation errors are accounted for in the construction of the ambiguity sets.
Armbruster and Delage~\cite{AmD15}
considered nonlinear 
PRO problems,
where the {\em random (relative) utility split} (RRUS)
scheme 
is used to elicit preference information for constructing an ambiguity set of nonparametric univariate utility functions.
The RRUS scheme is also applied in \cite{GXZ22}, who constructed the ambiguity set using moment-type conditions.
They propose a piecewise linear approximation of the utility functions in the ambiguity set and derive error bounds for the model approximation.
Hu and Mehrotra  \cite{HuM15} studied
 risk-averse utility optimization,
where the ambiguity set of utility is defined by the bounds on the marginal utility and the utility value, and some auxiliary constraints 
using Lebesgue integrable functions. 
Haskell et al.~\cite{HFD16} proposed a robust model which handles the ambiguity of a DM’s beliefs and taste.  
Vayanos et al.~\cite{VYMDR20} considered a preference robust multi-attribute recommender model  where the utility functions are linear and the marginal
utility functions 
may not be 
non-decreasing.
Wu et al.~\cite{wu2025bi} proposed a UPRO model for bi-attribute decision making problems where the utility function can be nonlinear.
All the above studies focus on single-stage decision making, where a single decision is made before observing the underlying uncertainty.
Liu et al.~\cite{liu2025preference2} considered a multistage expected utility maximization problem where the DM’s utility function at each stage depends on the historical path and the information on the true utility function is incomplete.

A second stream of PRO research,
PRRO, considers
problems in which the DM's risk preferences are ambiguous. 
The literature on PRRO primarily focuses on single-attribute decision-making problems involving 
ambiguous convex risk measures or distortion risk measures 
(see \cite{DeL17,DGX22,WaX20,GuX22,li2021inverse,bernard2026preference2}).  
A third stream,
PRCF, characterizes the DM's preferences through choice functions, which provide a
broader than expected utility or dual expected utility theory \citep{HXH22,wu2025robust2}.

\vspace{-0.1cm}
\subsection{Contributions}
\vspace{-0.1cm}

\textbf{Objectives and Settings.}
In this paper, we propose a modified polyhedral method (MPM) to elicit
general nonlinear univariate utility functions 
through adaptive queries for pairwise comparisons.
To this end,
we approximate a general nonlinear utility function using a piecewise linear function (PLF). The PLF is parameterized by a vector of utility increments between consecutive breakpoints. This representation allows us to characterize the set of feasible utility functions
as a polyhedron of increment vectors in a Euclidean space.
We envisage the
adaptive selection of questions 
as a step-wise procedure
within the framework of
the polyhedral method.
Specifically, we compute the analytic center and the longest axis of the inner Sonnevend ellipsoid \citep{Son85}
of the polyhedron.
We then use the intersections of the longest axis with the boundary of the polyhedron to formulate two optimization problems whose solutions provide a query with a choice between
a risky lottery and a certain lottery. 
Based on the DM's response,
we construct a cutting hyperplane that reduces the current polyhedron.
We further incorporate the elicited preference information into PRO by treating the analytic center of the feasible polyhedron as a nominal utility increment vector and introducing a weighting scheme to adjust the level of conservatism flexibly.

In developing the adaptive query design for nonlinear utility functions, we make the following assumptions:
(a) the nonlinear function is non-decreasing,
so that its PLF is characterized by a vector of increments;
(b) a DM's choice is driven by EUT,
where the utility function is deterministic 
and no prior distribution for the vector of increments is available;  
(c) no parametric assumptions, and no information about the  shape 
and Lipschitz modulus of  utility are available;
(d) the number of queries in each query is not fixed 
and  each query comprises two lotteries.
The first assumption is to simplify the exposition.
The second assumption is standard in the behavioral economics literature.
The third assumption avoids a potential misspecification bias
that arises from imposing an incorrect parametric form or utility shape.
Moreover,
our method does not impose a prespecified shape on the underlying utility function, making it particularly suitable when the locations of inflection points are unknown.
Finally, the number of queries is not fixed in advance. In addition, new breakpoints can be introduced adaptively during the elicitation process.
The main contributions of this paper can be summarized as follows.

{\bf Methodology in Preference Elicitation}.
We extend the polyhedral elicitation framework of \cite{THS04} from linear multivariate utility functions to nonlinear univariate utility functions.
We do so by using a PLF  to approximate the nonlinear utility function
and then reformulate the former in a linear form with respect to (w.r.t.) the vector of increments.
Since the PLF approximation may cause the generated cutting hyperplanes to deviate from their ideal cutting directions, 
we develop a strategy which adaptively 
update the set of blackpoints of PLF (Algorithm~\ref{alg:RUS}).
Unlike the RRUS scheme used to elicit utility,
the new method performs well without prior information about 
structure of the true utility function such as shape, Lipschitz modulus and the inflection point (Section~\ref{sec:information}).
We investigate the role of adaptive breakpoints (Section~\ref{sec:advantage}), quantify the reduction of the polyhedra
and the associated ambiguity set of utility functions under different metrics (Section~\ref{sec:quality}), and propose a variant of Algorithm~\ref{alg:RUS} (Algorithm~\ref{alg:RUS-fixed}) that focuses on preference elicitation with a prespecified set of breakpoints (Section~\ref{sec:variant_fixed}). The latter points to a potential direction for theoretically analyzing the relationship between the number of queries and the accuracy of the method.

{\bf Modeling in UPRO}.
We propose a new conservatism-control framework for UPRO (Section~\ref{Sec:Numberial}). Specifically, we use the analytic center of the polyhedron constructed from elicited preference information 
as the nominal increment vector and derive, through a dual reformulation, a weighting scheme that flexibly adjusts the level of conservatism.
Unlike conventional robustness-control mechanisms based on symmetric deviations from a nominal point in \cite{BeS04}, our framework accommodates asymmetric bounds on utility increments around the analytic center and assigns different weights to deviations in different directions. When the utility functions in the ambiguity set are piecewise linear, the resulting model admits a single mixed-integer linear programming (MILP) reformulation.

{\bf Numerical Results and Potential Extensions}.
We have undertaken extensive numerical experiments
on the proposed MPM and the conservative reduction of the UPRO model. The results show the effectiveness of MPM particularly in absence of prior knowledge of the utility shape or Lipschitz modulus.
They also demonstrate the advantages of the new UPRO model when applied to a portfolio optimization problem, achieving a balance between robustness and conservatism.
Finally, we discuss several promising extensions of the 
MPM including 
a hybrid elicitation scheme that integrates a bisection strategy to further improve the performance of Algorithm~\ref{alg:RUS} (Section~\ref{sec:hybrid_Bisection}),  generalization of  MPM to accommodate multivariate nonlinear utility functions (Section~\ref{sec:multi}),
and an application of its core principles to dynamic pricing problems (Section~\ref{sec:dynamic_pricing}).

\vspace{-0.1cm}
\subsection{Organization of the Paper}
\vspace{-0.1cm}

The rest of the paper is organized as follows.
Section~\ref{sec:Preliminaries}
introduces preliminaries for the original polyhedral method,
the expected utility theory,
and the PLF of nonlinear utility.
Section~\ref{sec:models} formulates the polyhedron of utility increment vectors based on the DM's responses to pairwise comparisons
and constructs the corresponding ambiguity set and its analytic center.
Section~\ref{sec:MPM} develops polyhedral methods for univariate nonlinear utilities and 
derives its key properties. Section~\ref{sec:strategies} 
develops a modified polyhedral algorithm and investigates its efficiency.
Section~\ref{Sec:Numberial} 
investigates the tractable reformulations of the resulting preference robust model 
with an adjustable level of conservatism and reports numerical results.
Section \ref{sec:extension} discusses extensions of the proposed method to the elicitation of multiattribute utility 
and its application to dynamic pricing.
Finally, 
Section~\ref{sec:Conclusion} outlines potential directions for future research.
All proofs are relegated to the Electronic Companion.

\vspace{-0.1cm}
\subsection{Notation}
\vspace{-0.1cm}

For any integer $N\geq 1$, we write $[N]$ for set $\{1,\cdots,N\}$ and
use $|S|$ to denote the cardinality of a set $S$.
We use a boldface lowercase letter to denote a column vector
and a boldface uppercase (capital) letter to denote a matrix. 
To ease the exposition, we write $(x_1,\cdots,x_n)$ 
for an $n$-dimensional column vector ${\bm x}$ 
and $({\bm a}, {\bm b})$ for $({\bm a}^{\top}, {\bm b}^{\top})^{\top}$ where $\top$ denotes transpose.
Given any matrix ${\bm A} = (a_{ij})_{i\in[m],j\in[n]} \in \R^{m\times n}$,
we let ${\bm a}_i^{\top}$ be its $i$-th row vector.
Let ${\bm e}$ denote a vector of ones of appropriate dimension.
${\cal S}_n:=\{{\bm v}\in \R^{n}_+:{\bm e}^\top {\bm v}\leq 1\}$ is an $n$-simplex.
We use 
$\mathds{1}_{S}(x)$
to represent an indicator function of a set $S$, 
that is, $\mathds{1}_{S}(x)=1$
for $x\in S$ and $0$ otherwise.
 
\vspace{-0.15cm}
\section{Preliminaries}
\label{sec:Preliminaries}
\vspace{-0.15cm}

In this section, we 
briefly review the original polyhedral method and
the EUT, 
and then 
introduce piecewise linear 
utility functions.
These preliminaries will be used to develop the MPM for adaptive 
preference elicitation strategies
in Section~\ref{sec:strategies}.

\vspace{-0.1cm}
\subsection{Original Polyhedral Method
in Marketing}
\label{sec:OPMF}
\vspace{-0.1cm}

Our work is inspired by the polyhedral method proposed by 
\cite{THS04}.
Readers who are familiar with the polyhedral method
may move directly to
Section~\ref{sec:Preference functionals}.
The basic concept
of the method can be described as follows.
Consider choice-based conjoint analysis 
where respondents are presented
with two product profiles and asked to choose the profile they prefer.
The product profile is described
in terms of a binary vector ${\bm x}\in \R^n$ representing $n$ attributes and levels,
that is,
the feasible profiles 
lie within a subset of
$\{0,1\}^n$.
Assume that respondents'
preferences are characterized
by a vector of partworths ${\bm \beta}\in \R^n$ representing the relative importance assigned to product attributes.
Then the utility to product profile ${\bm x}\in  \{0, 1\}^n$
is denoted by 
\bgeqn 
U_{\bm x}:={\bm \beta}^{\top}{\bm x}+\epsilon_{\bm x},
\label{eq:Ux-market}
\edeqn
where  
$\epsilon_{\bm x}$ is a random idiosyncratic shock
to utility
and ${\bm \beta}$ is 
a vector of parameters
to be identified.
Partworth $\beta_i$ of attribute $i$  is included when $x_i=1$ and not included when $x_i=0$.
The polyhedral method introduced in
\cite{THS04}
begins by a
polyhedron-shaped 
region $\mathscr{U}_0:= \{{\bm \beta}\in \R^n:{\bm A}_0{\bm \beta} \leq {\bm b}_0\}$,
where ${\bm A}_0\in \R^{m\times n}$, ${\bm b}_0\in \R^m$.
$\mathscr{U}_0$ is the set which 
contains the true ${\bm \beta}^*$
and is constructed based on initially available information.
The polyhedral method 
sequentially designs pairs of product profiles
represented by 
${\bm x}^1$ and 
${\bm y}^1$
for the consumer to compare.
After observing consumer's choice, the polyhedron $\mathscr{U}_0$ shrinks 
to $\mathscr{U}_1$ 
by a cutting hyperplane while retaining the true parameter vector
${\bm \beta}^*$. This process is repeated until a prescribed stopping criterion is satisfied.
Let
$\mathscr{U}_{k-1}$ denote the
polyhedron 
after observing the 
answers
to
the previous $k-1$ questions.
Assume that there is no response error
(i.e., $\epsilon_{{\bm x}_k}={\epsilon}_{{\bm y}_k}$, see \cite[Section 4.1]{SaV19}).
Upon selecting the $k$-th profiles $({\bm x}^k,{\bm y}^k)$ and observing the answer,
the region 
is updated
by
\bgeqn
\label{eq:poly-cut-market}
\mathscr{U}_{k}
:=\mathscr{U}_{k-1}\bigcap 
\left\{
\begin{array}{ll}
& \{{\bm \beta}\in \R^{n}: {\bm \beta}^{\top}{\bm x}^k \geq {\bm \beta}^{\top}{\bm y}^k, \inmat{ if } {\bm x}^k\succeq {\bm y}^k\}\\
&  \{ {\bm \beta}\in \R^{n}: {\bm \beta}^{\top}{\bm x}^k\leq {\bm \beta}^{\top}{\bm y}^k, \inmat{ if } {\bm x}^k\preceq {\bm y}^k\}
\end{array}
\right.,
\edeqn
where ${\bm x}\succeq {\bm y}$ means that ${\bm x}$ is preferred to ${\bm y}$, or equivalently
$U_{\bm x}\geq U_{\bm y}$.
The overall procedure
is as follows.
Step 1. 
Compute the analytic center ${\bm c}$ and the longest axis of the polyhedron $\mathscr{U}_{k-1}$.
Step 2. 
Determine the two intersections ${\bm \beta}_1$ and ${\bm \beta}_2$ between the longest axis and the boundary of the polyhedron $\mathscr{U}_{k-1}$.
Step 3.
Select ${\bm x}^k$ and ${\bm y}^k$ by solving two problems:
\begin{subequations}
\begin{align}
& {\bm x}^k\in \arg\max \{{\bm \beta}_1^{\top}{\bm x}:{\bm x} \in \{0,1\}^n, {\bm c}^{\top}{\bm x} \leq D\},\label{eq:original_PM_a}\\
& {\bm y}^k\in \arg\max \{{\bm \beta}_2^{\top}{\bm x}:{\bm x} \in \{0,1\}^n, {\bm c}^{\top}{\bm x} \leq  D\},\label{eq:original_PM_b}
\end{align}
\end{subequations}
where $D>0$ is randomly generated 
until ${\bm x}^k\neq {\bm y}^k$
(see \cite{THS04,SaV19}).
Step 4. Present the query to the consumer. 
Update the current polyhedron according to the consumer's response
by retaining the half-space containing the true parameter.

\vspace{-0.1cm}
\subsection{VNM's Expected Utility} 
\label{sec:Preference functionals}
\vspace{-0.1cm}
 
Let $(\Omega,{\cal F},\mathbb P)$ be a probability space and let
${\cal L}^0$ denote the set of all random variables mapping from $\Omega$ to $\R$.
The preference functional based on the expected utility theory is defined by $V_u:{\cal L}^0\to \R$,
\bgeqn
\label{eq:EUT}
\inmat{(EUT)}\quad V_u(\xi):=\bbe_{\mathbb{P}}[u(\xi(\omega))]=\int_{-\infty}^{\infty}u(x)dF_{\xi}(x),
\edeqn
where $u:[\underline{x},\bar{x}] \to \R$ is a nondecreasing utility function,
$F_{\xi}$ is the cumulative distribution function of random variable $\xi:\Omega \to \R$.
If we confine our discussion to the case where a DM's preference can be described by a utility function with
$u(\underline{x})=0$, $u(\bar{x})=1$,
then we make the following assumption.
\begin{assumption}
[EUT Decision Maker]
\label{ass:EUT}
Let ``$\succeq$" be a binary relation which signifies the DM's preference over ${\cal L}^0
$ and $V_{u}$ is given by \eqref{eq:EUT}.
There exists a non-decreasing utility function
$u^*:[\underline{x},\bar{x}]\to \R$ 
such that for any $\xi_1,\xi_2\in {\cal L}^0 $,
$$
\xi_1\succeq \xi_2 \;
\Longleftrightarrow
V_{u^*}(\xi_1)\geq V_{u^*}(\xi_2).
$$
\end{assumption}
We refer to 
$V_{u^*}:{\cal L}^0\to \R$ the preference functional 
under EUT with the true utility $u^*$. We focus
on the case where $u^*$ exists but information about $u^*$ is incomplete from the modeler's perspective.

\vspace{-0.1cm}
\subsection{Piecewise Linear Approximation of the Utility Function}
\label{sec:PLA-U}
\vspace{-0.1cm}

In this paper, we propose to use the polyhedral method to elicit the true utility function.
Since the existing polyhedral method relies on the linear representation of utility, it cannot be applied directly to nonlinear utility functions. We therefore approximate the nonlinear utility function by a PLF, which admits a linear representation in terms of its increment vector and thereby facilitates the application of the polyhedral method.
Specifically,
given a set of breakpoints ${\cal X}$ with 
cardinality 
$N=|{\cal X}|$,
let $u_N^*:[\underline{x},\bar{x}]\to \R$ be 
a PLF approximation
of $u^*$
with $u_N^*(x_i)=u^*(x_i)$, 
for $i\in [N]$.
Our thought is to elicit $u^*$ by eliciting the PLF $u_N^*$.
To do this,
we must present sufficient conditions to ensure that $u_N^*$ 
approximates to $u^*$ well,
and we make the following assumption.

\begin{assumption}[Lipschitz Continuity of the Utility Function]
\label{ass:Lip_u_1}
The true utility function $u^*$ is Lipschitz continuous over $[\underline{x},\bar{x}]$  
but the Lipschitz modulus $L$
is not necessarily known.
\end{assumption}
This assumption is needed when we quantify the difference between 
$u_N^*$ and $u^*$.
Specifically,  
under Assumption~\ref{ass:Lip_u_1},
we can easily derive 
$|u^*(x)-u_N^*(x)|\leq L\max_{i\in[N-1]}\{x_{i+1}-x_i\}$.
Consequently,
$|u^*(x)-u_N^*(x)|$
is small as long as $\max_{i\in [N-1]}\{x_{i+1}-x_i\}$ is small.
The assumption is satisfied by many utility functions of practical interest,
e.g.,~the power form 
$u(x) = x^a$,
where $0<a<1$, 
on any domain bounded away from zero,
and the exponential form 
$u(x)= \lambda(1-e^{-\beta x})/\beta$ for $x<0$,
$u(x)= (1-e^{-\alpha x})/\alpha$ for $x\geq 0$,
where $\alpha>0$, $\beta>0$,
$\lambda\geq 1$.
In the literature on PRO,
construction of  
an ambiguity set of utility functions
often requires complete information 
about $L$ \citep{GXZ22},
which is undesirable in this setup since $u^*$ is unknown. 
By contrast, this paper does not require 
knowledge of the modulus $L$, thereby broadening its applicability.

\begin{definition}
[Increment-Based PLFs, \cite{HZXZ22}]
\label{def:g}
Let  ${\cal X}:=\{x_1,\cdots,x_N\}$ with $\underline{x}=x_1<\cdots<x_N=\bar{x}$ being fixed breakpoints.
Define the
PLF $u_N:[\underline{x},\bar{x}]\to [0,1]$,
\bgeqn
\label{eq:PLA_utility}
u_N(x):={\bm v}_{[N-1]}^{\top} {\bm g}(x),
\edeqn
where ${\bm v}_{[N-1]}:=(v_1,\cdots,v_{N-1})\in \R^{N-1}$ is 
the vector of utility increments,
and $v_i$ represents the utility increment
over the interval $[x_i,x_{i+1}]$, for $i\in [N-1]$,
${\bm g}(x_1)=0$ and
\bgeqn
\label{eq:g(x)-incre} 
{\bm g}(x)=\Big(\underbrace{1,\cdots,1}_{i-1},\frac{x-x_i}{x_{i+1}-x_i},\underbrace{0,\cdots,0}_{N-1-i}\Big)
\in \R^{N-1}, \inmat{ for }x\in (x_{i},x_{i+1}],
~~i\in [N-1].
\edeqn
\end{definition}

\begin{figure}[!htbp]
\begin{minipage}[t]{\linewidth}
\centering
\includegraphics[width=2.5in]{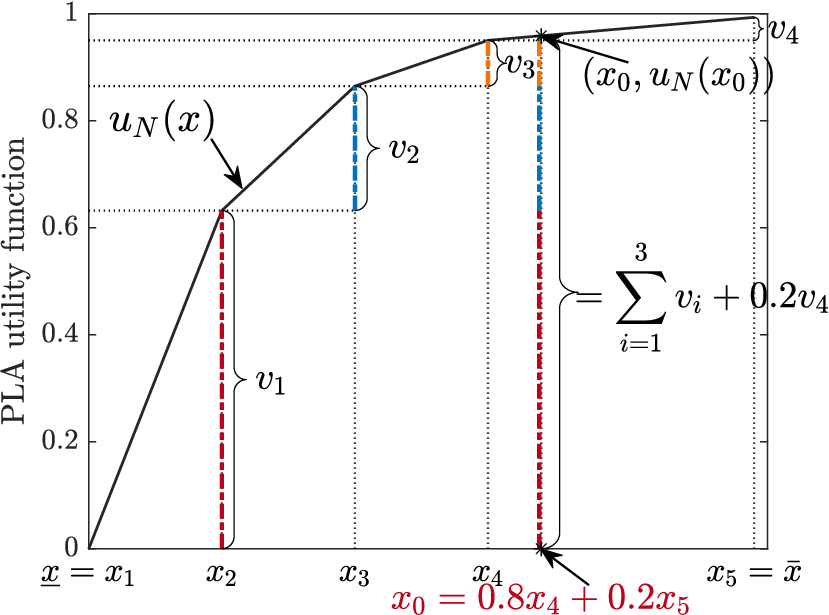}
\end{minipage}
\caption{PLF approximation of the normalized utility function
$u^*(x)=1-e^{-10x}$ over $[-0.5,0]$,
${\cal X}=\{x_1,\cdots,x_5\}$,
$x_0=0.8*x_4+0.2*x_5$.
${\bm g}(x_0)
=(1,1,1,0.2)$,
$u_N(x_0)=(v_1,v_2,v_3,v_4)^\top (1,1,1,0.2)=\sum_{i=1}^3 v_i+0.2v_4$.}
\label{fig:PLFs}
\vspace{-0.2cm}
\end{figure}

We confine our study to the case that $u_N$ is normalized over $[\underline{x},\bar{x}]$
and thus ${\bm e}^\top {\bm v}_{[N-1]}=1$.
Then $v_{N-1}=1-{\bm e}^\top {\bm v}_{[N-2]}$,
where ${\bm v}_{[N-2]}=(v_1,\cdots,v_{N-2})$.
In this set-up,
each $u_N(x)$
is uniquely determined by 
${\bm v}_{[N-2]}$ 
given a specified set of breakpoints 
${\cal X}$.
In Figure~\ref{fig:PLFs},
we depict 
a PLF $u_N$ given in \eqref{eq:PLA_utility}.

\begin{remark}
[Piecewise Linear  Utility $u_N(x)$ in \eqref{eq:PLA_utility} 
versus Linear Utility $U_{\bm x}$ in \eqref{eq:Ux-market}]
The PLF $u_N(x)$
defined in \eqref{eq:PLA_utility} 
resembles the deterministic component
of $U_{\bm x}$ in \eqref{eq:Ux-market} mathematically, where
${\bm v}_{[N-1]}$ corresponds to  ${\bm \beta}$
and ${\bm g}(x)$ corresponds to ${\bm x}$.
Like $U_{\bm x}$ which is uniquely determined by 
${\bm \beta}$, here $u_N(x)$ is uniquely determined by ${\bm v}_{[N-1]}$. 
Both parameter vectors take continuous values, although ${\bm v}_{[N-1]}$
is additionally constrained to be nonnegative and to sum to one.
This structural
similarity enables us to 
extend the polyhedral method in \eqref{eq:poly-cut-market} to
reduce progressively the feasible set of
the utility increment vector ${\bm v}_{[N-1]}$ in $u_N(x)$.
The main difference is
 ${\bm g}(x)$. 
In $U_{\bm x}$,
the component of a profile ${\bm x}\in \{0,1\}^m$ 
represents 
the inclusion of the corresponding attribute in
($x_i=1$) or not ($x_i=0$).
By contrast, for
$x\in (x_i,x_{i+1}]$,
the $i$-th component of ${\bm g}(x)$ varies continuously over $(0,1]$,
while its preceding components equal one and its subsequent components equal zero.
\end{remark} 

\vspace{-0.5cm}
\section{Preparation for the Modified Polyhedral Method}
\label{sec:models}
\vspace{-0.15cm}

By the increment-based representation introduced in Section~\ref{sec:PLA-U}, eliciting the PLF
$u_N^*$ is 
equivalent to eliciting its utility increment vector
${\bm v}^*_{[N-2]}\in {\cal S}_{N-2}$ 
representing the first $N-2$ coordinates of ${\bm v}_{[N-1]}^*$,
where ${\cal S}_{N-2}$ denotes the $(N-2)$-simplex.
The feasible set of increment vectors is characterized by two types of constraints:
those implied by the initial assumptions on ${\bm v}$ and those generated from 
lottery comparisons.
A vector of increments satisfying all these constraints is said to be admissible.
We next introduce
the initial polyhedron of vectors ${\bm v}\in \R^{N-2}$ and 
characterize how it is progressively updated based on the DM's responses to pairwise lottery,
which will be used in the MPM in Section~\ref{sec:strategies}.

\vspace{-0.1cm}
\subsection{Initial Polyhedron}
\vspace{-0.1cm}

Before designing queries for the DM, we express
the initial feasible polyhedron
of ${\bm v}\in \R^{N-2}$  in standard form,
which facilitates the computation of its analytic center in the subsequent analysis:
\bgeqn
\label{eq:initial_P}
{\cal V}_N^0:=\left\{{\bm v}\in \R^{N-2}: {\bm A} {\bm v} \leq {\bm b}\right\},
\edeqn
where
${\bm A}:=({\bm a}_1,\cdots,{\bm a}_{N-1})^{\top}\in \R^{(N-1)\times (N-2)}$,
${\bm b}:=(b_1,\cdots,b_{N-1})\in \R^{N-1}$, 
$ {\bm a}_1:={\bm e}$, $b_{1}:=1$,
${\bm a}_{j}:= -{\bm 1}_{j-1}$,
$b_{j}:=0$, 
for $j\in [N-1]\backslash{\{1\}}$,
where ${\bm e}:=\overbrace{(1,\cdots,1)}^{N-2}$
and ${\bm 1}_j:=(0,\cdots,0,\overbrace{1}^{j},0,\cdots,
0)\\
{\in \R^{N-2}}$.
Set ${\cal V}_N^0$ is a simplex in $\R^{N-2}$.
In  \eqref{eq:initial_P},
constraint ${\bm a}_1^{\top} {\bm v}\leq b_1$ represents $\sum_{j=1}^{N-2}v_j\leq 1$, which
ensures the nonnegativity of the increment $v_{N-1}=1-\sum_{j=1}^{N-2}{v}_{j}$.
This yields
$u_N(\bar{x})=1$.
The remaining constraints
${\bm a}_{j}^{\top} {\bm v}\leq b_{j}$
enforce the nonnegativity 
of increments $v_j$ for $j\in [N-1]\backslash{\{1\}}$ and hence
the nondecreasing property
of the utility function.
The initial polyhedron 
${\cal V}_N^0$ is constructed based on the available information about the DM's preferences prior to asking the DM 
to answer pairwise comparison queries. Each ${\bm v}$ 
corresponds to a feasible PLF consistent with the initial preference information. 
Our objective is to progressively reduce this feasible set ${\cal V}_{N}^0$ through pairwise comparison queries.
The main issue is how to design the queries such that the polyhedra
are cut efficiently.

\vspace{-0.1cm}
\subsection{Structures of Lotteries for Pairwise Comparison}
\vspace{-0.1cm}

In eliciting the
DM's preferences,
each query consists of
a risky lottery and a certain lottery.

\begin{definition}
[Pairwise Comparison of Lotteries under EUT]
Let $\{x,1-p;y\}$ 
denote the lottery which yields outcome  $x$ with probability  $1-p$ and outcome $y$ with probability $p$ (without loss of generality, we assume $y>x$).
The $m$-th query 
takes the form of a comparison between 
lotteries
\bgeqn
\label{eq:AmBm}
{\bm A}_m=\{r_1^m,1-p^m;r_3^m\}
\quad 
\inmat{and} \quad
{\bm B}_m=r_2^m.
\edeqn
The DM's response takes one of the following two forms:
(a) 
${\bm B}_m$ is preferred to ${\bm A}_m$, 
written ${\bm A}_m \preceq {\bm B}_m$;
(b) 
${\bm A}_m$ is preferred to ${\bm B}_m$, written 
${\bm A}_m \succeq {\bm B}_m$.
 In the case of indifference, either alternative may be recorded as the
DM's response.
Under Assumption~\ref{ass:EUT},
the preference information 
can be translated into expected utility values: 
in case (a)
$\bbe[u^*({\bm A}_m)]\leq \bbe[u^*({\bm B}_m)]$; 
in case (b), $\bbe[u^*({\bm A}_m)] \geq \bbe[u^*({\bm B}_m)]$,
where $u^*$ is the true utility function which captures 
DM's preferences.
\end{definition}
\vspace{-0.3cm}

Let 
\bgeqn
\label{eq:def_G}
G_{\bm A}(r_1^m,r_3^m,p^m):=(1-p^m){\bm g}(r_1^m)+p^m{\bm g}(r_3^m)
\quad 
\inmat{and}
\quad
G_{\bm B}(r_2^m):={\bm g}(r_2^m),
\edeqn
where ${\bm g}$ is defined as in \eqref{eq:g(x)-incre},
and ${\bm R}$ denote the affine recovery operator defined by
\bgeqn
\label{eq:def_R}
{\bm R}{\bm v}:=({\bm v},1-{\bm e}^{\top}{\bm v})={\bm v}_{[N-1]}.
\edeqn
By equation
\eqref{eq:PLA_utility}, $u_N(r_i^m)=({\bm R}{\bm v})^{\top}{\bm g}(r_i^m)$. 
The expected utilities of
${\bm A}_m$ and ${\bm B}_m$ are, respectively,
\bgeqn
\label{eq:expectedutility}
\bbe[u_N({\bm A}_m)]=({\bm R}{\bm v})^{\top}G_{\bm A}(r_1^m,r_3^m,p^m)
\quad
\inmat{and}
\quad
\bbe[u_N({\bm B}_m)]=({\bm R}{\bm v})^{\top}G_{\bm B}(r_2^m).
\edeqn
It follows from \eqref{eq:expectedutility} that the DM's response induces a linear constraint on ${\bm v}$:
(a) If ${\bm B}_m$ is preferred, 
then $({\bm R}{\bm v})^{\top}G_{\bm A}(r_1^m,r_3^m,p^m)\leq ({\bm R}{\bm v})^{\top}G_{\bm B}(r_2^m)$;
(b) If ${\bm A}_m$ is preferred, 
then $({\bm R}{\bm v})^{\top}G_{\bm A}(r_1^m,\\
r_3^m,p^m)\geq ({\bm R}{\bm v})^{\top}G_{\bm B}(r_2^m)$.
Let ${\cal V}_{N}^{m-1}$ denote the polyhedron of feasible increment vectors after observing the DM's  responses to the 
first $(m-1)$-queries.

\begin{definition}[Iterative Construction of the Polyhedron]
After collecting the DM's answers to the $m$-th query,
${\cal V}_{N}^{m-1}$ is updated with
\bgeqn
\label{eq:choice-xkyk-0}
{\cal V}_{N}^{m}:={\cal V}_{N}^{m-1} \bigcap 
\left\{{\bm v}\in \R^{N-2}:
h_m   \cdot({\bm R}{\bm v})^{\top}(G_{\bm A}(r_1^m,r_3^m,p^m)-G_{\bm B}(r_2^m))\leq 0\right\},
\edeqn
where $h_m=1$ if the DM chooses ${\bm B}_m$, and $h_m=-1$
if the DM chooses ${\bm A}_m$.
\end{definition}
Note that each ${\bm v}$ uniquely determines a PLF $u_N$.
Thus ${\cal V}_N^m$ induces a feasible set of PLFs consistent with the
observed preference information.

\vspace{-0.1cm}
\subsection{Analytic Center and Sonnevend's Ellipsoid}
\label{sec:AC}
\vspace{-0.1cm}

The modified polyhedral algorithm in this paper relies 
on the analytic center and  Sonnevend's ellipsoid of  ${\cal V}_{N}^{m-1}$.
For completeness, we first recall these
two notions for a general polyhedron.

\begin{definition}
[Analytic Center]
\label{def:AC}
The analytic center of the polyhedron
$$
\mathcal{P}:=\left\{{\bm v} \in  \R^{N-2}:\tilde{\bm a}_i^\top {\bm v} \leq \tilde{b}_i, i\in [I]\right\}
$$
is defined as the optimal 
solution of the minimization problem: 
$\min_{{\bm v}\in \R^{N-2}} \Phi({\bm v})$,
where the logarithmic barrier function is defined as
$\Phi({\bm v}):=-\sum_{i=1}^I \log(\tilde{b}_i-\tilde{\bm a}_i^\top {\bm v})$.
\end{definition}
If we assume that 
$\|\tilde{\bm a}_i\| = 1$, then the slack $\tilde{b}_i-\tilde{\bm a}_i^\top {\bm v}$ 
equals the Euclidean distance from
${\bm v}$ to the hyperplane ${\cal H}_i = \{{\bm v}: \tilde{\bm a}_i^\top {\bm v} = \tilde{b}_i\}$.
Under this normalization, 
the analytic center is the point that maximizes the product of
distances to the defining hyperplanes $\tilde{\bm a}_i^\top {\bm v} = \tilde{b}_i$ (see e.g., \cite{BoV04}).
Given the analytic center, we next introduce the Sonnevend's ellipsoid 
associated with $\mathcal{P}$.

\begin{definition}
[Sonnevend's Inner Ellipsoid, 
\cite{Son85}]
\label{def:inner-ellipsoid}
Given the analytic center ${\bm c}$ of polyhedron $\mathcal{P}$, the Sonnevend's inner ellipsoid of 
$\mathcal{P}$ is defined as
$$
E_{\mbox{inner}}:=\left\{{\bm v}\in \R^{N-2}:
({\bm v}-{\bm c})^{\top}{\bm H}({\bm v}-{\bm c})\leq 1
\right\},
$$
where ${\bm H}:=\nabla^2\Phi({\bm c})
=\sum_{i=1}^I\frac{1}{(\tilde{b}_i-\tilde{\bm a}_i^\top {\bm c})^2}\tilde{\bm a}_i^{\top}\tilde{\bm a}_i$
is the Hessian of the logarithmic barrier function evaluated at the analytic center.
\end{definition}
It follows from \cite{Son85} that 
$E_{\mbox{inner}} \subset \mathcal{P}$,
see also \cite[Section 8.5.3]{BoV04}.
By Definition~\ref{def:AC},
the analytic center of polyhedron ${\cal V}_{N}^{m-1}$,  denoted by ${\bm c}^{m}$, 
is the optimal solution in ${\bm v}$ to the following problem:
\vspace{-0.3cm}
\begin{subequations}
\label{eq:ana_center}
\begin{align}
\min\limits_{{\bm v},{\bm s}}~~& -\sum_{i=1}^{N-1+m-1} \log(s_i) \\
{\rm s.t.}\;\;\;\; & {\bm A}^{\top} {\bm v}+{\bm s}={\bm b}, {\bm s}\geq 0,\label{eq:one-dimensional-AC-b}\\
& ({\bm R}{\bm v})^{\top}(G_{\bm A}(r_1^k,r_3^k,p^k)-G_{\bm B}(r_2^k))+s_{N-1+k}=0  \; \inmat{if } {\bm A}_k \preceq {\bm B}_k,\; k\in[m-1], \label{eq:one-dimensional-AC-c}\\
& ({\bm R}{\bm v})^{\top}(G_{\bm A}(r_1^k,r_3^k,p^k)-G_{\bm B}(r_2^k))-s_{N-1+k}
=0\;  \inmat{if } {\bm A}_k  \succeq {\bm B}_k,\;k\in[m-1],  \label{eq:one-dimensional-AC-d}
\end{align}
\end{subequations}
where ${\bm s}:=(s_1,\cdots,s_{N-1+m-1})\in \R^{N-1+m-1}$,
${\bm A}$ and ${\bm b}$ are defined as in \eqref{eq:initial_P} and correspond to
initial information of ${\bm v}$.
Constraints \eqref{eq:one-dimensional-AC-c}-\eqref{eq:one-dimensional-AC-d} arise from the linear constraints of ${\bm v}$ based on the DM's answers to the previous $m-1$ queries.
Since the analytic center lies in the interior of ${\cal V}_N^{m-1}$,
the optimal solution $({\bm c}^m,{\bm s}^m)$ of problem \eqref{eq:ana_center}  satisfies that 
${\bm s}^m>0$, ${\bm A}{\bm c}^m<{\bm b}$.
Hence,
\bgeqn 
\label{eq:cm}
0<{c}^m_i<1, \; \forall i\in [N-2]
\quad \inmat{and}\quad   
0<c_{N-1}^m=1-{\bm e}^\top {\bm c}^m<1.
\edeqn

In the next section, 
we show how Sonnevend's inner ellipsoid can be used to design the
lotteries ${\bm A}_m$ and ${\bm B}_m$ in \eqref{eq:AmBm} 
and subsequently construct a cutting hyperplane for
${\cal V}_{N}^{m-1}$.

\vspace{-0.2cm}
\section{Design of Lotteries and Update of Polyhedra}
\label{sec:MPM}
\vspace{-0.2cm}

In this section, we move on to the details
as to how to implement iterative updates of the 
polyhedra of increments specified in \eqref{eq:choice-xkyk-0} so that the polyhedra are reduced efficiently.
This in turn requires us to design lotteries 
${\bm A}_m$ and ${\bm B}_m$, 
which amounts to optimally determining the parameters $r_1^m$, $r_2^m$, $r_3^m$, and $p^m$.
The process can be summarized as follows. 
First,
we obtain the analytic center of the polyhedron ${\cal V}_N^{m-1}$.
Second, 
as in the classical polyhedral method,
we obtain
the two intersections between the longest axis of the Sonnevend's ellipsoid $E_{\mbox{inner}}$ 
and the boundary
of the polyhedron,
denoted by ${\bm v}_1^m$ and ${\bm v}_2^m$ (see \cite{THS04}).
Finally, 
we identify $r_i^m$, for $i=1,2,3$ and $p^m$ by solving two optimization problems to be detailed in Section~\ref{sec:two_problem}.

\vspace{-0.2cm}
\subsection{ Design of Lotteries}
\label{sec:two_problem}
\vspace{-0.1cm}

The original polyhedral method
uses problems \eqref{eq:original_PM_a}-\eqref{eq:original_PM_b} to generate the $m$-th 
query by maximizing the linear utility functions over appropriate half-spaces. 
We follow the same idea
by replacing
the linear utility functions  in \eqref{eq:original_PM_a}-\eqref{eq:original_PM_b} with expected piecewise linear utility functions given in
\eqref{eq:expectedutility}.

\vspace{-0.3cm}
\begin{definition}[Construction of 
($r_1^m$, $r_2^m$, $r_3^m$, $p^m$)]
Let ${\bm c}^m$ be the analytic center obtained by solving \eqref{eq:ana_center} 
and let ${\bm v}_1^m$, ${\bm v}_2^m$ be the intersections between the longest axis of  the Sonnevend's ellipsoid $E_{\mbox{inner}}$
and the boundary of the polyhedron
${\cal V}_{N}^{m-1}$. 
Define
$F:=\{(r_1,r_3,p)\in [\underline{x},\bar{x}]\times [\underline{x},\bar{x}]\times [0,1]:r_1\leq r_3\}$.
Solve
\begin{subequations}
\label{eq:gene_m_Q}
\begin{align}
& r_2^m\in \arg \max_{r_2}\left\{({\bm R}{\bm v}_1^m)^{\top}G_{\bm B}(r_2):
r_2\in [\underline{x},\bar{x}],\;({\bm R}{\bm c}^m)^{\top}G_{\bm B}(r_2) \leq D
\right\},\label{eq:gene_m_Q-1-v1-p}\\
&(r_1^m,r_3^m,p^m)\in  \arg\max_{r_1,r_3,p} \left\{({\bm R}{\bm v}_2^m)^{\top}G_{\bm A}(r_1,r_3,p):
\begin{array}{l}
(r_1,r_3,p)\in F\\
({\bm R}{\bm c}^m)^{\top}G_{\bm A}(r_1,r_3,p) \leq D
\end{array}
\right\},\label{eq:gene_m_Q-2-v1-p}
\end{align}
\end{subequations}
where $G_{\bm B}(\cdot)$ and $G_{\bm A}(\cdot,\cdot,\cdot)$ are defined as in \eqref{eq:def_G} and 
$D\in (0,1]$ is a parameter.
\end{definition}
\vspace{-0.3cm}

\begin{figure}[!ht]
\vspace{-0.5cm}
\begin{minipage}[t]{0.33\linewidth}
\centering
\includegraphics[width=2.1in]{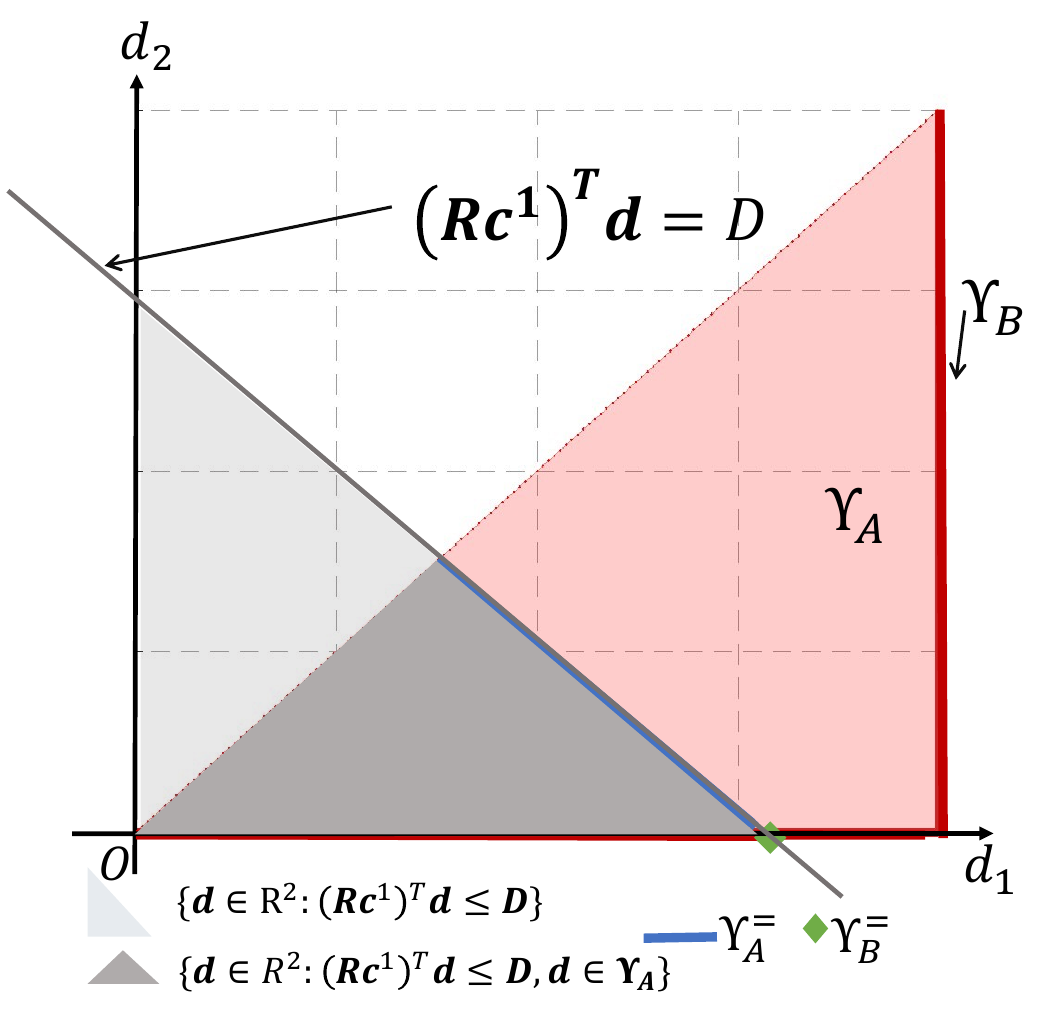}
\text{\footnotesize{\textbf{(a)}}}
\end{minipage}
\hspace{-0.22cm}
 \begin{minipage}[t]{0.33\linewidth}
\centering
\includegraphics[width=2.19in]
{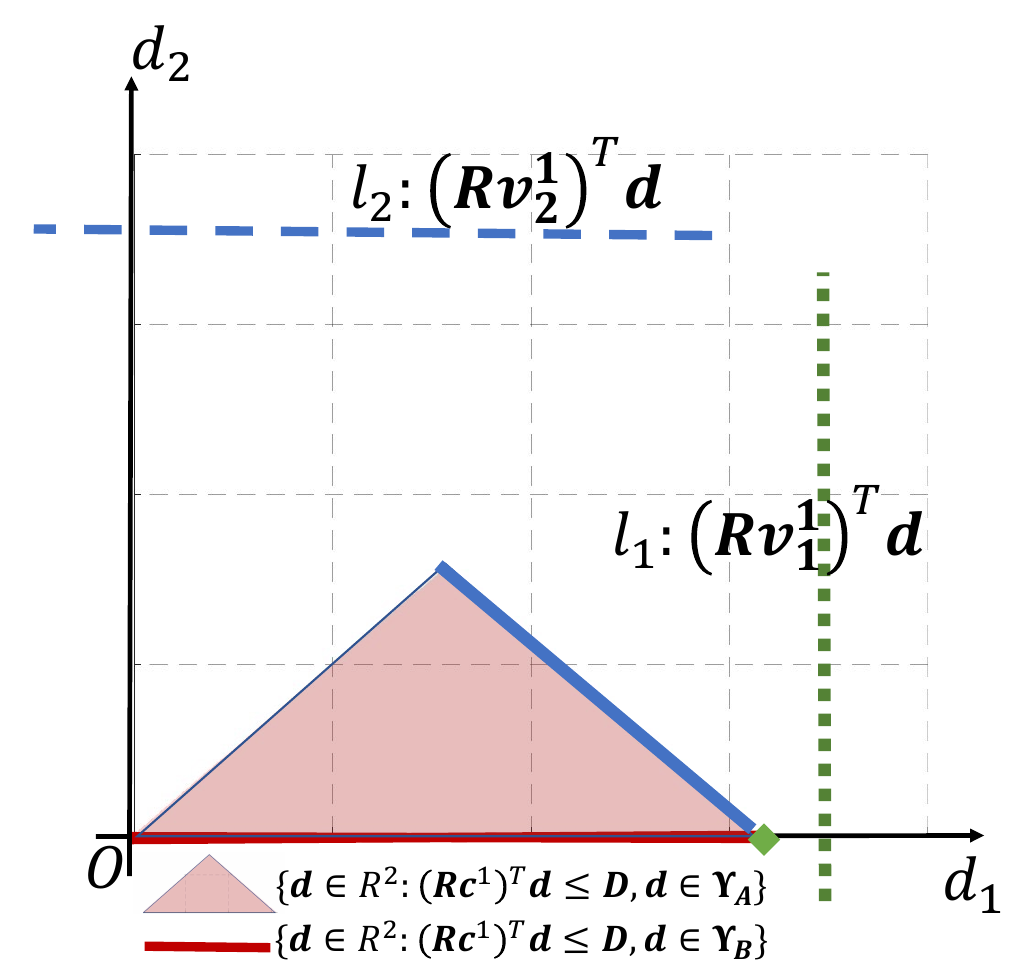}
\text{\footnotesize{\textbf{(b)}}}
\end{minipage}
\begin{minipage}[t]{0.33\linewidth}
\centering
\includegraphics[width=2.2in]{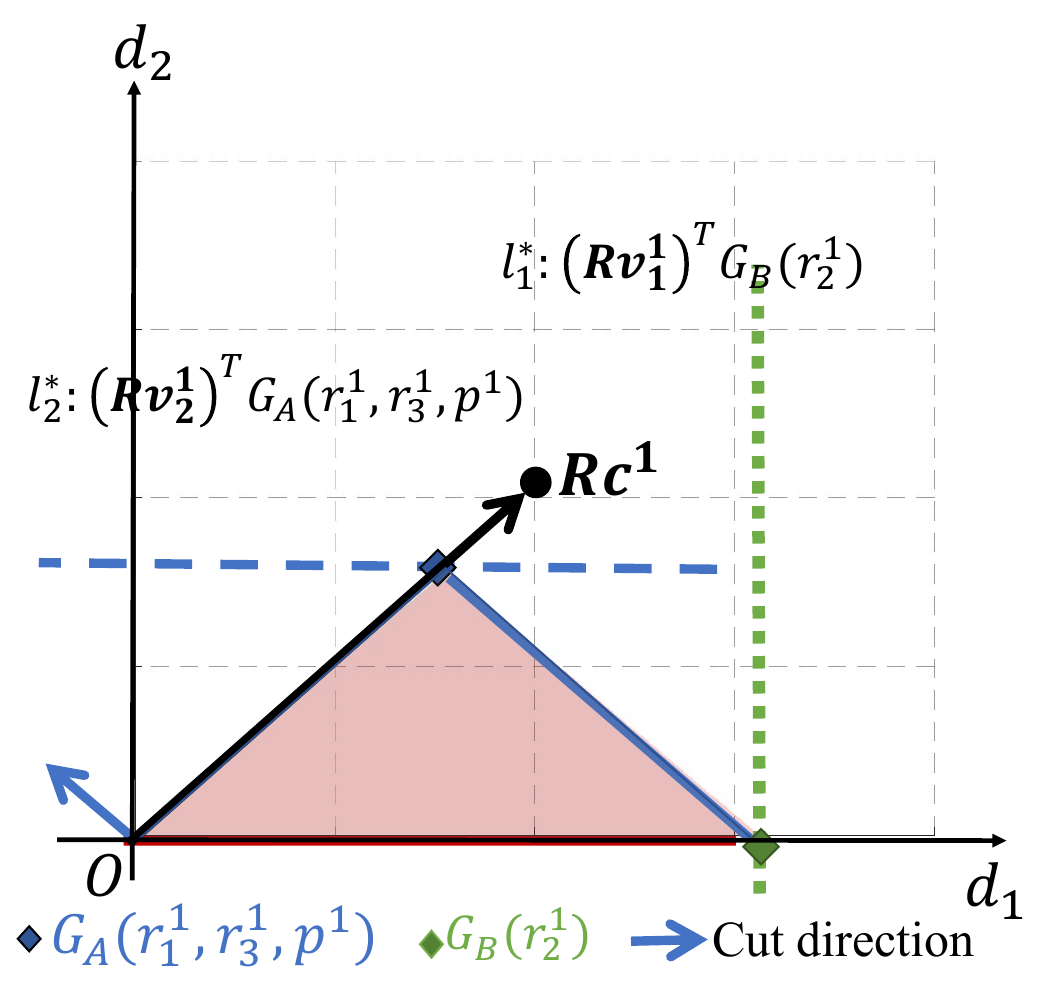}
\text{\footnotesize{\textbf{(c)}}}
\end{minipage}
\caption{\footnotesize{
(a) Illustration of the feasible sets
$\Upsilon_{\bm A}$, $\Upsilon_{\bm B}$, 
$\Upsilon_{\bm A}^{=}$
and $\Upsilon_{\bm B}^{=}$.
The gray area represents $\{{\bm d}\in \R^2:({\bm R}{\bm c}^1)^\top {\bm d}\leq D\}$;
the red area represents $\Upsilon_{\bm A}$ defined as in \eqref{eq:Upsilon_A};
and the red segments represents $\Upsilon_{\bm B}$ defined as in \eqref{eq:Upsilon_B}.
The blue segment represents $\Upsilon_{\bm A}^=$ and the green diamond 
represents $\Upsilon_{\bm B}^=$, respectively, 
as given in \eqref{eq:UpsilonAB}.
(b) Illustration of feasible sets and objective functions of problems \eqref{eq:G_A_G_B-rfmlt}.
The triangular region 
is the intersection of $\{{\bm d}\in \R^2:({\bm R}{\bm c}^1)^\top {\bm d}\leq D\}$ and $\Upsilon_{\bm A}$, 
and hence is the feasible set $\{{\bm d}\in \R^2: ({\bm R}{\bm c}^1)^\top {\bm d}\leq D,{\bm d}\in \Upsilon_{\bm A}\}$ of problem in \eqref{eq:G_A_G_B-rfmlt-b}.
The red segment on $d_1$-axis is the intersection of $\{{\bm d}\in \R^2:({\bm R}{\bm c}^1)^\top {\bm d}\leq D\}$ and $\Upsilon_{\bm B}$,
and hence is the feasible set $\{{\bm d}\in  \R^2: ({\bm R}{\bm c}^1)^\top {\bm d}\leq D, {\bm d}\in \Upsilon_{\bm B}\}$ 
 of problem in \eqref{eq:G_A_G_B-rfmlt-a}.
The dashed blue and green lines are the contours of the functions 
$l_2({\bm d}) := ({\bm R}{\bm v}_2^1)^\top {\bm d}$
and $l_1({\bm d}) := ({\bm R}{\bm v}_1^1)^\top {\bm d}$ over $\R^2$,
respectively.
(c) Illustration of the optimal solutions and the resulting
cut direction $G_{\bm A}(r_1^1,r_3^1,p^1)-G_{\bm B}(r_2^1)$ for the fist
cut hyperplane.
Problem \eqref{eq:G_A_G_B-rfmlt-b}
maximizes
$l_2({\bm d})$
over the feasible set $\{{\bm d}\in \R^2: ({\bm R}{\bm c}^1)^\top {\bm d}\leq D, {\bm d}\in \Upsilon_{\bm A}\}$,
with the maximum being attained  
at the blue diamond point ${\bm d}_{\bm A}^*=G_{\bm A}(r_1^1,r_3^1,p^1)$.
Similarly, problem \eqref{eq:G_A_G_B-rfmlt-a} 
maximizes $l_1({\bm d})$ over $\{{\bm d}\in \R^2:({\bm R}{\bm c}^1)^\top {\bm d}\leq D, {\bm d}\in \Upsilon_{\bm B}\}$,
with the maximum being attained at the green diamond point ${\bm d}_{\bm B}^*=G_{\bm B}(r_2^1)$.
}}
\label{fig:polyhedron_2-2}
\vspace{-0.5cm}
\end{figure}

We first provide details on the reformulations of problems \eqref{eq:gene_m_Q-1-v1-p}--\eqref{eq:gene_m_Q-2-v1-p} by introducing \(\bm d\) to represent the values of \(G_{\bm A}(r_1,r_3,p)\) and \(G_{\bm B}(r_2)\), respectively. We then use a specific example to visualize the resulting feasible sets and objective functions.
Let
\bgeqn 
\label{eq:Upsilon_A}
\Upsilon_{\bm A} := \{G_{\bm A}(r_1,r_3,p): (r_1,r_3,p)\in F
\}
= \{(1-p)
{\bm g}(r_1) + p{\bm g}(r_3):
(r_1,r_3,p)\in F
\}
\edeqn 
denote the image 
of $G_{\bm A}$ over $F$
and 
\bgeqn 
\label{eq:Upsilon_B}
\Upsilon_{\bm B}:= \{G_{\bm B}(r_2): r_2\in [\underline{x},\bar{x}]\} = 
\{{\bm g}(r_2): r_2\in [\underline{x},\bar{x}]\} 
\edeqn 
denote the image 
of $G_{\bm B}$ over $[\underline{x},\bar{x}]$,
see Figures~\ref{fig:polyhedron_2-2}.
We consider problems
\begin{subequations}
\begin{align}
&\max\limits_{{\bm d}\in \Upsilon_{\bm B}} ~\{({\bm R}{\bm v}_1^m)^{\top}{\bm d}: 
({\bm R}{\bm c}^m)^\top {\bm d}\leq D\}\label{eq:G_A_G_B-rfmlt-a},\\
&\max\limits_{{\bm d}\in \Upsilon_{\bm A}}~\{({\bm R}{\bm v}_2^m)^{\top}{\bm d}:
({\bm R}{\bm c}^m)^\top {\bm d}\leq D\}\label{eq:G_A_G_B-rfmlt-b}.
\end{align}
 \label{eq:G_A_G_B-rfmlt}
\end{subequations} 
A clear benefit of reformulating
\eqref{eq:gene_m_Q}
as \eqref{eq:G_A_G_B-rfmlt} is 
the optimization is carried out directly over the images $\Upsilon_{\bm A}$
and $\Upsilon_{\bm B}$, with both the objective functions and the additional constraints being linear in ${\bm d}$.
In the subsequent analysis, we
use 
both  \eqref{eq:G_A_G_B-rfmlt} and \eqref{eq:gene_m_Q} interchangeably, 
depending on which formulation is more convenient.
Define $G_{\bm B}^{-1}:\Upsilon_{\bm B}\to [\underline{x},\bar{x}]$ and 
$G_{\bm A}^{-1}:\Upsilon_{\bm A}\to F$ by
\bgeqn
\label{eq:GAB-1}
G_{\bm B}^{-1}({\bm d}):=\{r_2\in [\underline{x},\bar{x}]: G_{\bm B}(r_2)={\bm d}\}
\quad \inmat{and}\quad
G_{\bm A}^{-1}({\bm d}):=\{(r_1,r_3,p)\in F: G_{\bm A}(r_1,r_3,p)={\bm d}\}. \quad 
\edeqn
Since $G_{\bm B}$ is one-to-one over $[\underline{x},\bar{x}]$,
its inverse is  
single-valued over $\Upsilon_{\bm B}$,
whereas $G_{\bm A}^{-1}(\cdot)$ is 
generally
set-valued  over $\Upsilon_{\bm A}$.
Thus, in principle, we may
solve the programs in \eqref{eq:G_A_G_B-rfmlt-a}-\eqref{eq:G_A_G_B-rfmlt-b}
and then 
recover
the corresponding optimal solutions 
for $r_2^m$ and $(r_1^m,r_3^m,p^m)$ 
of problems \eqref{eq:gene_m_Q-1-v1-p} and 
\eqref{eq:gene_m_Q-2-v1-p}
through
$G_{\bm A}^{-1}$ and $G_{\bm B}^{-1}$,
respectively. 
However, we do not use this implementation in our numerical experiments.

Next, we use a specific example to illustrate the feasible sets and objective functions of problems \eqref{eq:G_A_G_B-rfmlt-a} and \eqref{eq:G_A_G_B-rfmlt-b}.
Specifically, for problem \eqref{eq:G_A_G_B-rfmlt-a}, 
the feasible set and objective function are 
$\{{\bm d}\in \Upsilon_{\bm B}:({\bm R}{\bm c}^1)^\top {\bm d} \leq D\}$  
and
$({\bm R}{\bm v}_1^1)^\top {\bm d}$ of problem \eqref{eq:G_A_G_B-rfmlt-a},
respectively,
whereas for problem \eqref{eq:G_A_G_B-rfmlt-b}, 
they are
$\{{\bm d}\in \Upsilon_{\bm A}:({\bm R}{\bm c}^1)^\top {\bm d} \leq D\}$ and 
$({\bm R}{\bm v}_2^1)^\top {\bm d}$.
We also illustrate the results 
of Proposition~\ref{prop:EDn=0_2},
by showing that problems \eqref{eq:G_A_G_B-rfmlt-a}-\eqref{eq:G_A_G_B-rfmlt-b}
admit optimal solutions in
$\Upsilon_{\bm B}^=$ and $\Upsilon_{\bm A}^=$, 
respectively,
in the following
two cases: (a) $N=3$, ${\bm g}(\cdot)\in \R^2$, ${\bm v}\in \R$, and ${\bm R}{\bm v}\in \R^2$, 
as shown in
Figure~\ref{fig:polyhedron_2-2};
and 
(b)$N=4$, ${\bm g}(\cdot)\in \R^3$, ${\bm v}\in \R^2$, and ${\bm R}{\bm v}\in \R^3$.
Consider case (a), where
${\bm R}{\bm c}^1=(0.5,0.5)$,
$({\bm v}_1^1)^\top{\bm d}=(v_{1,1}^1,0)^\top {\bm d}$, $({\bm v}_2^1)^\top {\bm d}=(0,v_{2,2}^1)^\top {\bm d}$,
and
$\Upsilon_{\bm B}=[0,1]\times \{0\}\cup \{1\}\times [0,1]$,
while
$\Upsilon_{\bm A}=\{(d_{A,1},d_{A,2}):d_{A,2}\leq d_{A,1},0\leq d_{A,1}\leq 1,0\leq d_{A,2}\leq 1\}$.
For $D=0.4$,
an 
optimal solution to problem \eqref{eq:G_A_G_B-rfmlt-a} is ${\bm d}_{B}^*=(0.8,0)={\bm g}(r_2^m)$,
where
$j_{r_2^m}=1$. 
In this problem $v_{1,1}^1>0$,
${\bm d}_{B}^*\in \Upsilon_{\bm B}^=$, and $r_2^m\in {\cal F}_{\bm B}^=$.
Similarly, an optimal solution to problem
\eqref{eq:G_A_G_B-rfmlt-b} is ${\bm d}_{A}^*=(0.4,0.4)=(1-0.5){\bm g}(r_1^m)+0.5{\bm g}(r_3^m)$,
where $r_1^m=\underline{x}$, $r_3^m=\bar{x}$, 
$j_{r_1^m}=1$,
and $j_{r_3^m}=2$.
In this case,
$v_{2,2}^1>0$,
${\bm d}_{A}^*\in \Upsilon_{\bm A}^=$, and $(r_1^m,r_3^m,p)\in {\cal F}_{\bm A}^=$.

\vspace{-0.15cm}
\subsection{Construction of Cut Hyperplane} 
\vspace{-0.2cm}

With $r_1^m, r_2^m, r_3^m$ and $p^m$
being identified in the preceding subsection, we are ready to define the $m$-th cut of polyhedron ${\cal V}_N^{m-1}$ formally and discuss its properties:
{\em Choice balance} and {\em Postchoice symmetry}.

\vspace{-0.3cm}
\begin{definition}
[Cut Hyperplane] For given ${\bm A}_m$
and ${\bm B}_m$, the 
cut hyperplane of polyhedron ${\cal V}_{N}^{m-1}$ in the space of  
$\R^{N-2}$ is given by
\bgeqn
\label{eq:cut_hyper}
l_{\rm cut}^m({\bm v}):=
\bbe[u_{N_{m-1}}({\bm A}_m)]-\bbe[u_{N_{m-1}}({\bm B}_m)]
=({\bm R}{\bm v})^{\top}(G_{\bm A}(r_1^m,r_3^m,p^m)-G_{\bm B}(r_2^m))=0.
\edeqn
\end{definition}
\vspace{-0.3cm}

\begin{criterion}[Choice Balance]
The newly generated cutting hyperplane should be as close as possible to the analytic center of the current polyhedron.
\end{criterion}
\vspace{-0.3cm}

The next theorem states
that the cut hyperplane 
always passes through the analytic center ${\bm c}^m$ under some conditions,
that is, this strategy satisfies the choice balance. 

\begin{theorem}
[Choice Balance]
\label{prop:EDn=0-2}
Assume the setting and conditions of Proposition~\ref{prop:EDn=0_2}.
Then the cut hyperplane 
passes through  the analytic center,
i.e.,
\bgeqn
\label{eq:Ch-Balance}
l_{\rm cut}^m({\bm c}^m)=({\bm R}{\bm c}^m)^{\top}(G_{\bm A}(r_1^m,r_3^m,p^m)-G_{\bm B}(r_2^m))=0,
\edeqn 
where $G_{\bm A}(r_1^m,r_3^m,p^m)$ and  $G_{\bm B}(r_2^m)$ are defined as in \eqref{eq:def_G} and ${\bm R}$ is defined as in \eqref{eq:def_R}.
\end{theorem}

\noindent{\bf Proof}.
Since the optimal solutions $r_2^m$ and $(r_1^m,r_3^m,p^m)$ of problems \eqref{eq:gene_m_Q-1-v1-p} and \eqref{eq:gene_m_Q-2-v1-p},
respectively,
satisfy the equalities 
in \eqref{eq:bar_r123-2},
we have
$({\bm R}{\bm c}^m)^{\top}G_{\bm A}(r_1^m,r_3^m,p^m)=D$ and 
$({\bm R}{\bm c}^m)^{\top}G_{\bm B}(r_2^m)=D$.
Therefore,
$$
l_{\rm cut}^m({\bm c}^m)=({\bm R}{\bm c}^m)^{\top}(G_{\bm A}(r_1^m,r_3^m,p^m)-G_{\bm B}(r_2^m))=0.
$$
The proof is complete.
\hfill \Box

Theorem~\ref{prop:EDn=0-2} guarantees choice balance, helping avoid cuts that eliminate only a small portion of the polyhedron and thereby enhancing elicitation efficiency.
By the optimality of
${\bm v}_1^m$ and ${\bm v}_2^m$ in the query construction,
Theorem~\ref{prop:EDn=0-2} 
implies that
$
l_{\rm cut}^m({\bm v}_1^m) \leq 0$, and
$l_{\rm cut}^m({\bm v}_2^m)\geq 0,
$
so ${\bm v}_1^m$ and ${\bm v}_2^m$ lie
on the two sides of the hyperplane.
To see this, we note that
\bgeq
l_{\rm cut}^m({\bm v}_2^m)
&=&({\bm R}{\bm v}_2^m)^\top \Big(G_{\bm A}(r_1^m,r_3^m,p^m)-G_{\bm B}(r_2^m)\Big)\\
&=&({\bm R}{\bm v}_2^m)^\top \Big(G_{\bm A}(r_1^m,r_3^m,p^m)
- G_{\bm A}(r_2^m,r_2^m,p^m)
\Big)
\geq 0,
\edeq
which implies, by \eqref{eq:Ch-Balance}, that    
$l_{\rm cut}^m({\bm v}_2^m)=({\bm R}{\bm v}_2^m-{\bm R}{\bm c}^m)^\top \Big(G_{\bm A}(r_1^m,r_3^m,p^m)-G_{\bm B}(r_2^m)\Big)\geq 0$.
Note that ${\bm R}{\bm v}_2^m-{\bm R}{\bm c}^m$ and ${\bm R}{\bm v}_1^m-{\bm R}{\bm c}^m$
are parallel and point in opposite directions,
i.e., there exists a positive constant $\beta$ such that
$({\bm R}{\bm v}_1^m-{\bm R}{\bm c}^m)=-\beta({\bm R}{\bm v}_2^m-{\bm R}{\bm c}^m)$,
therefore,
\bgeq
l_{\rm cut}^m({\bm v}_1^m) &=&({\bm R}{\bm v}_1^m)^\top \Big(G_{\bm A}(r_1^m,r_3^m,p^m)-G_{\bm B}(r_2^m)\Big)   \\
&\overset{\eqref{eq:Ch-Balance}}{=}& ({\bm R}{\bm v}_1^m-{\bm R}{\bm c}^m)^\top \Big(G_{\bm A}(r_1^m,r_3^m,p^m)-G_{\bm B}(r_2^m)\Big)
\\
&=& -\beta({\bm R}{\bm v}_2^m-{\bm R}{\bm c}^m)^\top \Big(G_{\bm A}(r_1^m,r_3^m,p^m)-G_{\bm B}(r_2^m)\Big)\leq 0.
\edeq
Next, we select a query that satisfies the postchoice symmetry criterion as closely as possible.

\vspace{-0.2cm}
\begin{criterion}[Postchoice Symmetry]
\label{eq:criterion_2}
The newly generated cutting hyperplane should be 
as close as possible to orthogonal
to the longest axis of the current polyhedron.
\end{criterion}
\vspace{-0.2cm}

This criterion is motivated by the observation that the longest axis represents the direction along which the current polyhedron exhibits the largest geometric spread. By generating a cutting hyperplane that is as orthogonal as possible to this axis, the cut targets this most extended direction directly, thereby reducing the geometric imbalance of the resulting polyhedron and avoiding excessively long and narrow feasible regions in subsequent iterations.

To implement this criterion, instead of using a fixed
$D$ in problems \eqref{eq:gene_m_Q-1-v1-p} and \eqref{eq:gene_m_Q-2-v1-p},
we generate a sequence $\{D_s\}_{s=1}^S$ over $(0,1]$,
such as $D_s=\frac{s}{S}$, for $s\in [S]$. 
For each fixed 
$s\in [S]$,
we solve problems \eqref{eq:gene_m_Q-1-v1-p} and \eqref{eq:gene_m_Q-2-v1-p} 
with $D=D_s$
to obtain a sequence $\{r_1^{m,s},r_2^{m,s},r_3^{m,s},p^{m,s}\}_{s=1}^{S}$.
For 
$s\in[S]$, let
$
{\bm d}^{m,s}
:=
G_{\bm A}(r_1^{m,s},r_3^{m,s},p^{m,s})
-
G_{\bm B}(r_2^{m,s})$.
We select 
$s_0\in [S]$ such that 
the normal vector of the resulting cutting hyperplane,
${\bm d}^{m,s}_{[N_{m-1}-2]}-d^{m,s}_{N_{m-1}-1}{\bm e}$,
is
most closely aligned with the longest-axis direction:
\begin{eqnarray}
\label{eq:postchoice-1}
s_0\in\arg\max_{s\in[S]}
\frac{
\left|
({\bm d}^{m,s}_{[N_{m-1}-2]}-d^{m,s}_{N_{m-1}-1}{\bm e})^\top
({\bm v}_1^m-{\bm v}_2^m)
\right|
}{
\|{\bm d}^{m,s}_{[N_{m-1}-2]}-d^{m,s}_{N_{m-1}-1}{\bm e}\|
\,
\|{\bm v}_1^m-{\bm v}_2^m\|}.
\end{eqnarray}
For simplicity,
 we will write $r_i^m$ and $p^m$ for $r_i^{m,s_0}$ and $p^{m,s_0}$ respectively, 
 for $i=1,2,3$.
This procedure requires solving $S$ pairs of problems, and therefore the choice of $S$ 
affects its computational efficiency.

\vspace{-0.3cm}
\section{Modified Polyhedral Algorithm}
\label{sec:strategies}
\vspace{-0.2cm}

Based on the discussions in the preceding two sections, we are ready to 
present detailed algorithmic procedures for the 
MPM described in this section.
We begin with a simple example illustrating
the strategy developed in
\eqref{eq:gene_m_Q-1-v1-p}-\eqref{eq:gene_m_Q-2-v1-p}, and \eqref{eq:postchoice-1}, 
followed by a formal 
description of the modified polyhedral algorithm.

\vspace{-0.2cm}
\subsection{An Illustrative and Simulating Example}
\vspace{-0.2cm}

We begin with a simulating example.

\vspace{-0.2cm}
\begin{example}
\label{exa:concave}
Consider the true utility function $u^*(x)=1-e^{-10x}$, $x\in [-0.5,0.5]$,
which
can be normalized by $\frac{u^*(x)-u^*(-0.5)}{u^*(0.5)-u^*(-0.5)}$.
We use a PLF
$u_N(x)=
({\bm R}{\bm v})^{\top}{\bm g}(x)$ 
with $N=4$ breakpoints 
$\{
-0.5,-0.48596,0,0.5\}$,
${\bm v}=(v_1,v_2)$, ${\bm R}{\bm v}=(v_1,v_2,1-{\bm e}^\top {\bm v})$.
We begin with an initial polyhedron
$
{\cal V}_4^0 :=\left\{{\bm v}\in \R^2:
v_1+v_2\leq 1, v_i\geq 0, i=1,2 \right\}.
$
The analytic center of ${\cal V}_N^0$ is ${\bm c}^1=(1/3,1/3)$ and ${\bm R}{\bm c}^1=(1/3,1/3,1/3)$.
The two intersections between the longest axis 
and the boundary of polyhedron ${\cal V}_4^0$
are
${\bm v}^1=(2/3,0)$,
${\bm v}^2=(0,2/3)$ and
${\bm R}{\bm v}^1=(2/3,0,1/3)$,
${\bm R}{\bm v}^2=(0,2/3,1/3)$.

(i) {Illustration of Postchoice Symmetry}.
Let $D_s=\frac{s}{10}$, for $s\in [10]$. By solving problem \eqref{eq:gene_m_Q-1-v1-p}-\eqref{eq:gene_m_Q-2-v1-p},
we obtain the optimal strategies: 
\bgeqn
\label{eq:r123p}
\{(r_1^{m,s},r_2^{m,s},r_3^{m,s},p^{m,s})\}_{s=1}^{6}=\left\{
\hspace{-0.2cm}
\begin{array}{l}
(-0.5,-0.494,0,0.15),
(-0.5,-0.488,0,0.30),
(-0.5,-0.482,0,0.45)\\
(-0.5,-0.384,0,0.6),
(-0.5,-0.240,0,0.75),
(-0.5,-0.096,0,0.90)
\end{array}
\hspace{-0.2cm}
\right\}~~~
\edeqn
and the resulting 
vector $G_{\bm A}(r_1^m,r_3^m,p^m)-G_{\bm B}(r_2^m)$
in the definition of cut hyperplane  
$
l_{\rm cut}^m({\bm v}):=
\bbe[u({\bm A}_m)]-\bbe[u({\bm B}_m)]
=({\bm R}{\bm v})^{\top}(G_{\bm A}(r_1^m,r_3^m,p^m)-G_{\bm B}(r_2^m))
$
is 
\bgeq
\{l_{\rm cut}^{1,s}\}_{s=1}^{6}=
\left\{
\begin{array}{l}
(-0.15,0.15,0),
(-0.3,0.3,0),
(-0.45,0.45,0)\\
(-0.39,0.39,0),
(-0.24,0.24,0),
(-0.09,0.09,0)
\end{array}\right\}.
\edeq
Note that when $D_s\in \{\frac{7}{10},\cdots,1\}$,
the resulting cut coefficient vector is zero, and the corresponding cut is therefore discarded.
We can see that in the first step, the cut hyperplanes corresponding to $\{(r_1^{m,s},r_2^{m,s},r_3^{m,s},p^{m,s})\}_{s=1}^{6}$ in \eqref{eq:r123p} 
are parallel since the original set ${\cal V}_4^0$ is symmetric.
The objective
values 
in \eqref{eq:postchoice-1} are all equal to $1$ and thus 
any $s\in [S]$ can be selected.
(ii) Cut procedures.

\underline{Iteration 1.}
By solving problems \eqref{eq:gene_m_Q-1-v1-p} and \eqref{eq:gene_m_Q-2-v1-p},
we obtain $r_2^1=-0.384$, $(r_{1}^1,r_{3}^1,p^1)=(-0.5,0,0.6)$.
Then the cutting hyperplane defined as in \eqref{eq:cut_hyper} takes the form of  
\bgeq
l_{\rm cut}^1({\bm v})=({\bm R}{\bm v})^{\top}\left(G_{\bm A}(r_1^1,r_3^1,p^1)-G_{\bm B}(r_2^1)\right)
= -0.39v_1+0.39v_2.
\edeq
Since $\bbe[u^*({\bm A}_1)]= -58.97<-45.52=\bbe[u^*({\bm B}_1)]$,
the DM prefers ${\bm B}_1$ to ${\bm A}_1$. By 
\eqref{eq:choice-xkyk-0},   
$h_1=1$. The next polyhedron ${\cal V}_4^1$ 
(see \eqref{eq:choice-xkyk-0}) is the intersection of ${\cal V}_N^0$
and half-space $\{{\bm v}\in \R^2:  h_1\cdot l_{\rm cut}^1({\bm v}) \leq 0\}$,
i.e.,
${\cal V}_4^1={\cal V}_4^0\bigcap \left\{{\bm v}\in \R^{2}: 
-0.39v_1+0.39v_2\leq 0 \right\}.$
We can obtain ${\bm c}^2=(0.59,0.16)$ and  ${\bm v}_1^2=(0.90,0.10)$ and ${\bm v}_2^2=(0.23,0.23)$.

\underline{Iterations 2.} We obtain $r_2^2= -0.49$, $(r_{1}^2,r_{3}^2,p^2)=(-0.5,0.5,0.2)$, 
$l_{\rm cut}^2({\bm v})=0.14v_1-0.2v_2-0.2(1-v_1-v_2)$.
The expected utilities under $u^*$ yields
$h_2=-1$.
See the first two iterations in Figure~\ref{fig:cut12}.

\begin{figure}[htbp]
\vspace{-0.3cm}
\begin{minipage}[t]{0.95\linewidth}
\centering
\includegraphics[width=6.0in]{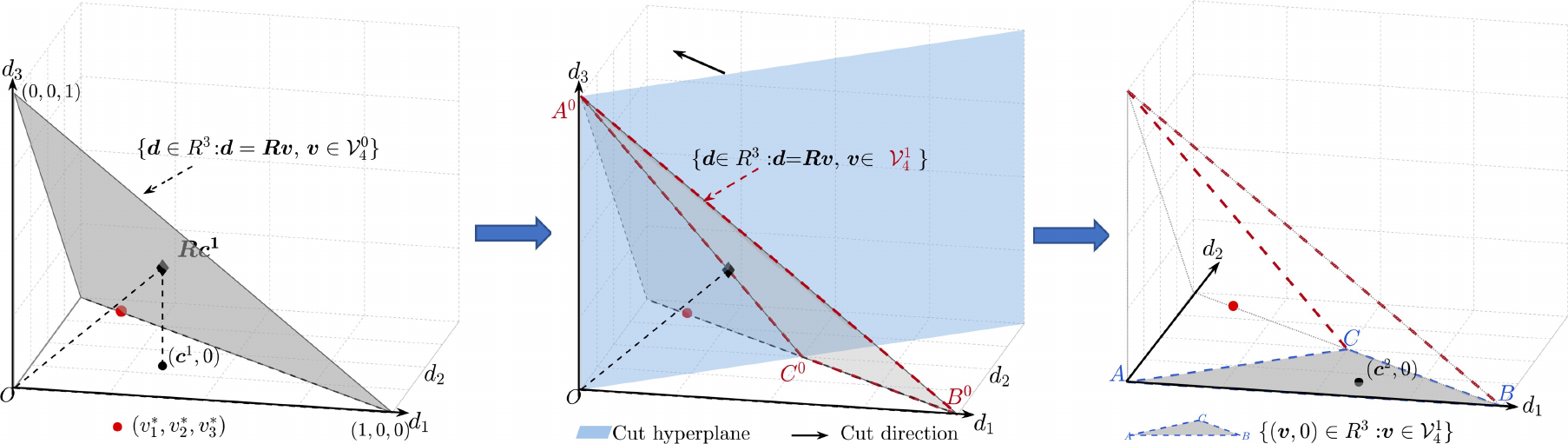}
\text{\footnotesize{\textbf{(a)}}}
\end{minipage}
\hfill
\begin{minipage}[t]{0.95\linewidth}
\centering
\includegraphics[width=6.0in]{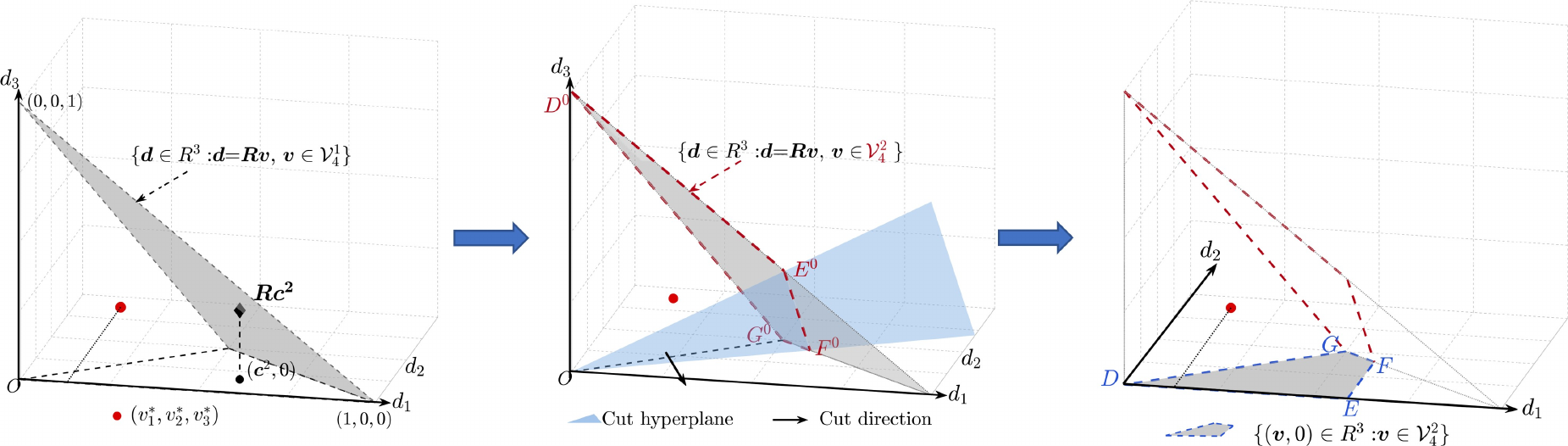}
\text{\footnotesize{\textbf{(b)}}}
\end{minipage}
\caption{\footnotesize{
(a) Visualization of the first cut and the resulting 
polyhedron ${\cal V}_4^1$. In the 
$d_1Od_2$-plane, 
the grey triangle with vertices $A$, $B$ and $C$ represents 
$\{({\bm d},0)\in \R^3: {\bm d}\in {\cal V}_4^1\}$
and the region enclosed
by the red dashed lines with vertices $A^0$, $B^0$ and $C^0$ 
represents
$\{{\bm d}\in \R^3: {\bm d}={\bm R}{\bm v},{\bm v}\in {\cal V}_4^1\}$ 
after the first cut hyperplane.
The red point is the full increment vector
$(v_1^*,v_2^*,v_3^*)=(0.18,0.81,0.01)$ with $v_i^*=\frac{u^*(x_{i+1})-u^*(x_i)}{u^*(\bar{x})-u^*(\underline{x})}$, for $i=1,2,3$,
where $u^*$ is the true utility.
(b) Visualization of the second cut and the resulting 
polyhedron ${\cal V}_4^2$.
The polyhedra
$\{({\bm v},0)\in \R^3: {\bm v}\in {\cal V}_4^2\}$
is enclosed by the red dotted lines with vertices $D$, $E$, $F$ and $G$
and the one with vertices $D^0$, $E^0$, $F^0$ and $G^0$ is the set 
$\{{\bm d}\in \R^3:{\bm d}={\bm R}{\bm v},{\bm v}\in {\cal V}_4^2\}$.}}
\label{fig:cut12}
\vspace{-0.5cm}
\end{figure}
\end{example}

An important takeaway from 
Example~\ref{exa:concave} is that 
even though the elicitation is feasible,
the two cut hyperplanes exclude the 
increment vector
${\bm v}^*$ 
associated with
the true utility function, see Figure~\ref{fig:cut12}. 
In the next subsection, we will 
explain why this occurs.

\vspace{-0.1cm}
\subsection{Direction Error of Cut Hyperplane}
\label{sec:direction error}
\vspace{-0.1cm}

The underlying reason behind the failure of the 
polyhedral cuts in Example~\ref{exa:concave} is the direction error of the cut hyperplanes. 
We next analyze this issue.

Let $u^*$ be the DM's true utility function, 
$u_N^*$ be the PLF approximation of $u^*$,
and let
${\bm v}^*$ be the increment vector associated with
$u_N^*$.
In the construction of ${\cal V}_{N}^{m}$ in \eqref{eq:choice-xkyk-0},
${\cal V}_{N}^{m}={\cal V}_{N}^{m-1}\cap {\cal M}^m$,
where
\bgeq
{\cal M}^m:=\left\{{\bm v}\in \R^{N-2}: u_N(\cdot)=({\bm R}{\bm v})^{\top}{\bm g}(\cdot),\;
h_m\cdot (\bbe[u_N({\bm A}_m)]-\bbe[u_N({\bm B}_m)])\leq 0\right\},
\edeq
with
$h_m\in \{1,-1\}$ indicating the DM's observed preference as defined in \eqref{eq:choice-xkyk-0}.
The phenomenon illustrated
in Figure~\ref{fig:cut12} can be described as follows:  
\underline{${\bm v}^*\in {\cal V}_{N}^{m-1}$ but  ${\bm v}^*\notin {\cal V}_{N}^{m}$}.
This means ${\bm v}^*\notin {\cal M}^m$,
that is,
$u_N^*(\cdot)=({\bm R}{\bm v}^*)^{\top}{\bm g}(\cdot)$ does not satisfy the constraint in ${\cal M}^m$, i.e.,
$
h_m \cdot (\bbe[u_N^*({\bm A}_m)]-\bbe[u_N^*({\bm B}_m)])> 0.
$
There are two cases:
(a) $h_m=1$ but $\bbe[u_N^*({\bm A}_m)]-\bbe[u_N^*({\bm B}_m)]>0$;
(b) $h_m=-1$ but $\bbe[u_N^*({\bm A}_m)]-\bbe[u_N^*({\bm B}_m)]<0$.
By the relation between $h_m$ and the sign of $\bbe[u^*({\bm A}_m)]-\bbe[u^*({\bm B}_m)]$,
the above two cases occur precisely when the expected-utility differences under $u^*$
and $u_N^*$
have opposite signs. Specifically,
(a) $\bbe[u^*({\bm A}_m)]-\bbe[u^*({\bm B}_m)]\leq 0$ but $\bbe[u_N^*({\bm A}_m)]-\bbe[u_N^*({\bm B}_m)]>0$;
(b) $\bbe[u^*({\bm A}_m)]-\bbe[u^*({\bm B}_m)]\geq 0$ but $\bbe[u_N^*({\bm A}_m)]-\bbe[u_N^*({\bm B}_m)]<0$.
Note that
\begin{subequations}
\begin{align}
&\bbe[u^*({\bm A}_m)]-\bbe[u^*({\bm B}_m)]=(1-p^m)u^*(r_1^m)+p^mu^*(r_3^m)-u^*(r_2^m), \label{eq:u^*}\\
&\bbe[u^*_N({\bm A}_m)]-\bbe[u^*_N({\bm B}_m)]=(1-p^m)u_N^*(r_1^m)+p^mu_N^*(r_3^m)-u_N^*(r_2^m).\label{eq:u^*_N}
\end{align}
\end{subequations}
Since $u_N^*$ is the PLF approximation of $u^*$,
we have 
$u_N^*(x_j)=u^*(x_j)$,
for all $x_j\in {\cal X}$,
where ${\cal X}$ is the set of breakpoints.
Let $r_1^m,r_3^m\in {\cal X}$. Then $u_N^*(r_i^m)=u^*(r_i^m)$, for $i=1,3$.
If, in addition, 
$r_2^m\notin {\cal X}$
and $u_N^*(r_2^m)>u^*(r_2^m)$,
then
by \eqref{eq:u^*}-\eqref{eq:u^*_N},
$
\bbe[u^*_N({\bm A}_m)]-\bbe[u^*_N({\bm B}_m)]<
\bbe[u^*({\bm A}_m)]-\bbe[u^*({\bm B}_m)],
$
case (b) may occur.
On the other hand,
if $u_N^*(r_2^m)< u^*(r_2^m)$,
then
by \eqref{eq:u^*}-\eqref{eq:u^*_N},
$
\bbe[u^*_N({\bm A}_m)]-\bbe[u^*_N({\bm B}_m)]>\bbe[u^*({\bm A}_m)]-\bbe[u^*({\bm B}_m)],
$
case (a) may occur.

Note that the direction error of the cut hyperplane differs
from the response error in the literature \citep{THG07}.
A \emph{direction error} occurs when the DM's preference between the two lotteries in a query, as represented by the true nonlinear utility function $u^*$, is inconsistent with the preference implied by its piecewise linear approximation $u_N^*$. Specifically, suppose that the DM prefers ${\bm B}_m$ to ${\bm A}_m$. Under the true utility function $u^*$, this implies 
$\mathbb{E}[u^*({\bm B}_m)] - \mathbb{E}[u^*({\bm A}_m)] > 0.$ However, under the approximated utility function $u_N^*$, it may occur that $\mathbb{E}[u_N^*({\bm B}_m)] - \mathbb{E}[u_N^*({\bm A}_m)] < 0.$ Thus, the preference directions implied by $u^*$ and $u_N^*$ are opposed.
As discussed earlier, this inconsistency arises because $\mathbb{E}[u_N^*({\bm A}_m)]$ may differ from $\mathbb{E}[u^*({\bm A}_m)]$ when the outcomes of ${\bm A}_m$ do not coincide with the breakpoints used in the piecewise linear approximation. Consequently, the piecewise linear
utility function may fail to represent the DM's true preferences accurately in certain situations.

\vspace{-0.1cm}
\subsection{The Algorithm}
\label{sec:algorithm}
\vspace{-0.1cm}

We begin 
with the set of breakpoints ${\cal X}_0=\{\underline{x},0.5*(\underline{x}+\bar{x}),0.25*\underline{x}+0.75*\bar{x},\bar{x}\}$.
Let $N_0:=|{\cal X}_0|$.
The initial polyhedron of the vector of increments is given by
$
{\cal V}_{N_0}^0:=\{{\bm v}\in \R^{N_0-2}: -{\bm v}\leq 0,{\bm e}^{\top}{\bm v}\leq 1\}\subset \R^{N_0-2}.
$
In this case, the piecewise linear utility functions can be written as 
$
u_{N_0}(x)=({\bm R}{\bm v})^{\top}{\bm g}(x),\;\; {\bm v}\in {\cal V}_{N_0}^0,
$
where ${\bm g}(x)$ is defined 
as in \eqref{eq:g(x)-incre} with ${\cal X}$ being replaced by ${\cal X}_0$ and $N$ replaced by $N_0$. The algorithm for the flexible polyhedral method is as follows.

\begin{algorithm}[!h]
\caption{Modified Polyhedral Method}
\label{alg:RUS}
\begin{algorithmic}[1]
\footnotesize 
\STATE\textbf{Input:}
Breakpoints set ${\cal X}_0$, number of queries $M$, $m=1$.
\REPEAT
\STATE 
Compute the analytic center ${\bm c}^m$ of ${\cal V}_{N_{m-1}}^{m-1}$ and the intersections ${\bm v}_1^m,{\bm v}_2^m$ between the longest axis and the boundary of ${\cal V}_{N_{m-1}}^{m-1}$.
\STATE
Solve problems \eqref{eq:gene_m_Q-1-v1-p} and \eqref{eq:gene_m_Q-2-v1-p} with a sequence of parameters $\{D_s\}_{s=1}^S$.
Select
a risky lottery ${\bm A}_m=(r_1^m,1-p^m;r_3^m)$ and a certain lottery ${\bm B}_m=r_2^m$ under postchoice symmetric criteria,
that is,
$r_i^{m}=r_i^{m,s_0}$, for $i\in [3]$, $p^m=p^{m,s_0}$,
where $s_0$ satisfies \eqref{eq:postchoice-1}.
Ask the DM to 
select between 
${\bm A}_m$ and 
${\bm B}_m$.
\STATE {\bf Update the Set of Breakpoints}.
Based on 
the lotteries ${\bm A}_m$ and ${\bm B}_m$,
update the 
set of breakpoints by
$
{\cal X}_{m}:={\cal X}_{m-1}\bigcup \{r_1^m,r_2^m,r_3^m\}
$
and set 
$N_m:=|{\cal X}_{m}|.
$
\STATE
{\bf Update the Previous $[m-1]$ Constraints}.
Use ${\cal X}_m$,
update matrices ${\bm P}$ and ${\bm Q}$ in (\ref{eq:definition-AB}),
and the
vector-valued function ${\bm g}(x)$ in (\ref{eq:g(x)-incre}) with $N=N_m$
by 
\bgeqn
\label{eq:g(t)-PRO_m}
{\bm g}^{{\cal X}_m}(x):=\Big(\underbrace{1,\cdots,1}_{i-1},\frac{x-x_i}{x_{i+1}-x_i},\underbrace{0,\cdots,0}_{N-1-i}\Big)
\in \R^{N_m-1}, \inmat{ for }x\in (x_{i},x_{i+1}],~~i\in [N_m-1].
\edeqn

Let
$
G_{\bm A}^{{\cal X}_m}(r_1^m,r_3^m,p^m):=(1-p^m){\bm g}^{{\cal X}_m}(r_1^m)+p^m {\bm g}^{{\cal X}_m}(r_3^m)
\inmat{ and }
G_{\bm B}^{{\cal X}_m}(r_2^m):={\bm g}^{{\cal X}_m}(r_2^m).
$
\STATE {\bf Update the Existing Ambiguity Set}.
${\cal V}_{N_{m-1}}^{m-1}$ is updated by
$$
{\cal V}_{N_m}^{m-1}:=\left\{{\bm v}\in \R^{N_m-2}: -{\bm v}\leq 0,{\bm e}^{\top}{\bm v}\leq 1,\;h_l\cdot ({\bm R}{\bm v})^{\top}\left(G_{\bm A}^{{\cal X}_m}(r_1^l,r_3^l,p^l)-G_{\bm B}^{{\cal X}_m}(r_2^l)\right)\leq 0, l\in [m-1]\right\},
$$ where $h_l=1$ if ${\bm B}_m$ is selected and $h_l=-1$ if ${\bm A}_m$ is selected.
\STATE
{\bf Update
the Resulting Ambiguity Set with New Constraint}. With the updated set ${\cal X}_m$
and vector-valued function ${\bm g}^{{\cal X}_m}(x)$,
update the
polyhedron
by 
\bgeqn
\label{eq:flexible-updating}
{\cal V}_{N_m}^{m} :={\cal V}_{N_m}^{m-1}\bigcap \left\{{\bm v}\in \R^{N_m-2}: h_m\cdot ({\bm R}{\bm v})^{\top}\left(G_{\bm A}^{{\cal X}_m}(r_1^m,r_3^m,p^m)-G_{\bm B}^{{\cal X}_m}(r_2^m)\right)\leq 0\right\},
\edeqn
where $h_l=1$ if ${\bm B}_m$ is selected and $h_l=-1$ if ${\bm A}_m$ is selected.
Update $m := m+1$.
\UNTIL $m=M$.
\hspace{-0.cm}\textbf{Output:}
Polyhedron ${\cal V}_{N_M}^M$.
\end{algorithmic}  
\label{alg:FPM}
\end{algorithm}

The intuition behind Algorithm~\ref{alg:FPM} is 
that by adding new points $r_1^m,r_2^m,r_3^m$ into the set of the breakpoints
${\cal X}_m$ at iteration $m$ 
whenever they are not already included,
we have $u_{N_m}^*(r_i^m)=u^*(r_i^m)$ for $i=1,2,3$. It follows from \eqref{eq:u^*}-\eqref{eq:u^*_N} that
$\bbe[u^*({\bm A}_m)]-\bbe[u^*({\bm B}_m)]
=\bbe[u^*_{N_m}({\bm A}_m)]-\bbe[u^*_{N_m}({\bm B}_m)]$.
This ensures that $h_m$ and $\bbe[u^*_{N_m}({\bm A}_m)]-\bbe[u^*_{N_m}({\bm B}_m)]$ have the same sign. Since $u^*_{N_m}=({\bm R}{\bm v}^*)^\top{\bm g}$,
we have ${\bm v}^*\in{\cal M}^m$, thereby preventing the direction error discussed in Section~\ref{sec:direction error}.
After
implementing 
$M$ 
queries
as outlined in the algorithm,
we have ``narrowed" 
the 
polyhedron of ${\bm v}$ from ${\cal V}_{N_0}^0$ to ${\cal V}_{N_{M}}^M\subset {\cal V}_{N_{M}}^0$.
It is important to note that ${\cal V}_{N_0}^0\in \R^{N_0-2}$, whereas  ${\cal V}_{N_{M}}^0\in \R^{N_M-2}$.
Although they do not lie in the same space,
we can measure the convergence by calculating the 
metrics
of the feasible set ${\cal V}_{N_M}^M$ of increment vectors or the corresponding set ${\cal U}_{N_M}^M$ of utility functions.
We will 
discuss this
in detail in Section \ref{sec:quality}.
At this point,
we have a family of piecewise linear utility functions
\bgeqn
\label{eq:PLA_flexible}
{\cal U}_{N_M}^M:=\left\{u_{N_{M}}(x)=({\bm R}{\bm v})^{\top}{\bm g}^{{\cal X}_{M}}(x):\;\; {\bm v} \in {\cal V}_{{N_{M}}}^{M}\right\}.
\edeqn

In the numerical experiments, we restrict $p$ to the interval $[\underline{p}, \bar{p}] \subset (0,1)$ to avoid the extreme cases $p=0$ and $p=1$. The choice of the initial set of breakpoints ${\cal X}_0$ influences the performance of the algorithm. Since the number of breakpoints determines the dimensionality of problem \eqref{eq:gene_m_Q},
we consider a small initial set.
The following example illustrates the implementation of Algorithm~\ref{alg:RUS}.

\begin{example}
\label{ex:revi-4.1}
We consider the true S-shaped utility function
$u^*(x)=(2(e^{8x}-1))/{8}$, for $x\in [-0.5,0]$
and 
$u^*(x)=(1-e^{-3x})/{3}$ for $x\in[0,0.5]$.
We set  $[\underline{p},\bar{p}]=[0.05,0.95]$ and retain 
two decimal places, and
the breakpoints set ${\cal X}_0=\{-0.5,0,0.25,0.5
\}$. Then $N_0=4$,
and the initial polyhedron is 
${\cal V}_{4}^0=\{{\bm v}\in \R^2: v_i\geq 0,i\in[2],~~ {\bm e}^{\top}{\bm v}\leq 1\}$.
$D_s=\frac{s}{10}$, $s\in[10]$.

\underline{Iteration 1}.
By solving problems 
\eqref{eq:gene_m_Q-1-v1-p} and \eqref{eq:gene_m_Q-2-v1-p},
we obtain $(r_1^1,r_2^1,r_3^1)=(-0.5,-0.05,0.5)$ and $p^1=0.3$.
Since $\bbe[u^*({\bm A}_1)]=0.3*(-0.245)+0.7*0.259=0.108>-0.010=\bbe[u^*({\bm B}_1)]$.
The set of breakpoints is updated by 
${\cal X}_1=\{-0.5,-0.05,0,0.25,0.5\}$ with $N_1=|{\cal X}_1|=5$,
and 
the initial polyhedron is updated by
$
{\cal V}_{5}^0=\{{\bm v}\in \R^{3}:v_i\geq 0, i\in [3], {\bm e}^{\top}{\bm v} \leq 1\}.
$
Then by \eqref{eq:g(t)-PRO_m},
$
{\bm g}^{{\cal X}_1}(r_1^1)
=(0,0,0,0)$,
$
{\bm g}^{{\cal X}_1}(r_2^2)=
(0.9,0,0,0)$,
and 
$
{\bm g}^{{\cal X}_1}(r_3^2)
=(1,1,1,1).
$
\underline{Iterations 2 and 3}.
We obtain $(r_1^2,r_2^2,r_3^2)=(-0.05,0.34,0.5)$, $p^2=0.79$,
and
${\cal X}_2={\cal X}_1
\cup \{0.34\}$.
Next,
$(r_1^3,r_2^3,r_3^3)=(-0.5,-0.30,0.5)$ and $p^3=0.20$.
Updating the set of breakpoints  ${\cal X}_3$ and ${\bm g}^{{\cal X}_3}(t)$,
we can continue the procedure.
 Let 
\bgeqn
\label{eq:under_upper}
\underline{u}_i:=\min_{{\bm v}\in {\cal V}_{N_m}^m} ({\bm R}{\bm v})^{\top}{\bm g}^{{\cal X}_m}(x_i)
\quad
\inmat{and}
\quad 
\bar{u}_i:=\max_{{\bm v}\in {\cal V}_{N_m}^m} ({\bm R}{\bm v})^{\top}{\bm g}^{{\cal X}_m}(x_i), \inmat{ for }i\in [N_m],
\edeqn
and $\underline{u}:[\underline{x},\bar{x}]\to [0,1]$ be the PLF passing through 
 points $(x_i,\underline{u}_i)$, for $i\in [N_m]$,
and $\bar{u}:[\underline{x},\bar{x}]\to [0,1]$ be the PLF passing through $(x_i,\bar{u}_i)$, for $i\in [N_m]$.
$[\underline{u}_i,\bar{u}_i]$ specifies the
\underline{range} of the function values of
$u_{N_m}\in {\cal U}_{N_m}^m$ at breakpoints $x_i$ 
for $i\in [N_m]$. 

{\bf Impact of the Initial Breakpoints ${\cal X}_0$}.
We consider three cases:
(a) $N_0=4$ with ${\cal X}_0=\{-0.5,0,0.25,0.5\}$;
(b) $N_0=5$ with ${\cal X}_0=\{-0.5,-0.25,0,0.25,0.5\}$;
(c) $N_0=6$ with ${\cal X}_0=\{-0.5,-0.25,0,0.12,0.25,0.5\}$
with $[\underline{p},\bar{p}]=[0.05,0.95]$ and retaining
two decimal places.
Figure~\ref{fig:PLA_Flex_1} depicts 
$u_N^{\rm AC}$, 
the range of 
${\cal U}_{N_M}^M$, 
and the CPU time as 
$N_0$ and
$M$ 
change.
Observe that the width of the pointwise utility-function range is nonuniform across its domain. This behavior is driven by certain extreme vertices of the polyhedron ${\cal V}_{N_M}^M$. To illustrate, consider the initial polyhedron ${\cal V}_4^0$ associated with ${\cal X}_0=\{x_1,x_2,x_3,x_4\}$ and $N_0=4$. In this case, ${\cal V}_4^0$ consists of two vertices, $A=(1,0,0)$ and $B=(0,0,1)$, 
which correspond to the following piecewise linear utility functions: 
\bgeqn 
u_A(x)=\left\{ \begin{array}{ll}
\frac{x-x_1}{x_2-x_1} & \text{for } x\in[x_1,x_2]\\
1 & \text{for } x\in[x_2,x_4] \end{array}
\right.
\quad \text{and} \quad 
u_B(x)=\left\{ \begin{array}{ll} 0 & \text{for } x\in[x_1,x_3]\\ \frac{x-x_3}{x_4-x_3} & \text{for } x\in[x_3,x_4] \end{array} \right..
\edeqn 
For any $x\in[x_2,x_4]$, $\sup_{u\in{\cal U}_4^0} u(x)=u_A(x)=1,$ while for any $x\in[x_1,x_3]$, $ \inf_{u\in{\cal U}_4^0} u(x)=u_B(x)=0. $
Consequently, the lower bound of ${\cal U}_4^0$ is identically $0$ on $[x_1,x_3]$, whereas its upper bound is $1$ on $[x_2,x_4]$. 
Since $x_2,x_3\in [x_2,x_4]\cap [x_1,x_3]$, it follows from \eqref{eq:under_upper} that
$[\underline{u}_i,\bar{u}_i]=[0,1]$ for $i=2,3$.
Therefore, the utility function range attains its maximum possible width over substantial portions of the domain.

{\bf Impact of the Range $[\underline{p},\bar{p}]$}.
We test three examples with case (a) $[\underline{p},\bar{p}]=[0.06,0.94]$, (b) $[\underline{p},\bar{p}]=[0.05,0.95]$, 
(c) $[\underline{p},\bar{p}]=[0.04,0.96]$
 with ${\cal X}_0=\{-0.5,0,0.25,0.5\}$ and
retaining two decimal places.
Figure~\ref{fig:PLA_Flex_2} depicts 
$u_N^{AC}$, 
the range of 
${\cal U}_{N_M}^M$, 
and the CPU time as 
$[\underline{p},\bar{p}]$ and
$M$ 
change.
We observe from Figure~\ref{fig:PLA_Flex_2} that, during the first 10 queries, a greater 
range for
\(p\) does not necessarily yield a tighter utility range, as the early-stage shrinkage depends on the specific cuts generated by the adaptive procedure. As more queries are conducted, a larger
range of \(p\) provides greater flexibility in generating cuts. Consequently, at \(M=50\), a wider \([\underline{p},\bar{p}]\) may result in a smaller utility range. 
We set $[\underline{p},\bar{p}]=[0.05,0.95]$ to avoid extreme lottery
probabilities while allowing a broad range of $p$.

\begin{figure}[!h]
\vspace{-0.2cm}
\centering
\begin{minipage}[t]{0.32\linewidth}
\centering
\includegraphics[width=2.1in]{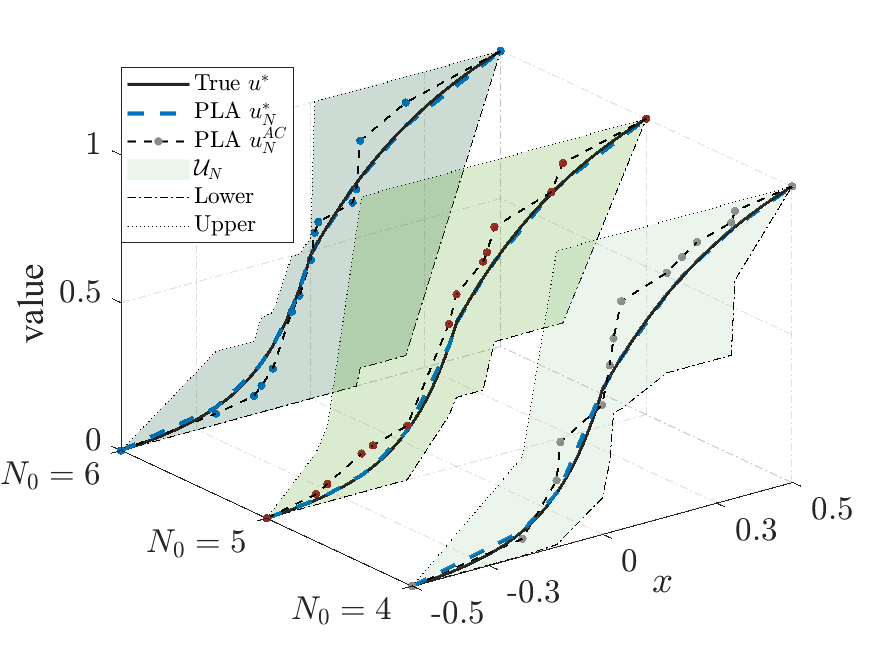}
\text{\footnotesize{\textbf{\hspace{-0.3cm}(a)}}}
\end{minipage}
\begin{minipage}[t]{0.32\linewidth}
\centering
\includegraphics[width=2.1in]{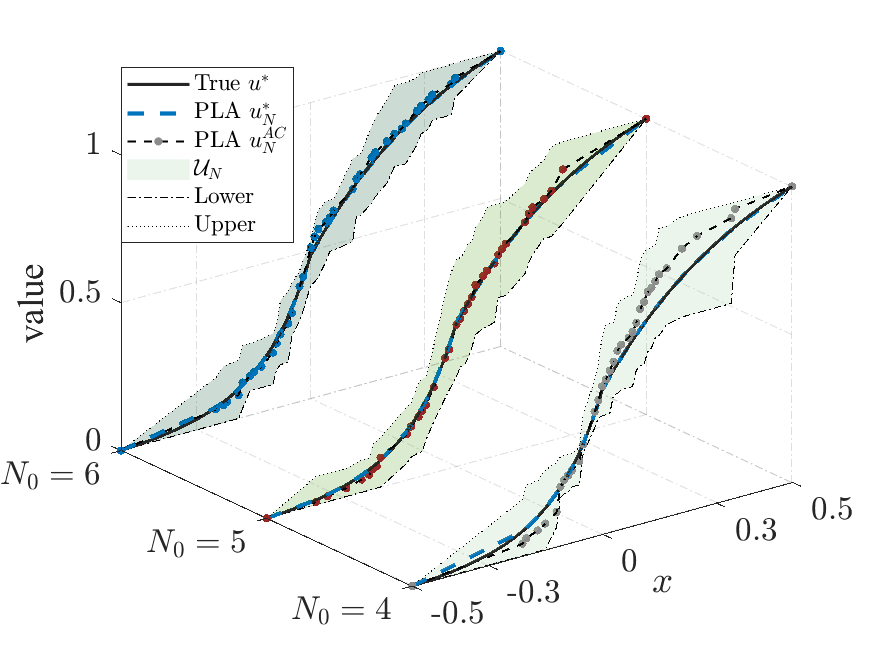}
\text{\footnotesize{\hspace{-0.3cm}\textbf{(b)}}}
\end{minipage}
\begin{minipage}[t]{0.32\linewidth}
\centering
\includegraphics[width=2.1in]{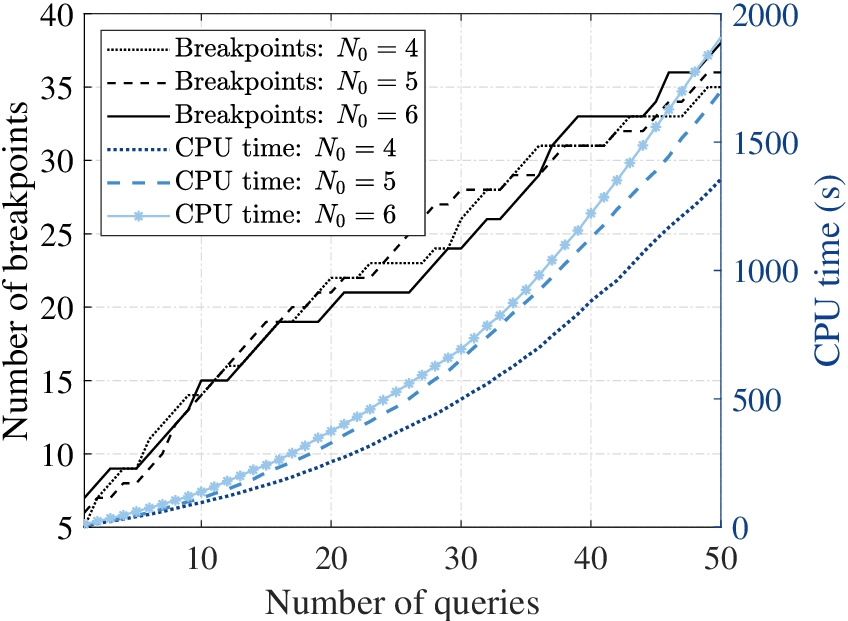}
\text{\footnotesize{\hspace{-0.2cm}\textbf{(c)}}}
\end{minipage}
\caption{\footnotesize{(a) Utility function range with $M=10$. 
(b) Utility function range with $M=50$.
The dotted curves in (a) and (b)
correspond to the piecewise linear utility functions
constructed with analytic centers
${\bm c}^{M}$,
the black curve represents the PLF $u_N^*$.
The green area specifies the range  $[\underline{u}_i,\bar{u}_i]$
for $i\in [N]$.
(c) CPU time and number of breakpoints. The blue curves depict the change of  
CPU time (seconds) as the number of queries $M$ increases.
The black curves  
show changes in
the number of breakpoints ($N_M=|{\cal X}_M|$) as $M$ varies.
}}
\label{fig:PLA_Flex_1}
\vspace{-0.4cm}
\end{figure}

\begin{figure}[!ht]
\centering
\begin{minipage}[t]{0.32\linewidth}
\centering
\includegraphics[width=2.1in]{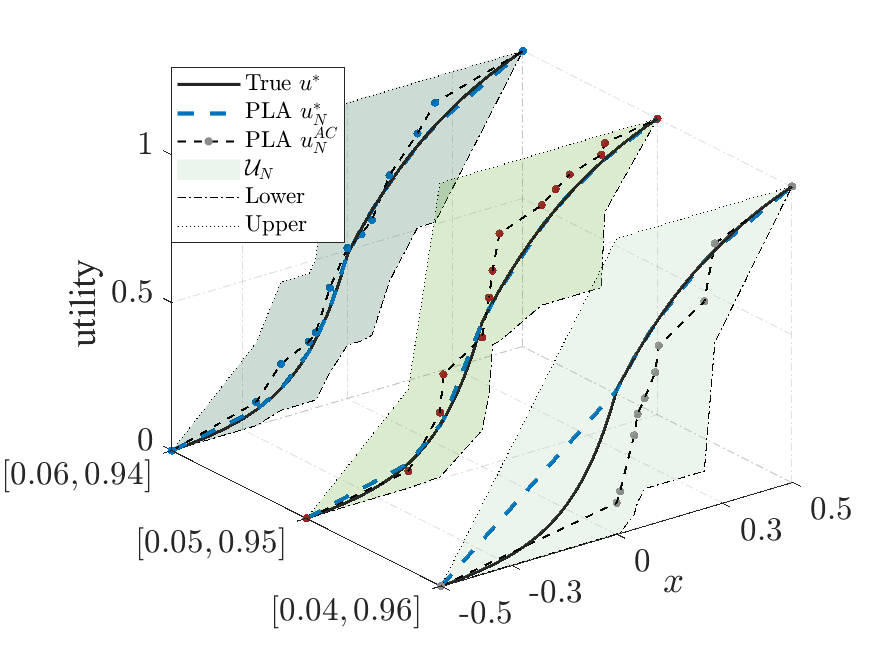}
\text{\footnotesize{\textbf{(a)}}}
\end{minipage}
\begin{minipage}[t]{0.32\linewidth}
\centering
\includegraphics[width=2.1in]
{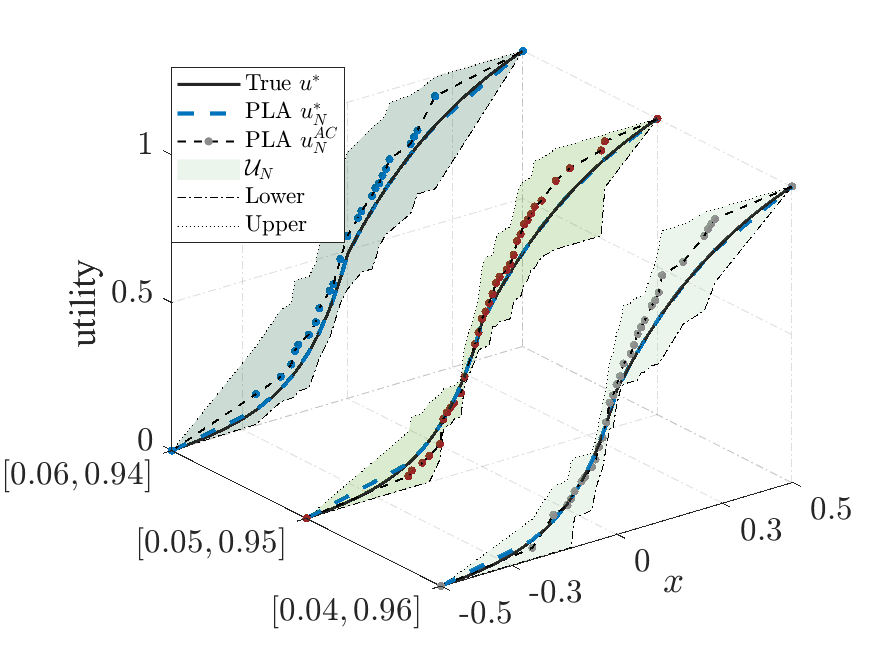}
\text{\footnotesize{\textbf{(b)}}}
\end{minipage}
\begin{minipage}[t]{0.32\linewidth}
\centering
\includegraphics[width=2.1in]
{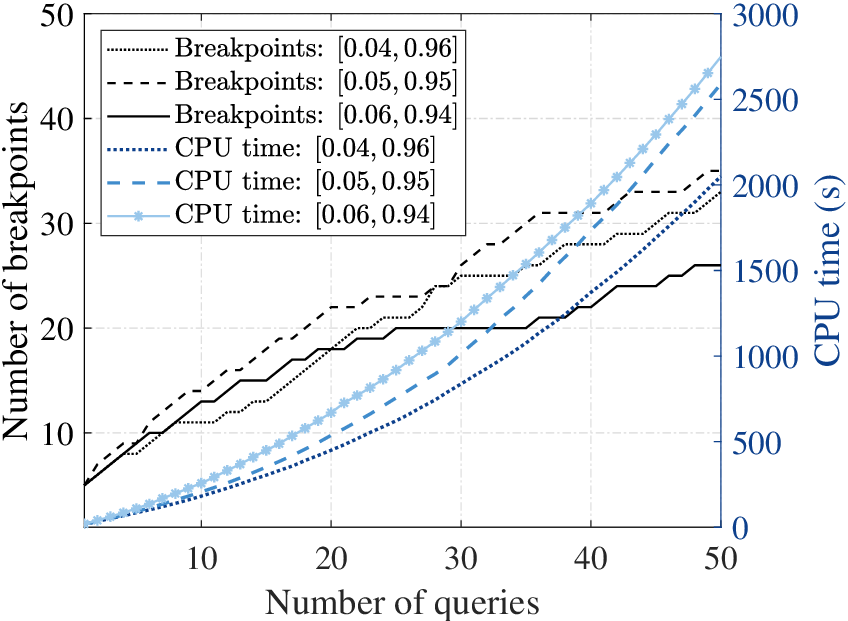}
\text{\footnotesize{\hspace{-0.2cm}\textbf{(c)}}}
\end{minipage}
\caption{\footnotesize{(a) Utility function range $M=10$. (b) Utility function range $M=50$. (c) CPU time and number of
breakpoints. The labels, legends, and curves have the 
same
interpretation as in Figure~\ref{fig:PLA_Flex_1}.
}}
\label{fig:PLA_Flex_2}
\vspace{-0.2cm}
\end{figure}

{\bf Impact of Retaining Decimal Places}.
We test three cases where the set of breakpoints is rounded to one decimal ($0.1$), two decimals ($0.01$), and three decimals ($0.001$),
with ${\cal X}_0=\{-0.5,0,0.25,0.5\}$
and $[\underline{p},\bar{p}]=[0.05,0.95]$. 
Figure~\ref{fig:PLA_Flex_3} depicts 
$u_N^{\rm AC}$, 
the range of 
${\cal U}_{N_M}^M$, 
and the CPU time as 
the decimals of $\underline{p}$ and $\bar{p}$
in $[\underline{p},\bar{p}]$ and $M$ 
change.
We observe from Figure~\ref{fig:PLA_Flex_3} that the utility range shrinks particularly rapidly when the breakpoints are rounded to one decimal place. In this case, at most 11 distinct breakpoints can be generated, so the dimension of the increment vector is at most \(N_M-2\leq 9\). Thus, for a large number of queries, such as \(M=50\), many cuts are imposed on a relatively low-dimensional polyhedron, resulting in a particularly narrow utility range. Retaining more decimal places allows more distinct breakpoints and hence a higher-dimensional polyhedron, for which the same number of queries results in less shrinkage. To balance the number of breakpoints and numerical precision, we retain
two decimal places in the subsequent experiments.

\begin{figure}[!h]
\centering
\begin{minipage}[t]{0.32\linewidth}
\centering
\includegraphics[width=2.1in]{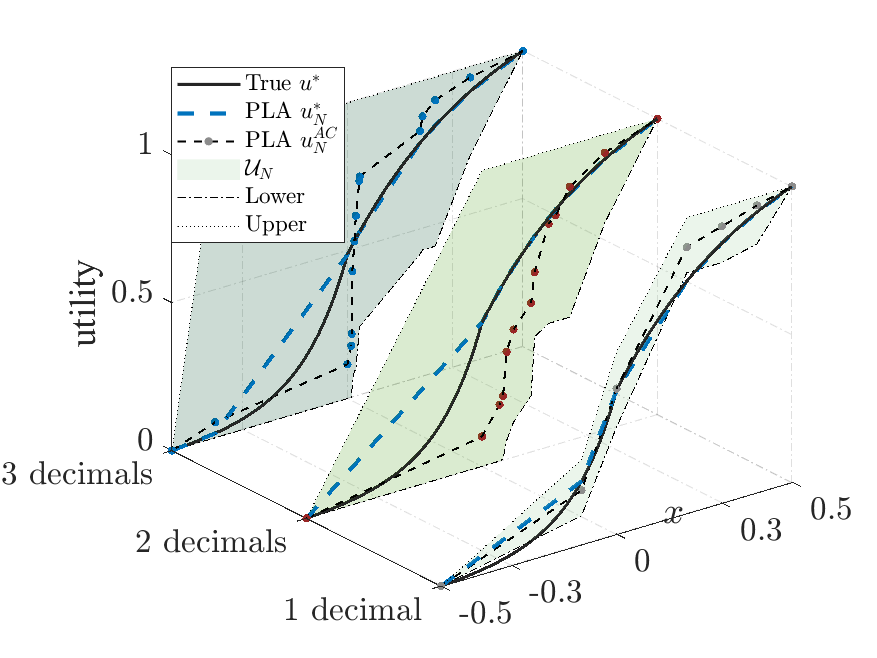}
\text{\footnotesize{\textbf{(a)}}}
\end{minipage}
\begin{minipage}[t]{0.32\linewidth}
\centering
\includegraphics[width=2.1in]
{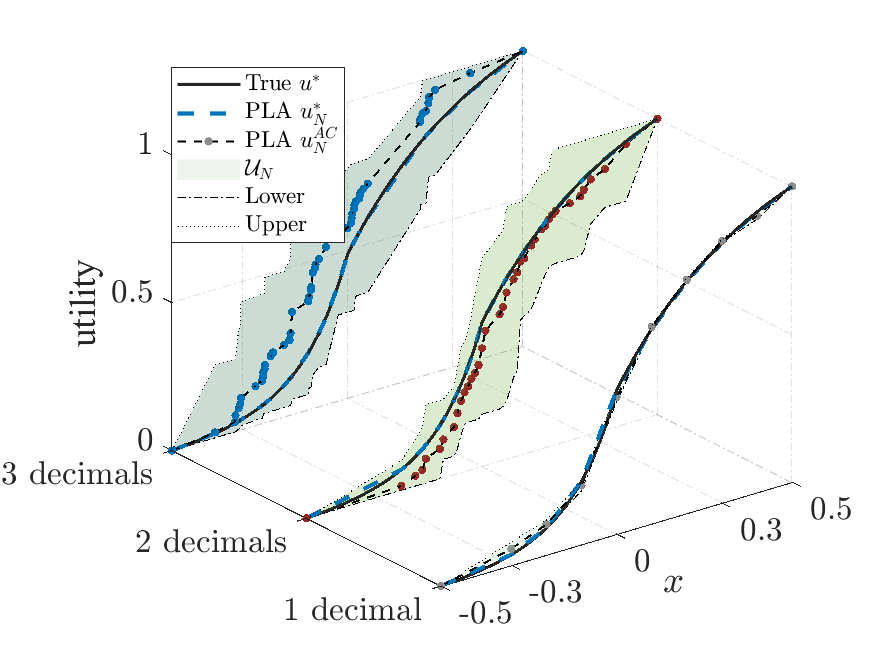}
\text{\footnotesize{\textbf{(b)}}}
\end{minipage}
\begin{minipage}[t]{0.32\linewidth}
\centering
\includegraphics[width=2.1in]
{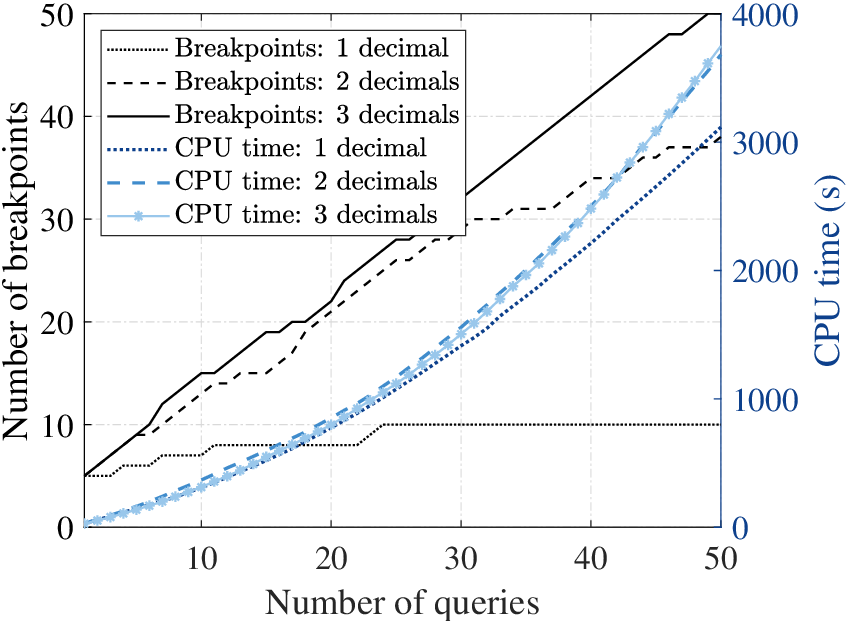}
\text{\footnotesize{\hspace{-0.2cm}\textbf{(c)}}}
\end{minipage}
\caption{\footnotesize{(a) Utility function range $M=10$. (b) Utility function range $M=50$. (c) CPU time and number of
breakpoints.
The labels, legends, and curves have the 
same
interpretation as in Figure~\ref{fig:PLA_Flex_1}.
}}
\label{fig:PLA_Flex_3}
\vspace{-0.2cm}
\end{figure}
\end{example}

The above discussion motivates us to further investigate the following issues: (a) how new breakpoints are generated; (b) how to effectively eliminate extreme points that underpin a wide range of utility functions; and (c) how the dimensionality of the increment vector affects the convergence behavior of the MPM. We address these issues within a broader framework that examines the practical applicability of the proposed MPM.

\vspace{-0.1cm}
\subsection{Practical Efficacy of the Proposed Method for Eliciting Preferences}
\label{sec:complexity}
\vspace{-0.1cm}

In this section, we provide insights into adaptive breakpoint generation in Algorithm~\ref{alg:RUS} and the advantages of MPM
over RRUS.

\vspace{-0.1cm}
\subsubsection{Adaptive Breakpoint Generation (Issue (a))}
\label{sec:advantage}
\vspace{-0.1cm}

We begin by analyzing 
the characteristics and effects of the newly generated breakpoints in Algorithm~\ref{alg:RUS} 
when information about the shape and Lipschitz modulus of the utility function is unavailable.
In Algorithm~\ref{alg:RUS}, new breakpoints may be generated by
$r_1^m$, $r_2^m$, or $r_3^m$. In our numerical experiments, however,
almost all generated $r_1^m$ and $r_3^m$ already belong to
${\cal X}_{m-1}$, even before
imposing this restriction.
Moreover, the equality
$
({\bm R}{\bm c}^m)^\top
G_{\bm A}(r_1^m,r_3^m,p^m)=D
$
in problem~\eqref{eq:gene_m_Q-2-v1-p} can generally be satisfied by
adjusting $p^m$ whenever $D$ lies within its admissible range $[\underline{p},\bar{p}]$.
Hence, newly generated breakpoints arise primarily from $r_2^m$.

Under Criterion~\ref{eq:criterion_2}
in Algorithm~\ref{alg:RUS}, $r_2^m$ is selected from candidate breakpoints $r_2^{m,s}$ satisfying 
\bgeqn 
\label{eq:candidate_breakpoints} 
({\bm R}{\bm c}^m)^\top G_{\bm B}(r_2^{m,s}) =D_s=\frac{s}{S}, \;\;\; s\in[S], 
\edeqn 
where $S=10$ in our numerical experiments. The candidate breakpoints therefore correspond to equally spaced utility levels of the piecewise linear utility approximation $({\bm R}{\bm c}^m)^\top{\bm g}(\cdot)$ constructed using the analytic center ${\bm c}^m$. As the number of queries $m$ increases, $({\bm R}{\bm c}^m)^\top{\bm g}(\cdot)$ converges to the piecewise linear approximation $u_N^*(\cdot)$ of the true utility function. Consequently, the breakpoints identified through \eqref{eq:candidate_breakpoints}  correspond approximately to equally spaced utility levels of $u_N^*$, namely $s/S$ for $s\in [S]$. Since equal utility increments span shorter intervals in regions where $u_N^*$ is steeper, more breakpoints are allocated to such regions, resulting in a closer approximation where the utility function varies more rapidly.
We use a simple example to illustrate this. 

\begin{example}[Adaptive Breakpoint Generation]
\label{exa:concave_breakpoints}
Consider the concave utility function in Example~\ref{exa:concave} 
and  ${\cal X}_0=\{-0.5,0,0.25,0.5\}$.
Figure~\ref{fig:Fig1-18} depicts 
representative
stages of the elicitation
process for $M=1,9,$ and $18$. As the elicitation proceeds, the
analytic-center approximation $({\bm R}{\bm c}^m)^\top{\bm g}(\cdot)$
approaches
the PLF
approximation $u_N^*(\cdot)$ of the true utility function, and the newly generated
breakpoints are
concentrated in its steeper region. In particular,
from $m=5$ onward, all subsequently generated breakpoints $r_2^m$ are
located at or to the left of $-0.04$. 
This numerical observation
illustrates how the MPM adaptively allocates
more breakpoints to regions where the utility function is steeper.

\begin{figure}[!htbp]
\vspace{-0.2cm}
\centering
\begin{minipage}[t]{0.32\linewidth}
\centering
\includegraphics[width=2.2in]{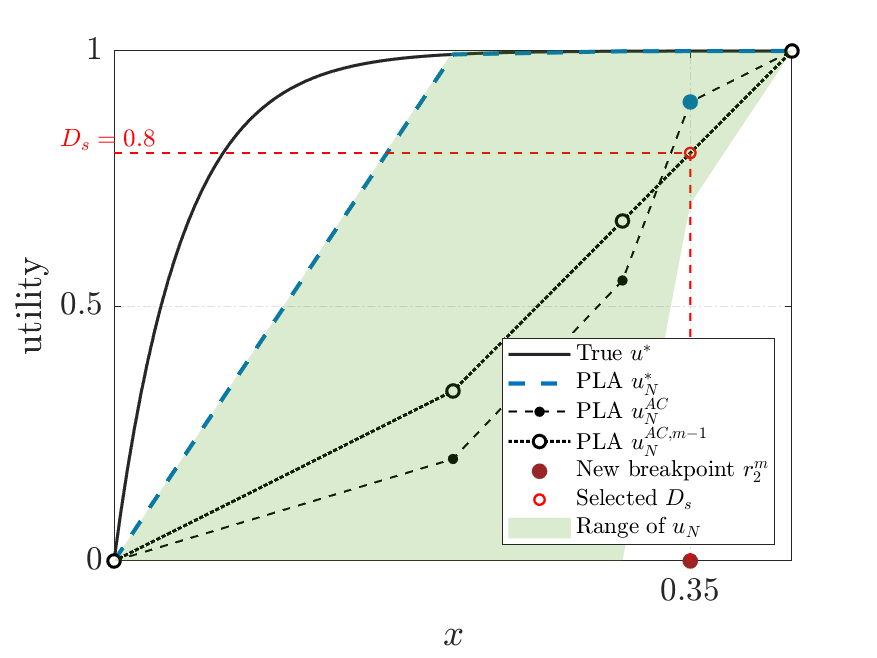}
\text{\footnotesize{\hspace{0.5cm}\textbf{(a)}}}
\end{minipage}
\hfill
\begin{minipage}[t]{0.32\linewidth}
\centering
\includegraphics[width=2.2in]{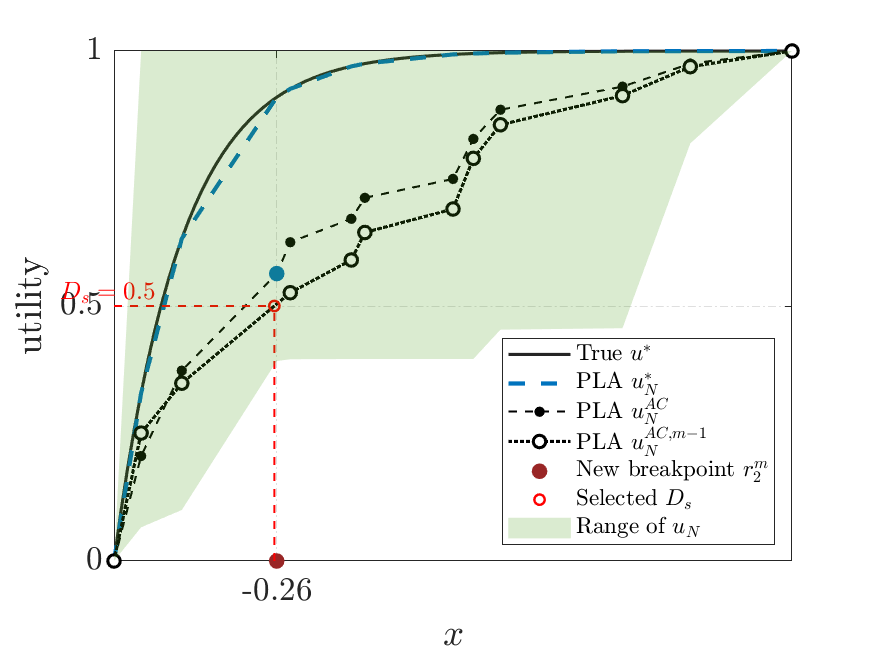}
\text{\footnotesize{\hspace{0.5cm}\textbf{(b)}}}
\end{minipage}
\hfill
\begin{minipage}[t]{0.32\linewidth}
\centering
\includegraphics[width=2.2in]{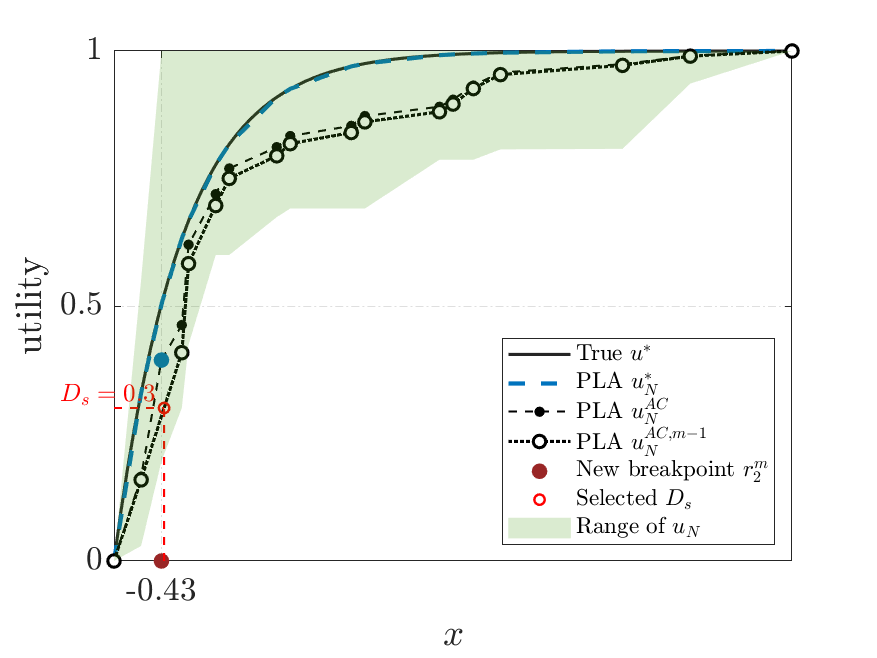}
\text{\footnotesize{\hspace{0.5cm}\textbf{(c)}}}
\end{minipage}
\caption{\footnotesize 
Evolution of the estimated utility function based on analytic center, utility function range, and 
newly generated breakpoints for the concave utility function with query number:
(a) $M=1$, (b) $M=9$, and (d) $M=18$.
The green area represents the pointwise utility-function range, and the red dashed 
lines indicate the level $D_s$ and its corresponding location, at which 
the red point marks the newly generated breakpoint.
}
\label{fig:Fig1-18}
\vspace{-0.5cm}
\end{figure}
\end{example}

\vspace{-0.1cm}
\subsubsection{Applicability under Limited Structural Information}
\label{sec:information}
\vspace{-0.1cm}

The main advantage of the MPM is that its implementation does not require prior knowledge of the shape or the Lipschitz modulus of the underlying utility function. This feature is particularly relevant when such structural information is unavailable or difficult to specify reliably.

To illustrate this advantage, we compare Algorithm~\ref{alg:RUS} with the RRUS scheme studied in \cite{AmD15}. At each iteration \(m\), RRUS randomly generates two outcomes \(r_1^m<r_3^m\) with $r_2^m=(r_1^m+r_3^m)/2$ and asks the DM to compare the deterministic outcome \(r_2^m\) with a binary lottery yielding \(r_1^m\) and \(r_3^m\) with probabilities \(1-p^m\) and \(p^m\), respectively. To determine \(p^m\), RRUS normalizes the utility values at \(r_1^m\) and \(r_3^m\) to \(0\) and \(1\), computes the lower and upper bounds on the feasible relative utility value at \(r_2^m\) under the current ambiguity set, and sets \(p^m\) to the midpoint of these two bounds. The resulting preference response is then used to shrink the ambiguity set.

As shown in the Electronic Companion of \cite{ZDX23_b}, the effectiveness of RRUS can be affected substantially by the availability of structural information on the utility function. In particular, when neither shape information nor a Lipschitz modulus is available, the feasible range of the relative utility at \(r_2^m\) can remain wide, limiting the informativeness of the resulting queries. 
By contrast, the MPM designs its queries adaptively based on the geometry of the current ambiguity set and does not require either type of structural information.
To highlight the differences in performance
of RRUS and Algorithm~\ref{alg:RUS}, we apply them to utility function in Example~\ref{exa:concave_breakpoints}.

\begin{example}[Impact of Structural Information: Continuation of Example~\ref{exa:concave_breakpoints}]
\label{exa:RRUS_comparison}
We compare
the performance of RRUS and Algorithm~\ref{alg:RUS} under different levels of structural information about the utility function.
Figure~\ref{fig:Fig_RRUS_PM50}  depicts
how the effectiveness of the methods depends on 
information about
shape and Lipschitz modulus.
Figures~\ref{fig:Fig_RRUS_PM50}(a) and (b) show that RRUS can reduce the ambiguity effectively when either the shape or Lipschitz modulus is known. 
Figure~\ref{fig:Fig_RRUS_PM50}(c) 
shows that Algorithm~\ref{alg:RUS}
can shrink the ambiguity set without such information.

\begin{figure}[!h]
\vspace{-0.2cm}
\centering
\begin{minipage}[t]{0.32\linewidth}
\centering
\includegraphics[width=2.0in]{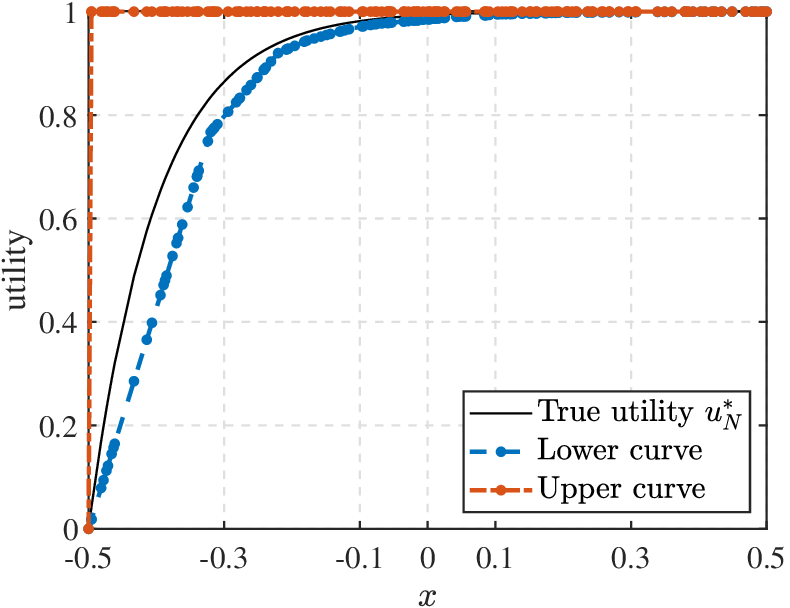}
\text{\footnotesize{\hspace{0.4cm}\textbf{(a)}}}
\end{minipage}
\hfill
\begin{minipage}[t]{0.32\linewidth}
\centering
\includegraphics[width=2.0in]{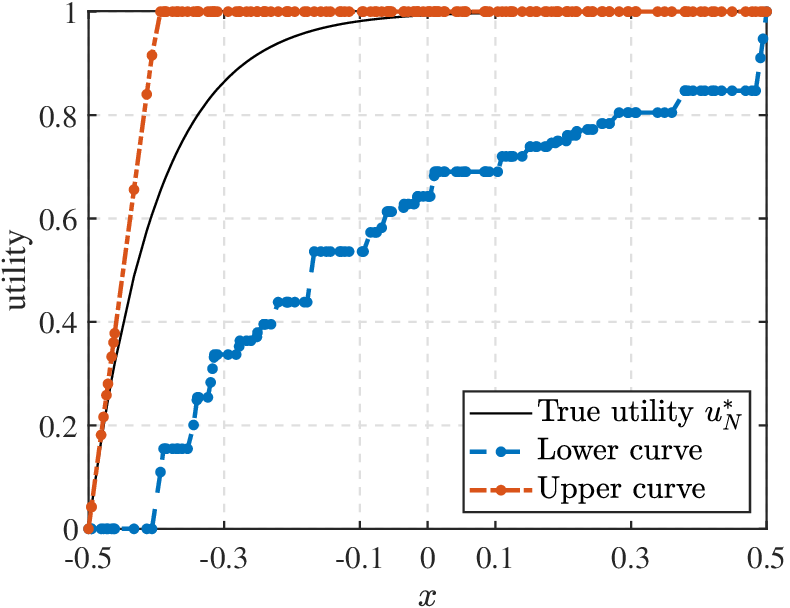}
\text{\footnotesize{\hspace{0.4cm}\textbf{(b)}}} 
\end{minipage}
\hfill
\begin{minipage}[t]{0.32\linewidth}
\centering
\includegraphics[width=2.0in]{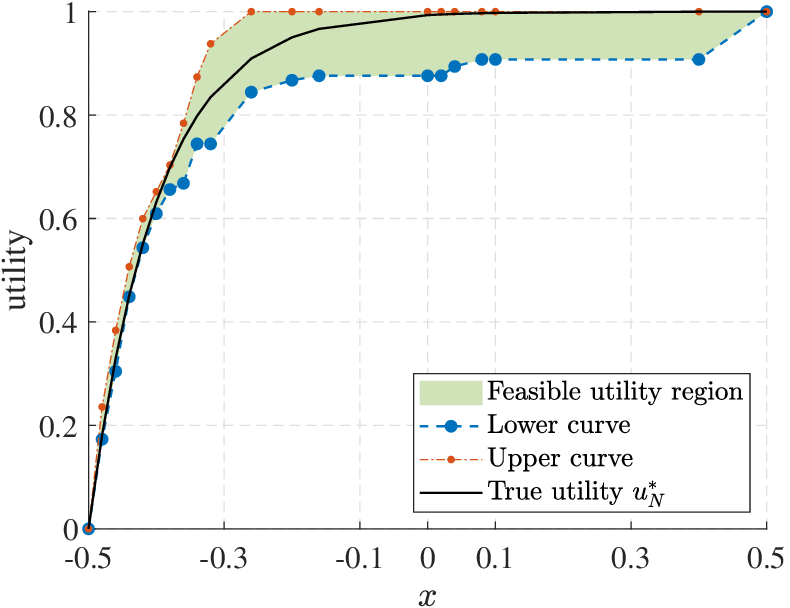}
\text{\footnotesize{\hspace{0.4cm}\textbf{(c)}}}  
\end{minipage}
\caption{\footnotesize{Comparison of the utility-function ranges under different structural information 
with $M=50$: 
(a) RRUS with known shape and unknown Lipschitz modulus;
(b) RRUS with known Lipschitz modulus $L$ and unknown shape;
(c) the MPM with neither shape nor Lipschitz modulus known.
The black curve represents the true utility function, while the lower and upper 
curves represent the resulting pointwise utility-function range.
}}
\label{fig:Fig_RRUS_PM50}
\vspace{-0.5cm}
\end{figure}
\end{example}

\vspace{-0.1cm}
\subsubsection{Quantifying the Reduction 
of Ambiguity Set}
\label{sec:quality}
\vspace{-0.1cm}

We now discuss how to measure 
reduction of the polyhedron and the corresponding ambiguity set of utility functions.  Specifically, we
quantify the size of ${\cal V}_{N_M}^M$ and the associated set of utility functions ${\cal U}_{N_M}^M$, as well as 
the quality of 
the approximate utility function constructed via 
the analytic center of ${\cal V}_{N_M}^M$.

In the literature of
multi-attribute linear utility function,
there are several measures to quantify the quality of estimation of $u_N^*$,
such as $D$-error,
root mean squared error, hit rate, and mean absolute error (see \cite{SaV19,McC02,KTG94}).
Here, 
we can measure the efficiency of the algorithm through the sup norm between  ${u}_{N_M}^*(\cdot)$ and the piecewise linear function $({\bm R}{\bm c}^M)^{\top}{\bm g}^{{\cal X}_M}(\cdot)$ 
and the radius of the polyhedron 
${\cal V}_{N_M}^{M}$
and the radius of set  ${\cal U}_{N_M}^{M}$.

Let 
$
u_{N_m}^{\rm AC}(\cdot):=({\bm R}{\bm c}^m)^{\top}{\bm g}^{{\cal X}_m}(\cdot) 
$
and 
$
u_{N_m}^*(\cdot)=({\bm R}{\bm v}^*)^{\top}{\bm g}^{{\cal X}_m}(\cdot),
$
where ${\bm c}^m\in \R^{N_m-2}$ is the analytic center of ${\cal V}_{N_m}^{m-1}$
and ${\bm v}^*\in \R^{N_{m}}$ is
the vector of 
increments
with $v_i^*=u^*(x_{i+1})-u^*(x_i)$, for $i\in [N_{m}-1]$.
The quality of the elicitation is measured from two aspects:
(a) 
The distance between $u_{N_m}^{\rm AC}(\cdot)$ and $u_{N_m}^*(\cdot)$:
$\dd_{\rm AC}
:= \sup_{i\in[N_m]} |({\bm R}{\bm c}^m)^{\top}{\bm g}^{{\cal X}_m}(x_i)- ({\bm R}{\bm v}^*)^{\top}{\bm g}^{{\cal X}_m}(x_i)|
= \sup_{i\in[N_m-2]} |\sum_{j=1}^i{c}^m_j -\sum_{j=1}^i{v}^*_j|$.
(b) 
The radius of polyhedron ${\cal V}_{N_m}^{m}$:
$
\dd_{\rm R_1}
:=\sup_{{\bm v}_1^m,{\bm v}_2^m\in {\cal V}_{N_m}^m}\|{\bm v}_1^m-{\bm v}_2^m\|_{\infty},
$
and the radius of set ${\cal U}_{N_m}^{m}$:
$\dd_{\rm R_2}
:=\sup_{x\in [\underline{x},\bar{x}]}|\underline{u}(x)-\bar{u}(x)|
= \sup_{i\in [N_m]}|\underline{u}_i-\bar{u}_i|,
$
where $\underline{u}$ and $\bar{u}$ are the PLFs passing through points $(x_i,\underline{u}_i)$
and $(x_i,\bar{u}_i)$ for $i\in [N]$ given in \eqref{eq:under_upper} respectively.

\begin{table}[!ht]
\vspace{-0.3cm}
\tiny{
\centering
\captionsetup{font=small}
\caption{The efficiency of the modified polyhedral method}
\vspace{-0.2em}
\renewcommand\arraystretch{1.1}
\begin{threeparttable}
\resizebox{0.75\linewidth}{!}
{
\begin{tabular}{c|ccccc}
\hline
\textbf{Polyhedral Method} & distance & $m = 10$ & $m=30$ & $m=50$ & $m=100$ \\
\hline
\multirow{4}{*}{\makecell{Flexible breakpoints\\ flexible $N_m$}} 
& $\dd_{\rm AC}$ & 0.191 & 0.117  & 0.056 &  0.054 \\
& $\dd_{\rm R_1}$& 1.000  & 0.516 & 0.417 & 0.154 \\
& $\dd_{\rm R_2}$& 0.994 &  0.592 & 0.424 & 0.293 \\
\hline
\end{tabular}
}
\end{threeparttable}
\vspace{-0.5cm}
\label{tab-result-N}
}
\end{table}

Table~\ref{tab-result-N} reports the three metrics/distances for Example~\ref{ex:revi-4.1}
with ${\cal X}_0=\{-0.5,0,0.25,0.5\}$, $N_0=4$,
$[\underline{p},\bar{p}]=[0.05,0.95]$, and a rounding precision of $0.01$.
All three distances decrease as $M$ increases, and the small $\dd_{{\rm AC}}$ confirms the good estimation quality of $u_{N_m}^{\rm AC}$. 
While $\dd_{{\rm R_2}}$ is comparable to $\dd_{{\rm R_1}}$ for $M\leq 50$, it decreases less for $M=100$. 
To further avoid wide utility ranges, we develop the following strategy to eliminate extreme points that underpin a wide range
of utility functions (issue (b)).

\vspace{-0.1cm}   
\subsubsection{Query Design for Reducing the Maximum Utility Range}
\label{sec:query_range}
\vspace{-0.1cm}

\begin{figure}[!htbp]
\vspace{-0.2cm}
\centering
\begin{minipage}[t]{0.32\linewidth}
\centering
\includegraphics[width=2.1in]{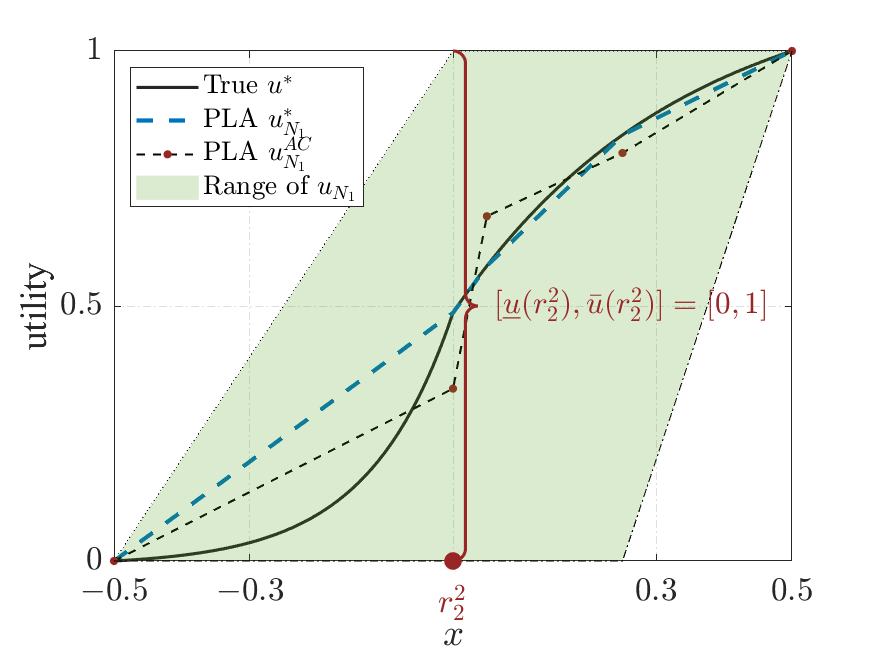}
\text{\footnotesize{\textbf{(a)}}}
\end{minipage}
\hfill
\begin{minipage}[t]{0.32\linewidth}
\centering
\includegraphics[width=2.1in]{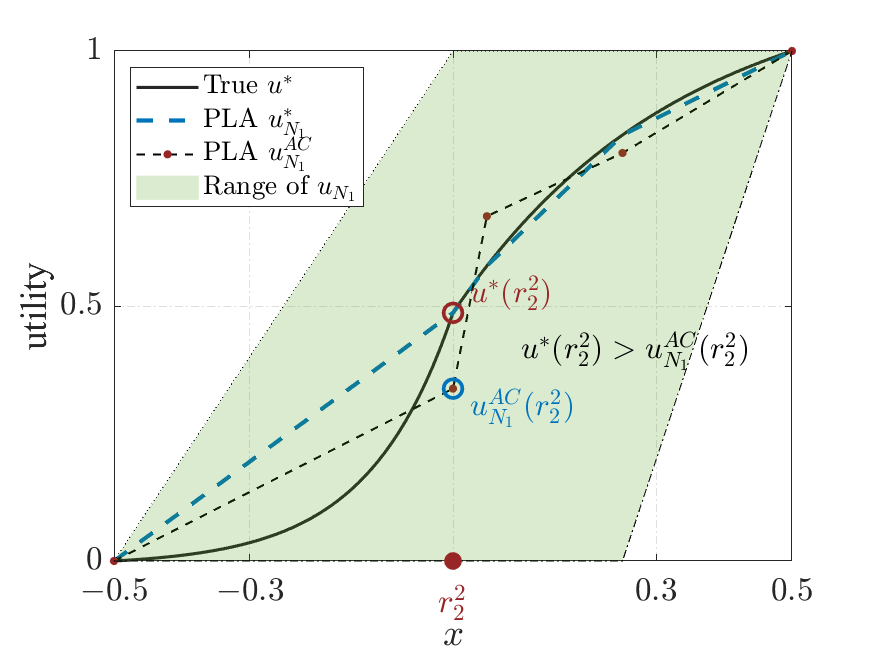}
\text{\footnotesize{\hspace{0.6cm}\textbf{(b)}}}
\end{minipage}
\hfill
\begin{minipage}[t]{0.32\linewidth}
\centering
\includegraphics[width=2.2in]{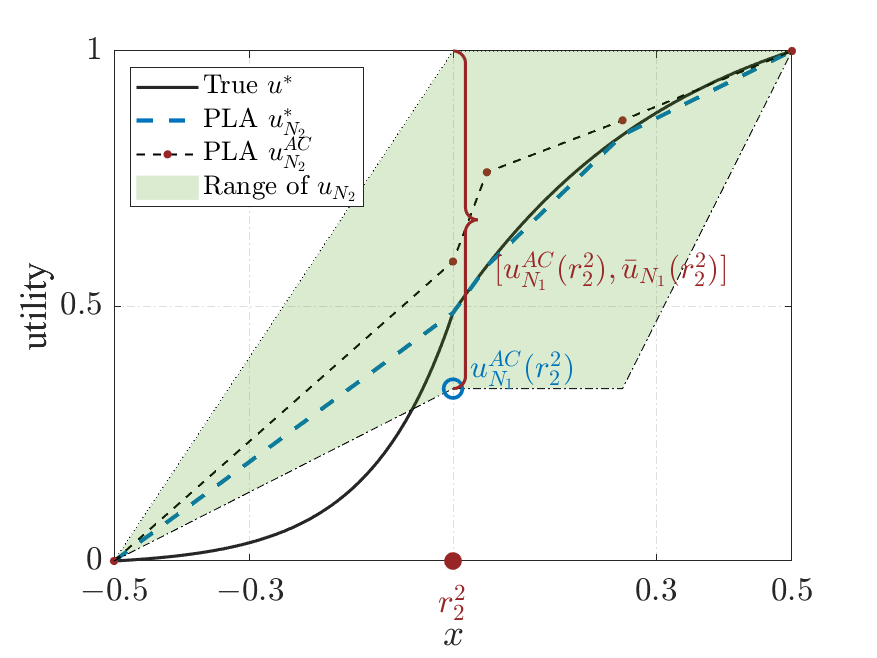}
\text{\footnotesize{\hspace{0.6cm}\textbf{(c)}}} 
\end{minipage}
\caption{Illustration of the query strategy for reducing the maximum utility range.
(a) The point $r_2^2$ is selected as a maximizer of the utility range
$[\underline{u}(r_2^2),\bar{u}(r_2^2)]$.
(b) The query probability is set to
$p^2=u_{N_1}^{\rm AC}(r_2^2)$.
(c) After incorporating the DM's response, the utility range at $r_2^2$
is reduced to
$[u_{N_1}^{\rm AC}(r_2^2),\bar{u}_{N_1}(r_2^2)]$.}
\label{fig:Fig_width}
\vspace{-0.35cm}
\end{figure}

We design an additional query to reduce the maximum utility range.
Let
$r_2^m\in\arg\max_{x\in{\cal X}_{m-1}}
\{\bar{u}(x)-\underline{u}(x)\},
$
where $\underline{u}$ and $\bar{u}$ are defined by the minimum and maximum
utility values over ${\cal U}_{N_{m-1}}^{m-1}$, respectively; see
\eqref{eq:under_upper}.
We set $r_1^m=\underline{x}$ and $r_3^m=\bar{x}$, and let
$
D=({\bm R}{\bm c}^m)^\top{\bm G}_{\bm A}(r_2^m).
$
Since
$
({\bm R}{\bm c}^m)^\top
{\bm G}_{\bm B}(r_1^m,r_3^m,p)=p,
$
we take $p^m=D$ so that the resulting cutting hyperplane satisfies
Theorem~\ref{prop:EDn=0-2}.
Thus, $p^m$ equals the value of the estimated utility
$u_{N_{m-1}}^{\rm AC}$ at $r_2^m$.
If the DM prefers $r_2^m$ to the risky lottery, then
$u(r_2^m)\geq p^m=u_{N_{m-1}}^{\rm AC}(r_2^m),$
and the utility range at $r_2^m$ is reduced from
$
[\underline{u}(r_2^m),\bar{u}(r_2^m)]$
to $[u_{N_{m-1}}^{\rm AC}(r_2^m),\bar{u}(r_2^m)]$.
Otherwise, it is reduced to
$
[\underline{u}(r_2^m),u_{N_{m-1}}^{\rm AC}(r_2^m)].
$
Since $r_2^m$ attains the maximum utility range, the query directly
reduces its current maximum width while preserving a cut through the
analytic center. See Figure~\ref{fig:Fig_width} for further illustration.

We now turn to discuss how the dimensionality of the increment vector affects the convergence behavior of the MPM,
which prompts us to investigate the complexity of the MPM (issue (c)).
The analysis in \cite[Section 3.2]{boyd2003localization} requires the
number of cuts to satisfy $M\geq N-2$, suggesting a possible role of the
dimension relative to the query number.
To apply this analysis to our framework, we consider a fixed-breakpoint variant (Algorithm~\ref{alg:RUS-fixed}) in Section~\ref{sec:variant_fixed}.
We further strengthen Criterion~\ref{eq:criterion_2}
by considering different combinations of $(r_1^m,r_3^m)$, which guarantees
that the angle between each resulting cutting hyperplane and the longest axis is bounded below by a positive constant.
Under this assumption, the analysis of
\cite{boyd2003localization}
may be adapted to derive a complexity result for Algorithm~\ref{alg:RUS-fixed}.

The above analysis motivates us to propose 
a modified query design strategy in  Section~\ref{Sec:Numberial}, where we alternately generate queries using Algorithm~\ref{alg:RUS} with the strengthened
Criterion~\ref{eq:criterion_2} (line~4 of Algorithm~\ref{alg:RUS-fixed})
and the strategy discussed in Section~\ref{sec:query_range}. 
This will enable us to exploit the advantages of (a) adaptively adding new breakpoints, (b) reducing the utility range 
at its widest point, and (c) controlling the rate of increase in the number of breakpoints.
Note that the preceding examples use different utility functions to highlight
specific features of the proposed elicitation method. In the following
PRO analysis, we conduct the preference elicitation for three utility
functions in a unified setting (see Figure~\ref{fig:feasible-se-3}, \ref{fig:feasible-91_range}, and \ref{fig:feasible-91_range3}) and evaluate the resulting PRO models.

\vspace{-0.2cm}
\section{Application in Preference Robust Optimization}
\label{Sec:Numberial}
\vspace{-0.2cm}

In this section, 
we apply the MPM developed
in the preceding section to decision-making problems where the DM's aim is to maximize the expected utility of rewards 
and where information about
the true utility function which characterizes the DM's preference is incomplete.
We consider the following 
one-stage maximin utility preference robust optimization model
\bgeqn
\label{eq:PRO-PLA-N-dis-2-mm}
\max_{{\bm z}\in 
Z }\min_{u_{N}\in {\cal U}_{N}^M} \bbe[u_{N}(f({\bm z},{\bm \xi}))],
\edeqn
where $Z\subset \R^n$ is a compact and convex set,
${\bm z}\in Z$ is a decision vector and
${\bm \xi}$ is a random vector with support $\Xi\subset \R^s$,
$f:Z\times \Xi \to \R$ is a continuous function representing reward,
${\cal U}_N^M$ is defined as in \eqref{eq:PLA_flexible}.
In this section, we write $N$ for $N_M$ to ease the exposition.
In practice, the
decision based on the worst-case expected utility 
may be overly
conservative,
particularly when the ambiguity set ${\cal U}_N^M$ is large.
To control the degree of conservatism, we introduce a mechanism for reducing the size of the ambiguity set ${\cal U}_N^M$.
Specifically, we 
rewrite
the inner minimization problem of \eqref{eq:PRO-PLA-N-dis-2-mm} as
\begin{subequations}
\label{eq:minmin_bounds-a0}
 \begin{align}
\displaystyle  \min_{{\bm v}_{[N-1]} \in \widehat{\cal V}_N^M} & \bbe[{\bm v}_{[N-1]}^\top{\bm g}(f({\bm z},{\bm \xi}))]\\
{\rm s.t.}\quad \; & {v}_i-{c}_i^{M}\in [-\underline{v}_i,\bar{v}_i], ~i\in [N-1] , 
 \label{eq:minmin_bounds-a-0b}
\end{align}
\end{subequations}
where $\widehat{\cal V}_N^M:=\left\{{\bm v}_{[N-1]}\in \R^{N-1}_+: {\bm v}_{[N-1]}={\bm R}{\bm v},{\bm v}\in {\cal V}_{N}^M\right\}$ and 
\begin{subequations}
\label{eq:bounds_vo-0}
\begin{align}
&\underline{v}_i:=-\min
\left\{{v}_i-{c}_{i}^{M}:
{\bm v}_{[N-1]}\in  \widehat{\cal V}_N^M \right\},\;i\in [N-1], \\
&\bar{v}_i:=\max
\left\{{v}_i-{c}_{i}^{M}: 
{\bm v}_{[N-1]}\in  \widehat{\cal V}_N^M \right\},\;i\in [N-1]. 
\end{align}
\end{subequations}
Here, ${\cal V}_N^M$ is defined in \eqref{eq:choice-xkyk-0}, and
$\widehat{\cal V}_N^M$ represents the same feasible set using the full
$(N-1)$-dimensional increment vector ${\bm v}_{[N-1]}={\bm R}{\bm v}$.
This representation allows the expected utility and  constraints in \eqref{eq:choice-xkyk-0}
to be expressed directly as linear functions of ${\bm v}_{[N-1]}$,
facilitating the subsequent dual reformulations in
Proposition~\ref{prop:reformulate_1}  and Theorem~\ref{thm:reformulate_n}. Moreover,
$\bm c_{[N-1]}^M=(\bm c^M,1-\bm e^\top\bm c^M)\in
\mathbb R^{N-1}$ denotes the full increment vector induced by the
analytic center $\bm c^M$ in the reduced space, where $c_i^M$ denotes
the $i$-th component of $\bm c_{[N-1]}^M$.
The quantities $\underline v_i$ and $\bar v_i$ 
signify
the range
of the $i$-th increment relative to $c_i^M$ over
${\mathcal V}_N^M$. Specifically,
$c_i^M-\underline v_i$ and $c_i^M+\bar v_i$ are the minimum and
maximum feasible values of $v_i$, respectively.
Constraints
\eqref{eq:minmin_bounds-a-0b}
are 
redundant for the inner minimization problem because the bounds $[-\underline v_i,\bar v_i]$ are derived from
$\widehat{\mathcal V}_N^M$. 
We next tighten these bounds to control the degree of conservatism. To this end, 
we introduce auxiliary variables
$y_i\in[0,1]$, $i\in[N-1]$, which scale the corresponding intervals
in \eqref{eq:minmin_bounds-a-0b}.
This leads to 
the following 
modified inner minimization problem: 
\begin{subequations}
\label{eq:min_expand-20}
\begin{align}
\inmat{(CR-A)}\quad \quad \displaystyle  \min_{{\bm v}_{[N-1]} \in \widehat{\cal V}_N^M ,{\bm y}}&
\bbe[{\bm v}_{[N-1]}^\top{\bm g}(f({\bm z},{\bm \xi}))] \\
{\rm s.t.}\;\;\quad &  {v}_i-{c}_i^{M}\in [-y_i\underline{v}_i,y_i\bar{v}_i], \;i\in [N-1], \label{eq:min_expand-20-a}\\
&0\leq y_i\leq 1, \; {\bm e}^\top{\bm y}=\Gamma, \label{eq:min_expand-20-b}
\end{align}
\end{subequations}
where  $\Gamma\leq N-1$ is a prescribed 
parameter and ``CR'' stands for conservatism reduction.
When $\Gamma= N-1$, \eqref{eq:min_expand-20} is equivalent to \eqref{eq:minmin_bounds-a0}. 
Here we are interested in the case that 
$\Gamma< N-1$, 
in which some of the intervals in \eqref{eq:min_expand-20-a}
are tightened, thereby potentially increasing the optimal value
relative to \eqref{eq:minmin_bounds-a0}.
This type of adjustment is commonly referred to as conservatism
reduction in robust optimization;
see e.g.~\cite{BeS04}.
Influenced
by the structure of the conservatism reduction scheme in \eqref{eq:min_expand-20}, we consider the following \underline{modification}:
\begin{subequations}
\label{eq:robust-pi0}
 \begin{align}
 \inmat{(CR-B)}\quad \quad \displaystyle  \min_{{\bm v}_{[N-1]} \in \widehat{\cal V}_N^M ,{\bm \Pi}\in {\cal S}}&
 \bbe[{\bm v}_{[N-1]}^\top{\bm g}(f({\bm z},{\bm \xi}))] \\
{\rm s.t.}\quad\quad &  {v}_i-{c}_i^{M}\in [-({\bm \Pi}_i \sigma) \underline{v}_i,({\bm \Pi}_i \sigma)\bar{v}_i], \; i\in [N-1], \label{eq:robust-pi0-b}
\end{align}
\end{subequations}
where 
${\cal S}:=\{{\bm \Pi}\in \R^{(N-1)\times(N-1)}:0\leq {\bm \Pi}_{i,j}\leq 1, \forall i,j\in [N-1],\;{\bm \Pi}{\bm e}={\bm e},
{\bm \Pi}^\top{\bm e}={\bm e}\}$ is 
the set of doubly stochastic matrices,
and
${\bm \Pi}_i$ is the $i$-th row of matrix ${\bm \Pi}$ for $i\in [N-1]$.
The vector $\bm\sigma$ contains weights
generated by the function
$\varsigma(t):=\gamma t^{\gamma-1}$, parameterized by $\gamma\in (0,1]$,
as follows:
\bgeqn
\label{eq:sigma0}
{\bm \sigma}=(\sigma_1,\cdots,\sigma_{N-1}):=\left(\varsigma\left({1}/{N}\right),\cdots,\varsigma\left({N-1}/{N}\right)\right) \slash\varsigma\left({1}/{N}\right).
\edeqn
The key difference
between \eqref{eq:min_expand-20} and \eqref{eq:robust-pi0}
can be summarized as follows:
\begin{itemize}
\item
The optimal $y_i$ in \eqref{eq:min_expand-20} takes either $0$ or $1$, as shown
below. When $y_i=1$, the original interval
$[-\underline v_i,\bar v_i]$ is retained, whereas when $y_i=0$, the
interval shrinks to $\{0\}$, fixing $v_i$ at $c_i^M$. Thus, the reduced feasible region is
$\widehat{\mathcal V}_N^M\cap \left\{
{\bm v}_{[N-1]}:
v_i-c_i^M
\in
[-y_i\underline v_i, y_i\bar v_i],
\ i\in[N-1]
\right\}$.
Although $u_N^*=({\bm R}{\bm v}^*)^\top {\bm g}\in {\cal U}_N^M$,
${\bm v}^*$ need not remain in this reduced region,
since ${\bm c}^M$ may not coincide with ${\bm v}^*$.

\item
In \eqref{eq:robust-pi0}, the reduced feasible region is
$
\widehat{\mathcal V}_N^M
\cap
\left\{
{\bm v}_{[N-1]}:
v_i-c_i^M
\in
[-(\bm\Pi_i{\bm\sigma})\underline v_i,
  (\bm\Pi_i{\bm\sigma})\bar v_i],
\ i\in[N-1]
\right\}.
$
Since $\bm\Pi_i{\bm\sigma}\in[0,1]$, the intervals can be reduced by
different proportions between these two extremes. This provides
a more gradual tightening and mitigates the risk of excluding ${\bm v}^*$
from the reduced set, although its inclusion is not
guaranteed. Details of this scheme are provided in Section~\ref{sec:CR_B}.
\end{itemize}

In the remainder of this section, we elaborate
two
conservatism reduction schemes   \eqref{eq:min_expand-20} and~\eqref{eq:robust-pi0}.

\vspace{-0.1cm}
\subsection{Conservatism Reduction Scheme in \eqref{eq:min_expand-20}}
\vspace{-0.1cm}

The conservatism  reduction scheme in \eqref{eq:min_expand-20}
can be explained through
the Lagrange dual problem of \eqref{eq:minmin_bounds-a0},
which is given in the following proposition.

\begin{proposition}
[Lagrange Dual  of \eqref{eq:minmin_bounds-a0}]
\label{prop:reformulate_1} 
Let $f({\bm z},{\bm \xi})={\bm z}^\top {\bm \xi}$, where 
${\bm \xi}$ follows a discrete distribution with $P({\bm \xi}={\bm \xi}^k)=\frac{1}{K}$ for $k\in [K]$. Then the
Lagrange dual  of \eqref{eq:minmin_bounds-a0} is
\begin{subequations}
\label{eq:inner_problem}
\begin{align}
\displaystyle \max_{\underline{\bm \lambda},\bar{\bm \lambda},\beta,{\bm \eta},{\bm \theta}} & ~({\bm c}_{[N-1]}^{M})^\top\left(\frac{1}{K}\sum_{k=1}^K{\bm g} ({\bm z}^\top {\bm \xi}^k)\right)
-\underline{\bm \lambda}^\top\underline{\bm v}
-\bar{\bm \lambda}^\top\bar{\bm v}
-{\bm \eta}^\top {\bm c}_{[N-1]}^{M} \nonumber \\
&
+\sum_{m=1}^M {\theta}_m h_m ({\bm c}_{[N-1]}^{M})^{\top}
\left(G_{\bm A} (r_{1}^{m},r_3^{m},p^{m})-G_{\bm B} (r_2^{m,s})\right)  \label{eq:inner_obj}\\
{\rm s.t.}\quad & 
\frac{1}{K}\sum_{k=1}^K{\bm g} ({\bm z}^\top {\bm \xi}^k)
-\underline{\bm \lambda}
+\bar{\bm \lambda} +\beta{\bm e}-{\bm \eta}+\sum_{m=1}^M {\theta}_m h_m \left(G_{\bm A} (r_{1}^{m},r_3^{m},p^{m})-G_{\bm B} (r_2^{m})\right)=0,
\label{eq:inner_cons}\\
&  \underline{\bm \lambda}\geq 0,\bar{\bm \lambda}\geq 0,{\bm \eta}\geq 0,{\bm \theta}\geq 0.
\end{align}
\end{subequations}
\end{proposition}

In the dual program \eqref{eq:inner_problem},
the bounds $\underline{\bm v}$ and
$\bar{\bm v}$ appear only in the objective through the term
$-\underline{\bm \lambda}^\top\underline{\bm v}
-\bar{\bm \lambda}^\top\bar{\bm v}=\sum_{i=1}^{N-1}(-\underline{\lambda}_i\underline{v}_i-\bar{\lambda}_i\bar{v}_i)$.
For any $i\in [N-1]$,
if $-\underline{\lambda}_i\underline{v}_i-\bar{\lambda}_i\bar{v}_i<0$,
then $\underline{\lambda}_i>0$ or $\bar{\lambda}_i>0$.
By complementary slackness, the former implies that the $i$-th lower-bound
constraint
$\widetilde{v}_i=v_i-{c}_i^{M}\geq -\underline{v}_i$ in \eqref{eq:minmin_bounds-a-0b}
is binding,
whereas the latter implies that the $i$-th upper-bound constraint
$\widetilde{v}_i=v_i-{c}_i^{M}\leq \bar{v}_i$ in \eqref{eq:minmin_bounds-a-0b} is 
binding.
Motivated by this observation, we replace the term
$-\underline{\bm \lambda}^\top\underline{\bm v}
-\bar{\bm \lambda}^\top\bar{\bm v}$ 
with 
\bgeqn
\label{eq:product_sort-2}
I_1(\underline{\bm \lambda},\bar{\bm \lambda},\Gamma)
:= \min_{0\leq y_i\leq 1, {\bm e}^\top{\bm y}=\Gamma} \sum_{i=1}^{N-1}y_i\Big(-\Big(\underline{\lambda}_i\underline{v}_i
+\bar{\lambda}_i \bar{v}_i\Big)\Big),
\edeqn
and consider the following problem:
\bgeqn
\label{eq:level_con_Gamma}
\displaystyle \max_{\underline{\bm \lambda}\geq 0,\bar{\bm \lambda}\geq 0,{\bm \eta}\geq 0,{\bm \theta}\geq 0,\beta} 
\left\{I_2(\underline{\bm \lambda},\bar{\bm \lambda},\beta,{\bm \eta},{\bm \theta})
+
I_1(\underline{\bm \lambda},\bar{\bm \lambda},\Gamma)
:\inmat{constraint } \eqref{eq:inner_cons}\right\}, 
\edeqn
where 
$I_2(\underline{\bm \lambda},\bar{\bm \lambda},\beta,{\bm \eta},{\bm \theta})$ denotes the objective function of problem \eqref{eq:inner_problem} excluding $-\underline{\bm \lambda}^\top\underline{\bm v}-\bar{\bm \lambda}^\top\bar{\bm v}$.
When $\Gamma$ is an integer, problem
\eqref{eq:product_sort-2} 
selects the $\Gamma$ most
negative terms among
$-\Big(\underline{\lambda}_i\underline{v}_i
+\bar{\lambda}_i \bar{v}_i\Big)$.
The next theorem states the  equivalence
of
\eqref{eq:min_expand-20} and \eqref{eq:level_con_Gamma}
through the Lagrange dual 
of \eqref{eq:min_expand-20} and \eqref{eq:product_sort-2}.

\begin{theorem}
[Reformulation of 
Problems
\eqref{eq:min_expand-20} and \eqref{eq:level_con_Gamma}]
\label{thm:reformulate_1} 
Consider the case that 
$f({\bm z},{\bm \xi})={\bm z}^\top {\bm \xi}$ and 
${\bm \xi}$ follows a discrete distribution with $P({\bm \xi}={\bm \xi}^k)=\frac{1}{K}$ for $k\in [K]$.  Then
problem \eqref{eq:min_expand-20} and \eqref{eq:level_con_Gamma}  
can be reformulated as:
\begin{subequations}
\label{eq:inner_problem-robust}
\begin{align}
\displaystyle \max_{\underline{\bm \lambda},\bar{\bm \lambda},\beta,{\bm \eta},{\bm \theta}, {\bm \tau}, \tau_0} & ~({\bm c}_{[N-1]}^{M})^\top\left(\frac{1}{K}\sum_{k=1}^K{\bm g} ({\bm z}^\top {\bm \xi}^k)\right)
-{\bm \eta}^\top {\bm c}_{[N-1]}^{M}-\sum_{i=1}^{N-1}\tau_i-\tau_0\Gamma\\
&
+\sum_{m=1}^M {\theta}_m h_m ({\bm c}_{[N-1]}^{M})^{\top}
 \left(G_{\bm A} (r_{1}^{m},r_3^{m},p^{m})-G_{\bm B} (r_2^{m,s})\right) \\
{\rm s.t.}\qquad & 
\frac{1}{K}\sum_{k=1}^K{\bm g} ({\bm z}^\top {\bm \xi}^k)
-\underline{\bm \lambda}
+\bar{\bm \lambda} +\beta{\bm e}-{\bm \eta}+\sum_{m=1}^M {\theta}_m h_m \left(G_{\bm A} (r_{1}^{m},r_3^{m},p^{m})-G_{\bm B} (r_2^{m})\right)=0,\\
&- \underline{\lambda}_i \underline{v}_i-\bar{\lambda}_i \bar{v}_i+\tau_i+\tau_0 \geq 0, i\in [N-1],\\
&\underline{\bm \lambda}\geq 0,\bar{\bm \lambda}\geq 0,{\bm \eta}\geq 0,{\bm \theta}\geq 0, {\bm \tau}\geq 0.
\end{align}
\end{subequations}
\end{theorem}

\noindent{\bf Proof}.
Problem \eqref{eq:inner_problem-robust} is the Lagrange dual of \eqref{eq:min_expand-20}. 
Alternatively, by deriving
the Lagrange
dual of \eqref{eq:product_sort-2} and substituting 
the resulting dual problem into
\eqref{eq:level_con_Gamma}, 
we also obtain \eqref{eq:inner_problem-robust}.
The details are omitted.
\hfill $\Box$

The equivalence shows that tightening the bound constraints on
$\bm v_{[N-1]}$ in the primal problem \eqref{eq:minmin_bounds-a0} has the same effect as assigning
weights to the  
bound terms in the dual problem \eqref{eq:inner_problem}.

\vspace{-0.1cm}
\subsection{Conservatism Reduction Scheme in \eqref{eq:robust-pi0}}
\label{sec:CR_B}
\vspace{-0.1cm}

Instead of selecting the $\Gamma$ most negative terms in \eqref{eq:product_sort-2}, 
we assign a weight in $[0,1]$ to each term $-\left(\underline{\lambda}_{i} \underline{v}_{i}+\bar{\lambda}_{i}\bar{v}_{i}\right)$ through $\bm\Pi\bm\sigma$.
The weight vector is defined in \eqref{eq:sigma0},
where $\varsigma(t)=\gamma t^{\gamma-1}$ with $\gamma\in (0,1]$.
For $\gamma\in(0,1)$, the function $t^\gamma$ is strictly concave
and increasing, and hence its derivative $\varsigma(t)$ 
is strictly
decreasing. Therefore,
$
\sigma_1> \sigma_2> \cdots >\sigma_{N-1}.$
As $\gamma$ approaches $1$, $t^\gamma$ becomes less concave and
approaches the linear function $t$, and the differences between
the
weights diminish. When $\gamma=1$, $\varsigma(t)=1$ and 
thus
$\sigma_1=\cdots=\sigma_{N-1}=1$.

If $\gamma=0.1$, $N=10$,
then 
${\bm \sigma}=(1,0.535,0.371,0.287,0.234,0.198,
0.172,0.153,0.137)$.
If  $\gamma=0.5$, $N=10$,
then ${\bm \sigma}=(1,0.707,0.576,0.499, 0.447,0.407,0.377,0.353,0.333)$.
For additional examples of ${\bm \sigma}$,
we observe that, as $\gamma$ increases, the differences between
consecutive weights decrease, which is consistent with the discussion above.
By
analogy with
\eqref{eq:level_con_Gamma}, we consider
the following problem:
\bgeqn
\max_{\underline{\bm \lambda}\geq 0,\bar{\bm \lambda}\geq 0,{\bm \eta}\geq 0,{\bm \theta}\geq0,\beta}  \left\{I_2(\underline{\bm \lambda},\bar{\bm \lambda},\beta,{\bm \eta},{\bm \theta})
+\displaystyle \min_{{\bm \Pi}\in {\cal S}} \left\{\sum_{i\in [N-1]} 
({\bm \Pi}_i{\bm {\bm \sigma}})
\cdot \Big(-\Big(\underline{\lambda}_i \underline{v}_i
+\bar{\lambda}_i\bar{v}_i \Big)\Big) \right\}
:\inmat{constraint } \eqref{eq:inner_cons}\right\}.\quad 
\label{eq:defi_new_control-a}
\edeqn
Note that
\bgeqn
\label{eq:sort}
 \min_{{\bm \Pi}\in {\cal S}} \left\{\sum_{i\in [N-1]} ({\bm \Pi}_i{\bm \sigma})\cdot \Big(-\Big( \underline{\lambda}_i\underline{v}_i
+\bar{\lambda}_i\bar{v}_i \Big)\Big)\right\}=\sum_{i=1}^{N-1}\sigma_i \left(-\left(\underline{\lambda}_{(i)} \underline{v}_{(i)}+\bar{\lambda}_{(i)}\bar{v}_{(i)}\right)\right),
\edeqn
where the parenthesized indices $(1),\ldots,(N-1)$ denote a permutation
of $1,\ldots,N-1$ such that
$
-\left(
\underline{\lambda}_{(1)}\underline v_{(1)}
+\bar{\lambda}_{(1)}\bar v_{(1)}
\right)
\leq \cdots \leq
-\left(
\underline{\lambda}_{(N-1)}\underline v_{(N-1)}
+\bar{\lambda}_{(N-1)}\bar v_{(N-1)}
\right).
$
With \eqref{eq:sort}, problem \eqref{eq:defi_new_control-a} becomes
\bgeqn
\label{eq:level_con_sigma}
\max_{\underline{\bm \lambda}\geq 0,\bar{\bm \lambda}\geq 0,{\bm \eta}\geq 0,{\bm \theta}\geq 0,\beta}  \left\{I_2(\underline{\bm \lambda},\bar{\bm \lambda},\beta,{\bm \eta},{\bm \theta})
+\displaystyle \sum_{i\in [N-1]} {\sigma}_i \Big(-\Big(\underline{\lambda}_{(i)}\underline{v}_{(i)}
+\bar{\lambda}_{(i)}\bar{v}_{(i)} \Big)\Big) 
:\inmat{constraint } \eqref{eq:inner_cons}\right\}.
\edeqn
Thus, 
larger
weights
are assigned to more negative 
terms
$-\left(\underline{\lambda}_{i} \underline{v}_{i}+\bar{\lambda}_{i}\bar{v}_{i}\right)$.
Similarly, we establish the equivalence between \eqref{eq:robust-pi0} and \eqref{eq:defi_new_control-a}
through the Lagrange dual representations of \eqref{eq:robust-pi0} and \eqref{eq:sort}.

\begin{theorem}
[Reformulations of \eqref{eq:robust-pi0}
and \eqref{eq:defi_new_control-a}]
\label{thm:reformulate_n}
Assume the setting and conditions of 
Theorem~\ref{thm:reformulate_1}.
Then problems \eqref{eq:robust-pi0}
and \eqref{eq:defi_new_control-a}
can be reformulated as 
\begin{subequations}
\label{eq:reformulation_conser}
\begin{align}
\displaystyle \max_{\underline{\bm \lambda},\bar{\bm \lambda},\beta,{\bm \eta},{\bm \theta},{\bm \Theta},
\underline{\bm \tau},\bar{\bm \tau}} & ({\bm c}_{[N-1]}^{M})^\top\left(\frac{1}{K}\sum_{k=1}^K{\bm g} ({\bm z}^\top {\bm \xi}^k)\right) +\sum_{m=1}^M {\theta}_m h_m  ({\bm c}_{[N-1]}^{M})^{\top}
\left(G_{\bm A} (r_{1}^{m},r_3^{m},p^{m})-G_{\bm B} (r_2^{m})\right)\nonumber\\
&-{\bm \eta}^\top {\bm c}_{[N-1]}^{M}-\sum_{i,j} {\bm \Theta}_{ij}-\underline{\bm \tau}^\top {\bm e}-\bar{\bm \tau}^\top {\bm e}\\
 {\rm s.t.} \qquad  
& \frac{1}{K}\sum_{k=1}^K{\bm g} ({\bm z}^\top {\bm \xi}^k)
-\underline{\bm \lambda}
+\bar{\bm \lambda} +\beta{\bm e}-{\bm \eta}+\sum_{m=1}^M {\theta}_m h_m \left(G_{\bm A} (r_{1}^{m},r_3^{m},p^{m})-G_{\bm B} (r_2^{m})\right)=0,\label{eq:reformulation_conser-b}\\
& (-\underline{\bm \lambda} \circ \underline{\bm v}-\bar{\bm \lambda} \circ \bar{\bm v}){\bm \sigma}^\top+\underline{\bm \tau}
{\bm e}^\top+{\bm e} \bar{\bm \tau}^\top+{\bm \Theta}\geq 0,\label{eq:reformulation_conser-c}\\
& \underline{\bm \lambda}\geq 0,\bar{\bm \lambda}\geq 0, {\bm \eta} \geq 0,{\bm \theta}\geq 0,
{\bm \Theta} \geq 0\label{eq:reformulation_conser-d}.
\end{align}
\end{subequations}
\end{theorem}

Combining the results above,
the PRO model \eqref{eq:PRO-PLA-N-dis-2-mm} under
the conservatism reduction 
scheme (CR-B)
can be formulated as
\bgeqn
\label{eq:PRO-PLA-N-dis-2}
\displaystyle  \max_{{\bm z}\in Z}\min_{{\bm v}_{[N-1]} \in \widehat{\cal V}_N^M ,{\bm \Pi}\in {\cal S}} \Big\{
\bbe[{\bm v}_{[N-1]}^\top{\bm g}(f({\bm z},{\bm \xi}))] :  {v}_i-{c}_i^{M}\in [-({\bm \Pi}_i \sigma) \underline{v}_i,({\bm \Pi}_i \sigma)\bar{v}_i], ~i\in [N-1]\Big\}.
\edeqn
By Theorem~\ref{thm:reformulate_n},
problem \eqref{eq:PRO-PLA-N-dis-2} can be reformulated as the following
single maximization problem:
\begin{align}
\label{eq:reformulation_conser-a}
 \displaystyle \max_{
 {\bm z}\in Z,
 \underline{\bm \lambda},\bar{\bm \lambda},\beta,{\bm \eta},{\bm \theta},{\bm \Theta},
 \underline{\bm \tau},\bar{\bm \tau}} & ({\bm c}_{[N-1]}^{M})^\top\left(\frac{1}{K}\sum_{k=1}^K{\bm g} ({\bm z}^\top {\bm \xi}^k)\right) +\sum_{m=1}^M {\theta}_m h_m  ({\bm c}_{[N-1]}^{M})^{\top}
 \left(G_{\bm A} (r_{1}^{m},r_3^{m},p^{m})-G_{\bm B} (r_2^{m})\right)\nonumber\\
 &-{\bm \eta}^\top {\bm c}_{[N-1]}^{M}-\sum_{i,j} {\bm \Theta}_{ij}-\underline{\bm \tau}^\top {\bm e}-\bar{\bm \tau}^\top {\bm e}\\
 {\rm s.t.} \quad\quad \,\,\,
& \inmat{constraints }\eqref{eq:reformulation_conser-b}, \eqref{eq:reformulation_conser-c},\eqref{eq:reformulation_conser-d},
\nonumber
\end{align}
where ${\bm g}(\cdot)$ is defined as in \eqref {eq:g(x)-incre}. 
Note that ${\bm g}(x)$ can be represented by a mixed-integer  linear system;
see Proposition~\ref{prop:g(x)-rfmlt}. Consequently, 
 \eqref{eq:reformulation_conser-a}  admits the following MILP reformulation:
\begin{subequations}
\label{eq:reformulation_conser-a1}
\begin{align}
 \displaystyle \max_{\substack{
 {\bm z}\in Z,
 \underline{\bm \lambda},\bar{\bm \lambda},\beta,{\bm \eta},{\bm \theta}\\
 {\bm \Theta},
 \underline{\bm \tau},\bar{\bm \tau},{\bm y}^k,{\bm w}^k}} & ({\bm c}_{[N-1]}^{M})^\top\left(\frac{1}{K}\sum_{k=1}^K {\bm P}{\bm y}^k+{\bm Q}{\bm w}^k\right) +\sum_{m=1}^M {\theta}_m h_m  ({\bm c}_{[N-1]}^{M})^{\top}
 \left(G_{\bm A} (r_{1}^{m},r_3^{m},p^{m})-G_{\bm B} (r_2^{m})\right)\nonumber\\
 &-{\bm \eta}^\top {\bm c}_{[N-1]}^{M}-\sum_{i,j} {\bm \Theta}_{ij}-\underline{\bm \tau}^\top {\bm e}-\bar{\bm \tau}^\top {\bm e}\\
 {\rm s.t.}  \qquad 
 &\frac{1}{K}\sum_{k=1}^K({\bm P}{\bm y}^k+{\bm Q}{\bm w}^k)
-\underline{\bm \lambda}
+\bar{\bm \lambda} +\beta{\bm e}-{\bm \eta}+\sum_{m=1}^M {\theta}_m h_m \left(G_{\bm A} (r_{1}^{m},r_3^{m},p^{m})-G_{\bm B} (r_2^{m})\right)=0,\\
 &  x_{j}w_{j}^{k}-y_{j}^{k}\leq 0, \;
    x_{j+1}w_{j}^{k}-y_{j}^{k} \geq 0, \;
    w_{j}^{k}\in \{0,1\},\;j\in [N-1],\;k\in [K],~~\\
& \sum_{j=1}^{N-1}w_{j}^{k}=1,\;\sum_{j=1}^{N}y_{j}^k=
{\bm z}^\top {\bm \xi}^k,
\;k\in [K], \\
& \inmat{constraints }
\eqref{eq:reformulation_conser-c},\eqref{eq:reformulation_conser-d},
\end{align}
\end{subequations}
where the matrices
${\bm P},{\bm Q}\in \R^{(N-1)\times (N-1)}$ are defined as in \eqref{eq:definition-AB}.

\vspace{-0.1cm}
\subsection{Numerical Tests}
\label{sec:tests}
\vspace{-0.1cm}

In this section,
we consider a specific portfolio optimization problem,
in which
an investor allocates a budget between
$n$ assets. Let $z_i$ denote
the proportion of the total budget allocated to asset $i$ for $i\in [n]$,
and let
$Z:=\{{\bm z}\in \R^n:{\bm z}\geq 0,{\bm e}^{\top}{\bm z}=1\}$
denote the feasible set.
Let $\xi_i$ be the random weekly return of asset $i$ for $i\in [n]$. 
The total return of the portfolio is then
$f({\bm z},{\bm \xi})={\bm z}^{\top}{\bm \xi}$.
We consider the formulation \eqref{eq:PRO-PLA-N-dis-2},
where the ambiguity set ${\cal U}_{N}^M$ is constructed by the MPM.
Using historical stock-return data, we investigate the performance of
the proposed PRO model when no prior
information about the DM's utility function is available.

\vspace{-0.1cm}
\subsubsection{Setup}
\label{sec:setup}
\vspace{-0.1cm}

We consider 10 financial indices:
S\&P 500 Index (GSPC),
CBOE Volatility Index (VIX),
Dow Jones Industrial Average Index (DJI),
Dow Jones Transportation Average Index (DJT),
Dow Jones Composite Average Inex (DJA),
Treasury Yield 5 Years Index (FVX),
Treasury Yield 10 Years Index (TNX),
Treasury Yield 30 Years Index (TYX),
NASDAQ Composite Index (IXIC),
Russell 2000 Index (RUT).
We use 24 monthly observations from January 2015 to December 2016\footnote{https://finance.yahoo.com},
so that $K=24$.
For evaluation purposes, we specify the true utility function as
$u^*(t)=(2*(e^{8 t}-1))/8$
for $t\in [-0.5,0]$ and
$u^*(t)=(1-e^{-3 t})/3$
for $t\in [0,0.5]$,
while treating the investor's utility function as unknown during 
the process of eliciting preference.
We investigate
how the optimal value
and the
worst-case utility function
evolve as more information about the investor's preferences
is elicited under different levels of conservatism.
All tests are conducted in MATLAB R2024b on a PC with 16 GB of RAM
and a 2.20 GHz Intel Core i7 processor. The optimization problems are
solved using Gurobi and YALMIP.

\vspace{-0.1cm}
\subsubsection{In-Sample Performance of (CR-B) under the S-Shaped Utility Function}
\label{sec:insample}
\vspace{-0.1cm}

\begin{figure}[!ht]
\vspace{-0.2cm}
\begin{minipage}[t]{0.32\linewidth}
\centering
\includegraphics[width=2.2in]{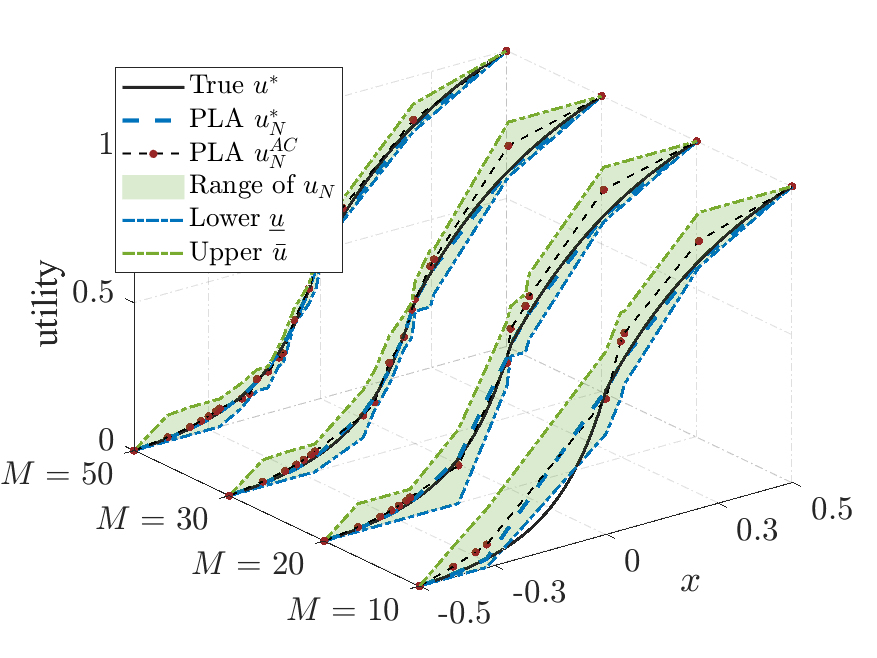}
\text{\footnotesize{\textbf{(A)}}}
\end{minipage}
\begin{minipage}[t]{0.32\linewidth}
\centering
\includegraphics[width=2.1in]{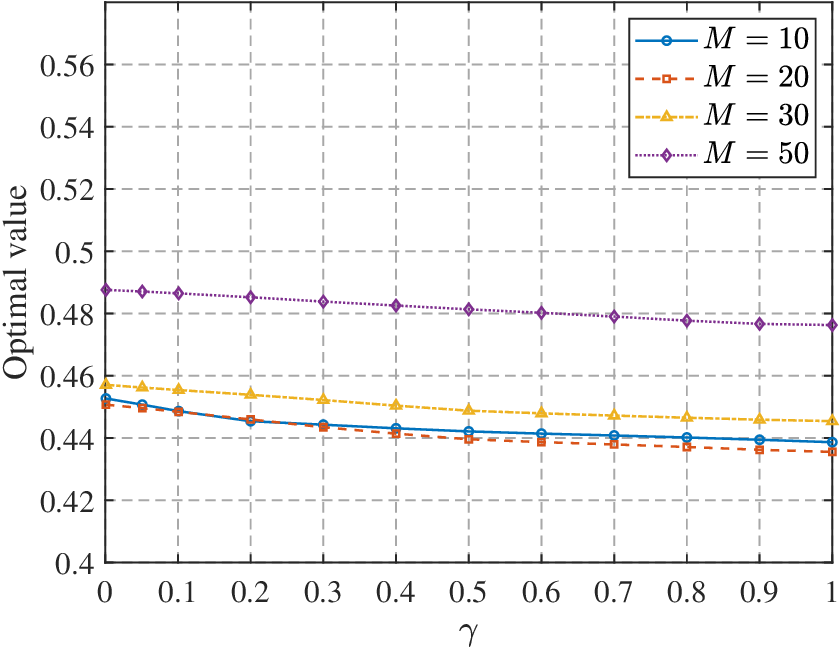}
\text{\footnotesize{\hspace{0.8cm}\textbf{(B1)}}}
\end{minipage}
\hfill
\begin{minipage}[t]{0.32\linewidth}
\centering
\includegraphics[width=2.1in]{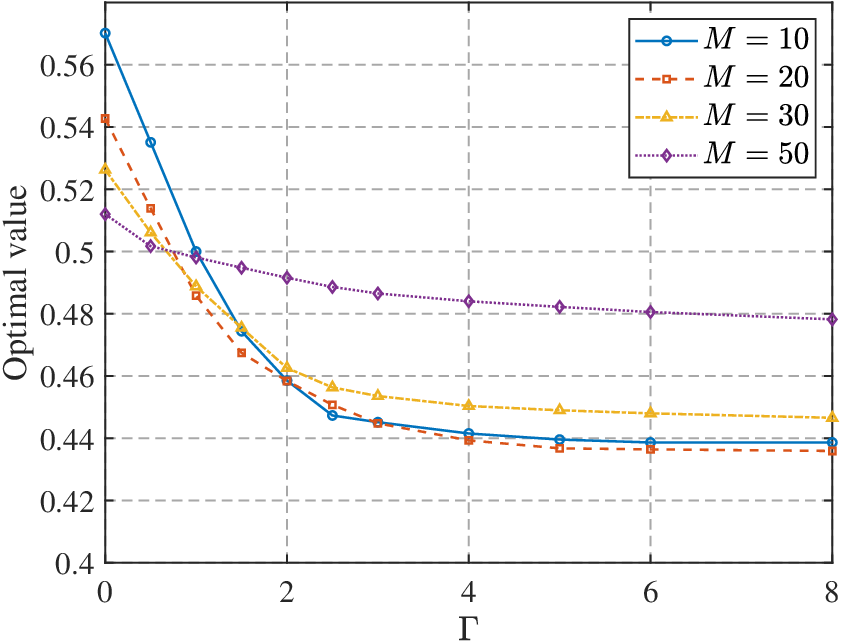}
  \text{\footnotesize{\hspace{0.8cm}\textbf{(B2)}}}
\end{minipage}
\caption{\footnotesize{
(A) Utility function range as $M$ increases.
(B1) Optimal values of $\rho_\gamma^B$as $\gamma$ varies.
(B2) Optimal values of $\rho_\Gamma^A$ as $\Gamma$ varies.
Both panels report results for the number of queries $M=10,20,30,50$.
}}
\label{fig:feasible-se-3}
\end{figure}

\begin{figure}[!h]
\vspace{-1em}
\begin{minipage}[t]{0.35\linewidth}
\centering
\includegraphics[width=2.6in]{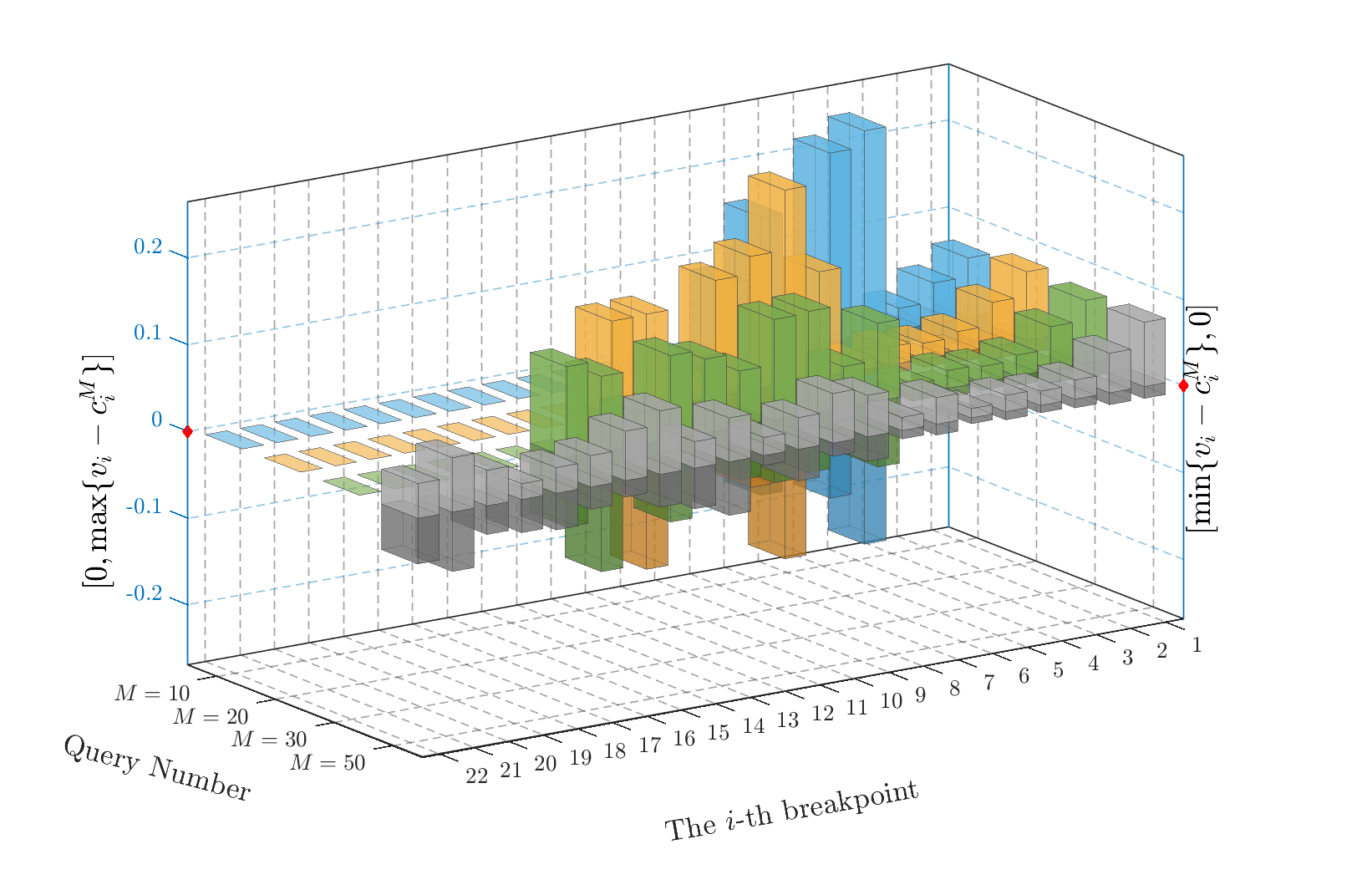}
\text{\footnotesize{\textbf{(A)}}}
\end{minipage}
\hfill
\begin{minipage}[t]{0.31\linewidth}
\centering
\includegraphics[width=2.0in]{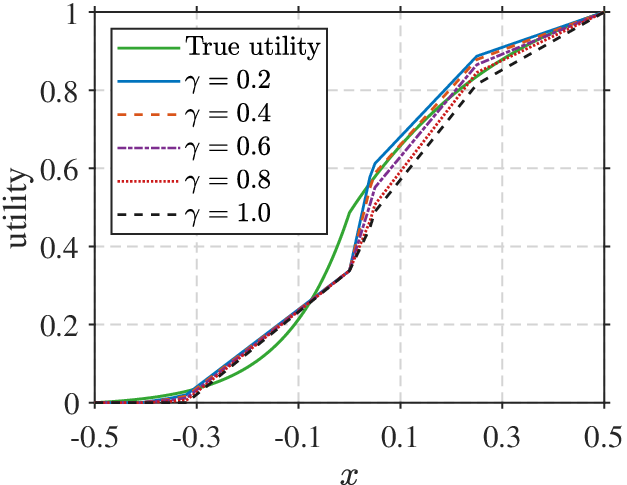}
\text{\footnotesize{\hspace{0.6cm}\textbf{(B1)}}}
\end{minipage}
\hfill
\begin{minipage}[t]{0.31\linewidth}
\centering
\includegraphics[width=2.0in]{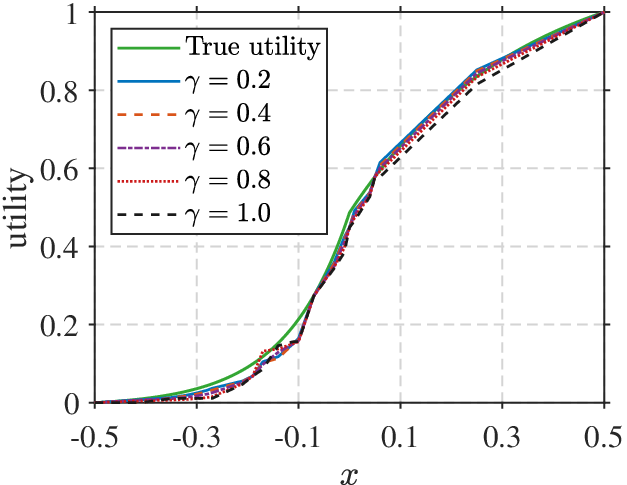}
\text{\footnotesize{\hspace{0.6cm}\textbf{(B2)}}}
\end{minipage}
\caption{\footnotesize{(A)
Increment value range. The interval bars of $v_{[N-1]}$ are 
defined as in \eqref{eq:bounds_vo-0}, 
that is,
$\underline{v}_i^M:=-\min
\left\{{v}_i-{c}_{i}^{M}:
{\bm v}_{[N-1]}={\bm R}{\bm v},{\bm v}\in {\cal V}_{N}^M\right\}$,
$\bar{v}_i^M:=\max\left\{{v}_i-{c}_{i}^{M}: {\bm v}_{[N-1]}
={\bm R}{\bm v},{\bm v}\in {\cal V}_{N}^M\right\}$.
(B1)-(B2) Worst-case utility functions under different values of $\gamma$ and number of queries: $M=10$ and $M=50$.
In each panel, the green curve represents the true utility function, while the
remaining curves represent the worst-case utility functions for
$\gamma\in\{0.2,0.4,0.6,0.8,1\}$.
}}
\vspace{-0.2cm}
\label{fig:feasible-se-2}
\end{figure}

We set ${\cal X}_0=\{-0.5,0,0.25,0.5\}$ with $N_0=4$,
$[\underline{p},\bar{p}]=[0.05,0.95]$, and with all breakpoints
rounded to two decimal places.
We construct the feasible set ${\cal V}_{N}^M$ of ${\bm v}$ using
the MPM in Section~\ref{sec:strategies}
and then define the set
$\widehat{\cal V}_{N}^M$ as in
\eqref{eq:minmin_bounds-a0}.
The aim
of the test is to examine the 
performance of the PRO 
model 
under the conservatism reduction scheme (CR-B) in 
\eqref{eq:PRO-PLA-N-dis-2}.
The number of queries $M$
takes values in $\{10,20,30,50\}$,
and the parameter $\gamma$ (which determines the vector of weights ${\bm \sigma}$ 
in \eqref{eq:PRO-PLA-N-dis-2}),
and thus 
controls
the level of conservatism of preference robust optimization 
model \eqref{eq:PRO-PLA-N-dis-2},
varies from 
$0.2$ to $1$.
The range of feasible set $\widehat{\cal V}_{N}^M$ can be visualized by the corresponding lower and upper utility values 
$
\underline{u}_i^M:=\min_{{\bm v}_{[N-1]}\in \widehat{\cal V}_{N}^M} {\bm v}_{[N-1]}^{\top}{\bm g}(x_i)$
and
$
\bar{u}_i^M:=\max_{{\bm v}_{[N-1]}\in \widehat{\cal V}_{N}^M} {\bm v}_{[N-1]}^{\top}{\bm g}(x_i), \inmat{ for }i\in [N].
$
Figure~\ref{fig:feasible-se-3}(a) 
depicts 
the region bounded by the
PLFs $\underline{u}:[\underline{x},\bar{x}]\to [0,1]$ and $\bar{u}:[\underline{x},\bar{x}]\to [0,1]$ passing through $(x_i,\underline{u}_i^M)$ and $(x_i,\bar{u}_i^M)$ for $i\in [N]$.
Compared with Figure~\ref{fig:PLA_Flex_1}(a) under the same parameter
setting, the utility range is further reduced, demonstrating the
effectiveness of the modified query design strategy.
Figure~\ref{fig:feasible-se-2}(b) displays the range $[-\underline{v}_i,\bar{v}_i]\subset [-1,1]$, for $i\in [N-1]$.
Figure~\ref{fig:feasible-se-2}(a)-\ref{fig:feasible-se-2}(c)
shows the change in
the worst-case utility function
with $\gamma\in \{0.2,0.4,0.6,0.8,1\}$ for 
$M=\{10,50\}$.

To facilitate the comparison of the two
schemes to reduce
(CR-A) in \eqref{eq:min_expand-20} and (CR-B) in \eqref{eq:robust-pi0}, we define the optimal value functions associated
with their respective inner problems. For (CR-A), let
\begin{equation}
\rho_{\Gamma}^{\mathrm{A}}(\bm z)
:=
\min_{\substack{
\bm v_{[N-1]}\in\widehat{\mathcal V}_N^M,\,
\bm y}}
\mathbb{E}
\left[
\bm v_{[N-1]}^\top g(f(\bm z,\bm\xi))
\right]
\quad
\text{s.t. } \eqref{eq:min_expand-20-a}\text{--}\eqref{eq:min_expand-20-b},
\label{eq:rho_Gamma-0}
\end{equation}
where $\Gamma$ controls the level of conservatism in (CR-A).
Similarly, for (CR-B), let
\begin{equation}
\rho_{\gamma}^{\mathrm{B}}(\bm z)
:=
\min_{\substack{
\bm v_{[N-1]}\in\widehat{\mathcal V}_N^M,\,
\bm\Pi\in\mathcal S}}
\mathbb{E}
\left[
\bm v_{[N-1]}^\top g(f(\bm z,\bm\xi))
\right]
\quad
\text{s.t. } \eqref{eq:robust-pi0-b},
\label{eq:rho_B}
\end{equation}
where $\gamma$ determines the weight vector $\bm\sigma$ in \eqref{eq:sigma0} and
thereby controls the level of conservatism in (CR-B).
The corresponding preference robust optimization problems can then be
written~compactly~as
\begin{equation}
\max_{\bm z\in Z}\rho_{\Gamma}^{\mathrm{A}}(\bm z)
\qquad\text{and}\qquad
\max_{\bm z\in Z}\rho_{\gamma}^{\mathrm{B}}(\bm z).
\label{eq:PRO-PLA-N-dis-Gamma}
\end{equation}

Figure~\ref{fig:feasible-se-3}(a) further shows that the optimal value under (CR-B) decreases approximately linearly with \(\gamma\). 
By contrast, Figure~\ref{fig:feasible-se-3}(b) shows a nonlinear decrease in the optimal value under (CR-A) as \(\Gamma\) increases.
To investigate the effect of the weight vector in (CR-B) in \eqref{eq:PRO-PLA-N-dis-2} further, 
we replace ${\bm \sigma}$ in \eqref{eq:sigma0} with ${\bm \sigma}\circ {\bm e}_s$,
where ${\bm e}_s=(\overbrace{1,\cdots,1}^{s},0,\cdots,0)\in \R^{N-1}$,
and subsequently, replace ${\bm \sigma}$ with ${\bm \sigma}\circ {\bm e}_s$ in a reformulation \eqref{eq:reformulation_conser-a}. Figure~\ref{fig:optimal values_compare}(a)-\ref{fig:optimal values_compare}(c) show the change of optimal values of problem \eqref{eq:reformulation_conser-a} with 
$s\in\{1,\cdots,19\}$ for fixed $\gamma\in \{0.1,0.5,0.9\}$ respectively.
Figure~\ref{fig:optimal values_compare2}(a)--(c) show the optimal values of problem
\eqref{eq:reformulation_conser-a} as $s$ changes for
$\gamma=0.001$, $0.2$, $0.4$, $0.6$, $0.8$, and $1$, where panels (a)-(c) correspond to $M=10$, $M=30$, and $M=50$,
respectively.

To further examine the performance of the MPM
and the UPRO model,
we have conducted numerical tests on two additional utility functions that exhibit different shapes and the out-of-sample performance of PRO model.
We begin by presenting the utility range, increment bounds, and changes in the optimal values of model (CR-B) for two additional utility functions, followed by an analysis of the out-of-sample performance of model (CR-B) and a comparison between MPM and RRUS.

\vspace{-0.1cm}
\subsubsection{In-Sample Performance of (CR-B) under Alternative S-Shaped and Concave Utilities}
\vspace{-0.1cm}

In this section, we consider two more utility functions to illustrate the performance of
the MPM
and the UPRO model (CR-B).

\noindent{\bf S-shaped Utility Function}.
We consider another type of S-shaped utility function
\bgeqn
\label{eq:S2_utility}
u_2(x)=
\begin{cases}
x^\alpha, & x\geq0,\\
-\lambda(-x)^\beta, & x<0,
\end{cases}
\qquad
\alpha=0.6,\quad\beta=0.48,\quad\lambda=2.25.
\edeqn

\begin{figure}[!h]
\vspace{-1em}
\begin{minipage}[t]{0.32\linewidth}
\centering
\includegraphics[width=2.2in]{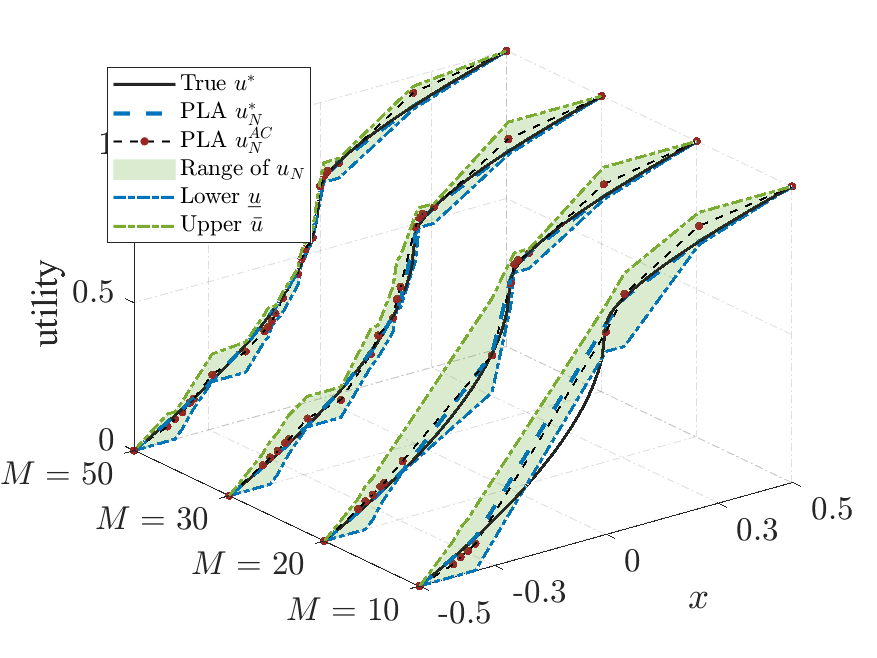}
\text{\footnotesize{\hspace{0.3cm}\textbf{(a)}}}
\end{minipage}
\hfill
\begin{minipage}[t]{0.32\linewidth}
\centering
\includegraphics[width=2.4in]{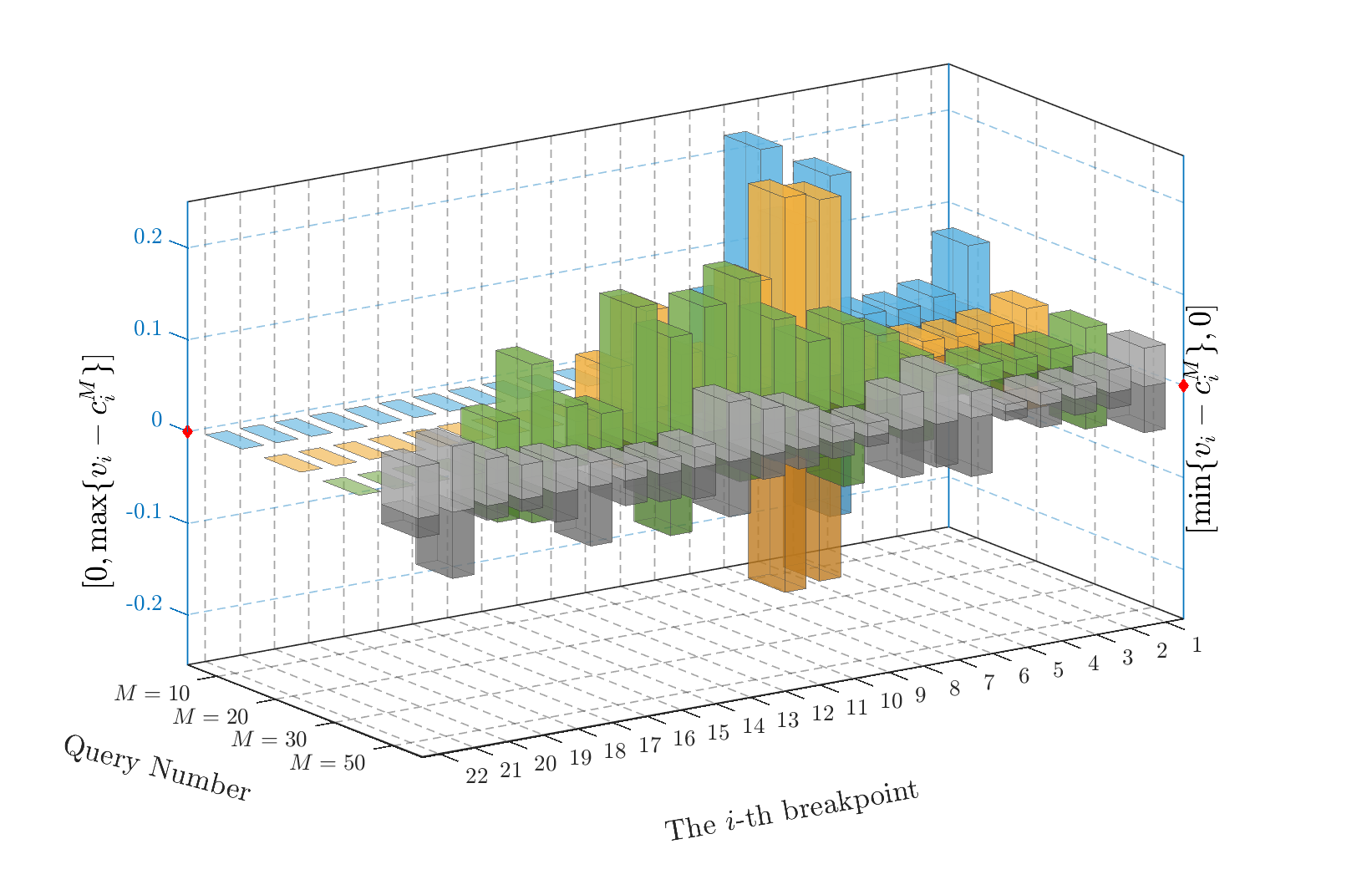}
\text{\footnotesize{\hspace{0.6cm}\textbf{(b)}}}
\end{minipage}
\hfill
\begin{minipage}[t]{0.32\linewidth}
\centering
\includegraphics[width=1.9in]{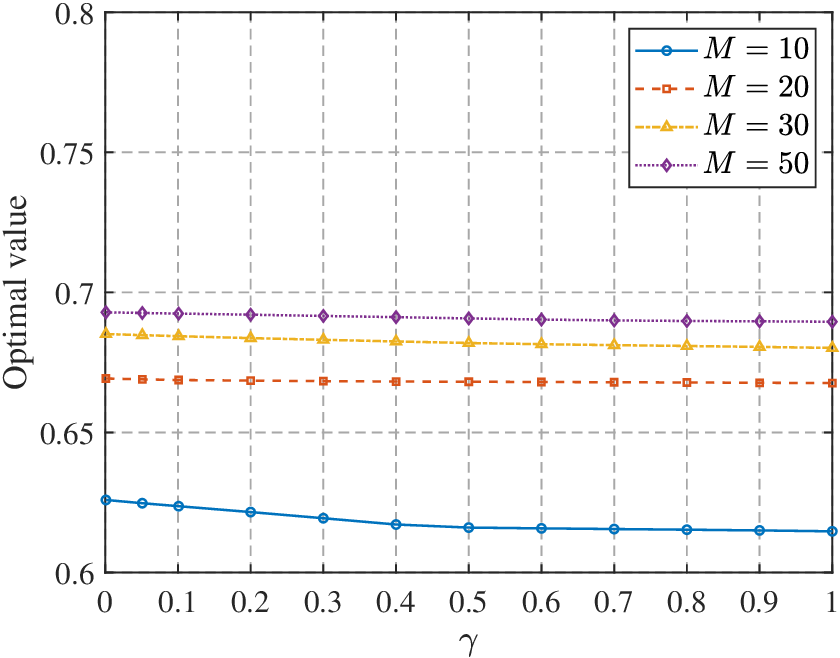}
\text{\footnotesize{\hspace{0.6cm}\textbf{(c)}}}
\end{minipage}
\caption{\footnotesize{ S-shaped utility function:
(a) The range of feasible utility functions as $M$ increases.
(b) The ranges of the utility increments $v_{[N-1]}$ for different values of $M$.
(c) The optimal values of (CR-B) for~different~values~of~$\gamma$.
}}
\label{fig:feasible-91_range}
\end{figure}

Figure~\ref{fig:feasible-91_range}(a) 
depicts the utility-value range bounded by the lower and upper PLFs, $\underline{u}$ and $\bar{u}$,
respectively. Figure~\ref{fig:feasible-91_range}(b)  
displays the corresponding increment-value
ranges $[-\underline{v}_i,\bar{v}_i]$ for $i\in[N-1]$.
Figure~\ref{fig:feasible-91_range}(c) further shows that the optimal value under (CR-B) decreases
approximately linearly with $\gamma$.

\noindent{\bf Concave Utility Function}.
We next consider a concave utility function
\bgeqn
\label{eq:concave_utility}
u_1(x)=1-\exp(-\gamma x),\qquad \gamma=10.
\edeqn

\begin{figure}[!h]
\vspace{-1em}
\begin{minipage}[t]{0.32\linewidth}
\centering
\includegraphics[width=2.2in]{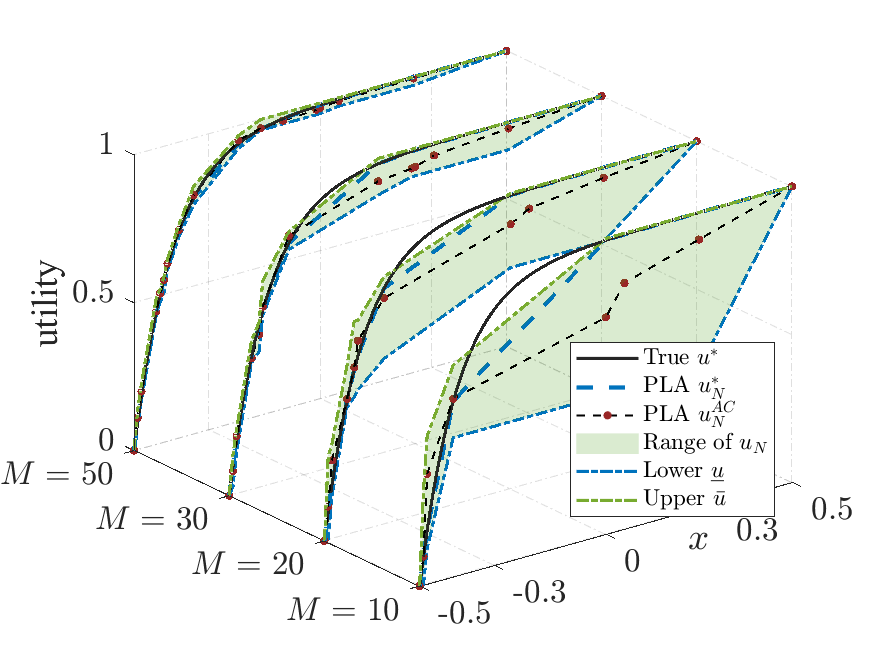}
\text{\footnotesize{\hspace{0.3cm}\textbf{(a)}}}
\end{minipage}
\hfill
\begin{minipage}[t]{0.32\linewidth}
\centering
\includegraphics[width=2.4in]{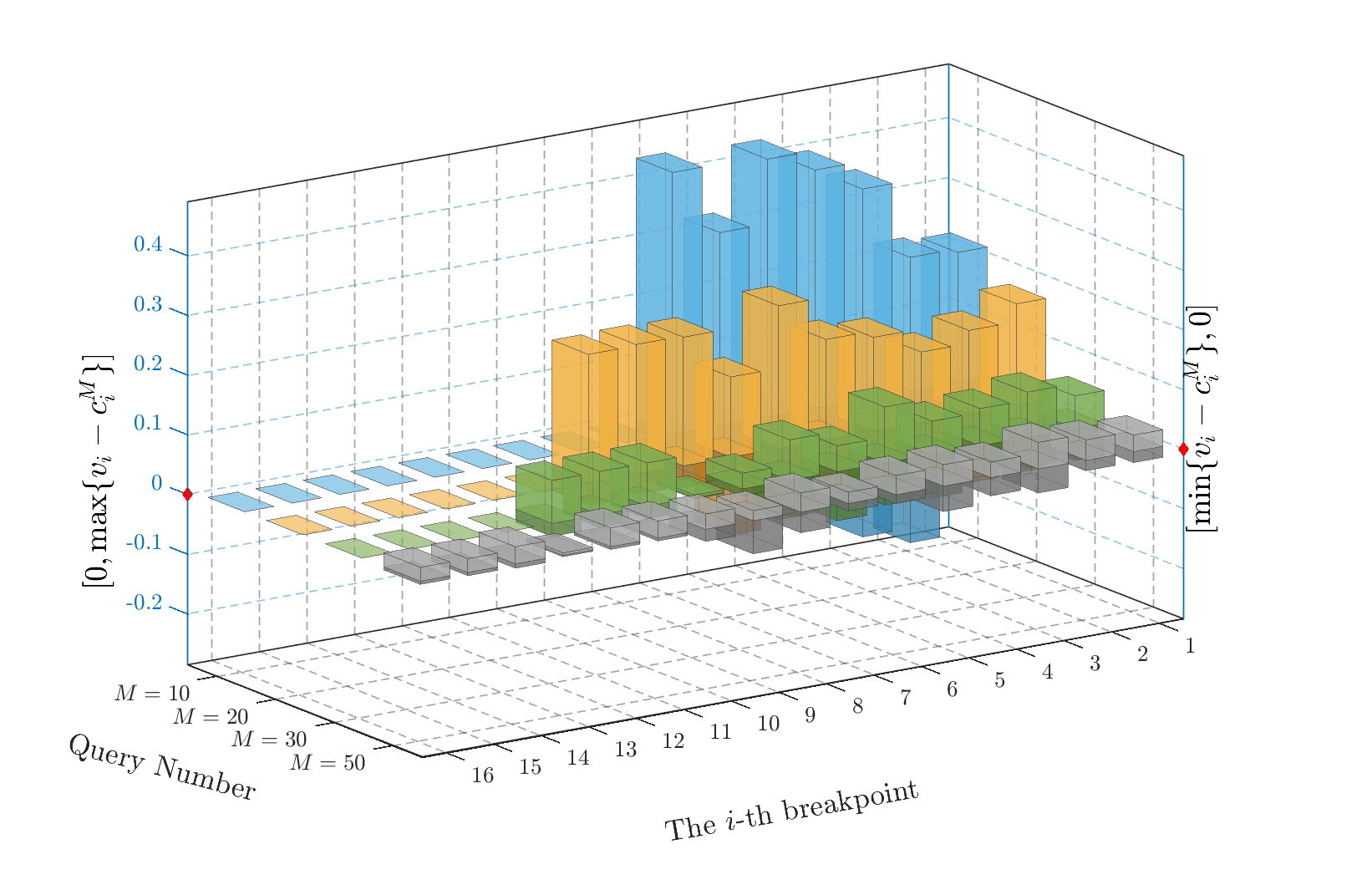}
\text{\footnotesize{\hspace{0.6cm}\textbf{(b)}}}
\end{minipage}
\hfill
\begin{minipage}[t]{0.32\linewidth}
\centering
\includegraphics[width=1.9in]{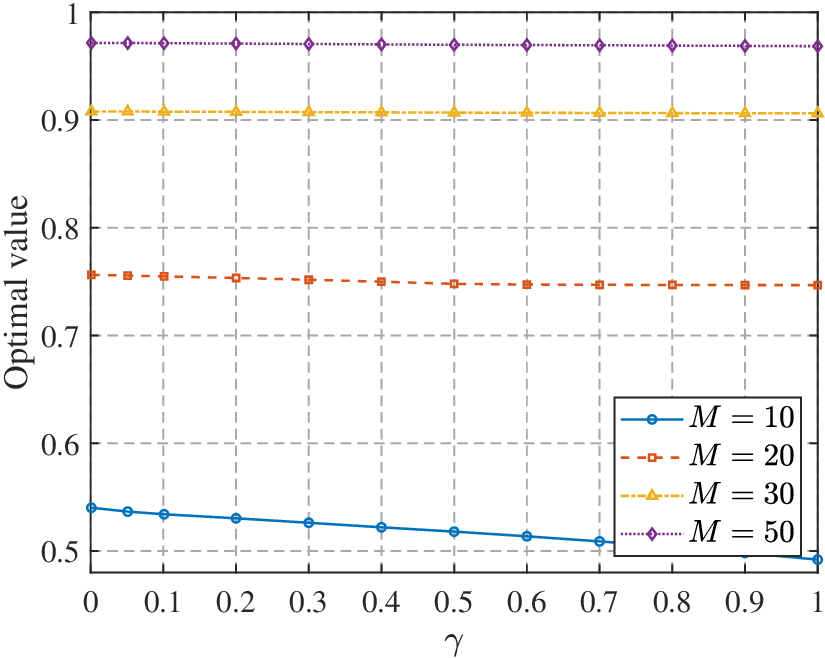}
\text{\footnotesize{\hspace{0.6cm}\textbf{(c)}}}
\end{minipage}
\caption{\footnotesize{ Concave utility function: 
(a) The range of feasible utility functions as $M$ increases.
(b) The ranges of the utility increments $v_{[N-1]}$ for different values of $M$.
(c) The optimal values of (CR-B) for different~values~of~$\gamma$.
}}
\label{fig:feasible-91_range3}
\vspace{-0.5cm}
\end{figure}

Figure~\ref{fig:feasible-91_range3}(a) 
depicts the utility-value range bounded by the lower and upper PLFs,
$\underline{u}$ and $\bar{u}$,
respectively. Figure~\ref{fig:feasible-91_range3}(b)  
displays the corresponding increment-value
ranges $[-\underline{v}_i,\bar{v}_i]$ for $i\in[N-1]$.
Compared with Figure~\ref{fig:Fig_accuracy_50}(a) under the same parameter setting, the utility range is further reduced, demonstrating the effectiveness of the modified query design strategy proposed at the end of Section~\ref{sec:query_range}.
Figure~\ref{fig:feasible-91_range3}(c) further shows that the optimal value under (CR-B) decreases
approximately linearly with $\gamma$.

\vspace{-0.1cm}
\subsubsection{Out-of-Sample Performance of (CR-B) under Different Utility Functions}
\label{sec:OoS_utility}
\vspace{-0.1cm}

We conduct a
rolling-window out-of-sample experiment using 10 financial indices in Section~\ref{Sec:Numberial}.
We use monthly observations from {January 2015 to December 2019}. We carry out out-of-sample tests using a rolling window of 24 months, so that $K=24$ in each optimization problem. Specifically, we first use the observations from {January 2015 to December 2016} to determine the optimal portfolio strategies under (CR-B), and evaluate the resulting strategies using the realized returns in {January 2017}. We then roll the estimation window forward by one month, using the observations from {February 2015 to January 2017} to determine the portfolio strategies for {February 2017}. This procedure is repeated on a monthly basis until the end of the out-of-sample period. 

\begin{figure}[!h]
\centering
\begin{minipage}[t]{0.32\linewidth}
\centering
\includegraphics[width=2.1in]{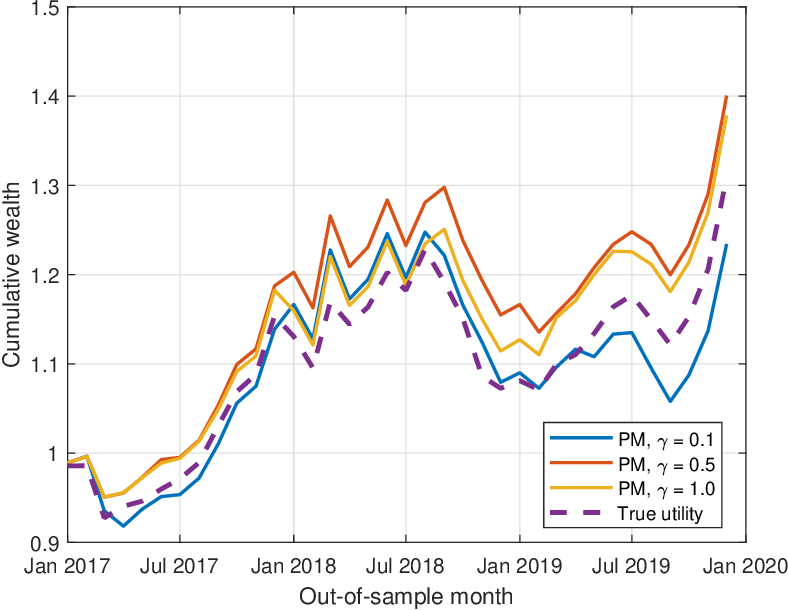}
\text{\footnotesize{\hspace{0.3cm}\textbf{(a)}}}
\end{minipage}
\hfill
\begin{minipage}[t]{0.32\linewidth}
\centering
\includegraphics[width=2.1in]
{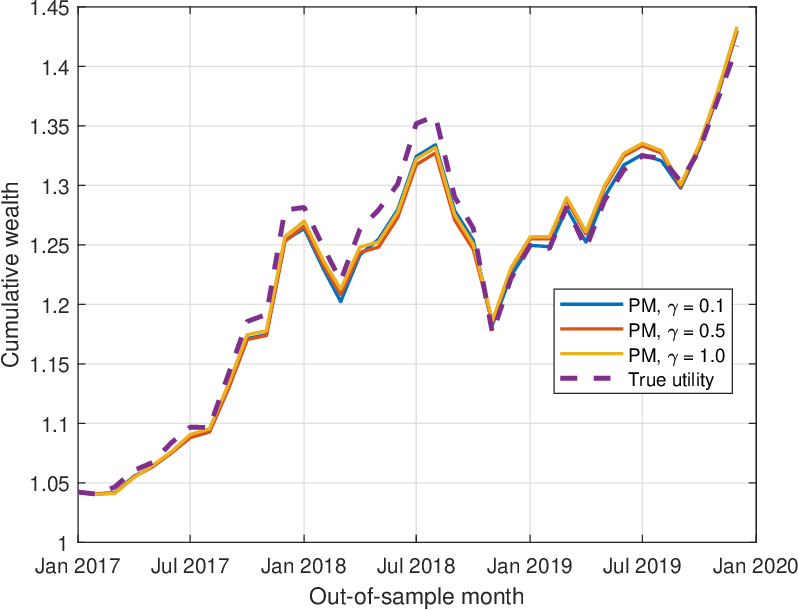}
\text{\footnotesize{\hspace{0.3cm}\textbf{(b)}}}
\end{minipage}
\hfill
\begin{minipage}[t]{0.32\linewidth}
\centering
\includegraphics[width=2.1in]{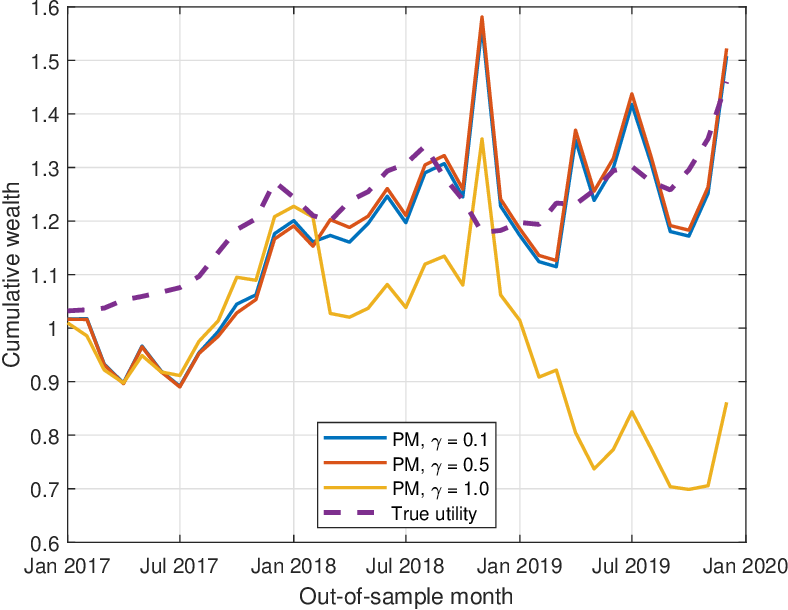}
\text{\footnotesize{\hspace{0.3cm}\textbf{(c)}}}
\end{minipage}
\caption{\footnotesize
Out-of-sample performance of the portfolio decisions obtained by the proposed polyhedral method under different values of $\gamma$, compared with the true-utility benchmark, for $M=50$ queries: (a) S-shaped utility in Section~\ref{Sec:Numberial}; (b) S-shaped utility \eqref{eq:S2_utility} (c) concave utility \eqref{eq:concave_utility}.
}
\label{fig:OOS_comparison2-4}
\vspace{-0.2cm}
\end{figure}

We consider the S-shaped utility function in  Section~\ref{Sec:Numberial}, together with the two additional utility functions in \eqref{eq:S2_utility} and \eqref{eq:concave_utility}, as discussed above. 
Figure~\ref{fig:OOS_comparison2-4} (a)-(b) depicts the out-of-sample performance of the PRO model under both utility functions with $M=50$ queries. 
We observe that the MPM yields stable portfolio performance across the considered values of $\gamma$ in both utility cases. 
Moreover, the cumulative wealth trajectories generated by the model (CR-B) remain close to the true-utility benchmark and exhibit relatively limited sensitivity to the choice of~$\gamma$.

\begin{figure}[!h]
\begin{minipage}[t]{\linewidth}
\centering
\includegraphics[width=2.2in]{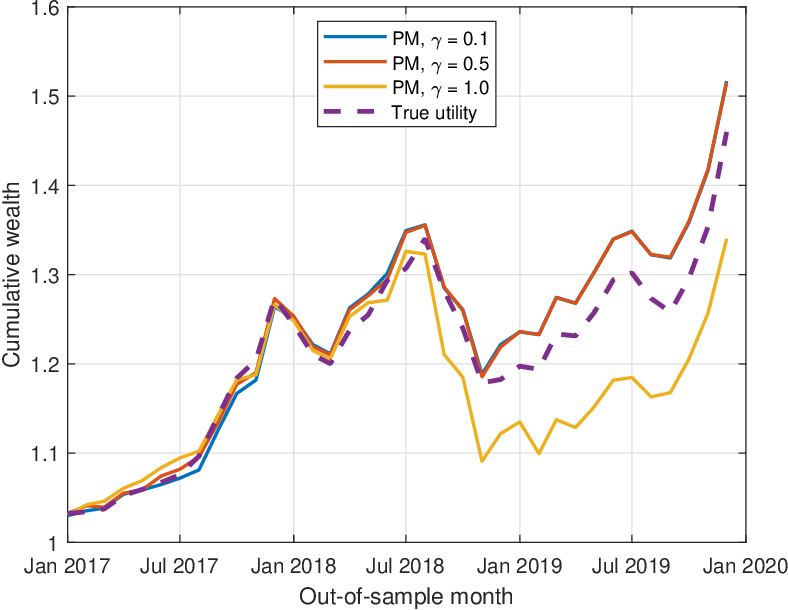}
\end{minipage}
\caption{\footnotesize
Concave utility: Out-of-sample performance of the portfolio decisions obtained by the proposed polyhedral method with the portfolio variance constraint under different values of $\gamma$, compared with the true-utility benchmark, for $M=50$ queries.
}
\label{fig:OOS_comparison2-4_2}
\vspace{-0.5cm}
\end{figure}

In Figure~\ref{fig:OOS_comparison2-4}(c), while the results for $\gamma=0.5$ and $\gamma=0.1$ remain relatively close to the true-utility benchmark, a more noticeable deviation is observed for $\gamma=1$. To investigate whether controlling portfolio risk leads to more stable out-of-sample performance, we further consider a risk-constrained portfolio setting by imposing an upper bound on the portfolio variance. Specifically, for each rolling window, we impose
$
{\bm x}^{\top}\widehat{\bm\Sigma}_t{\bm x}
\leq \bar{\sigma}_t^2,
$
where $\widehat{\bm\Sigma}_t$ denotes the sample covariance matrix
estimated from the returns in the corresponding rolling window, and
$\bar{\sigma}_t^2=2\sigma_{\min,t}^2$, with $\sigma_{\min,t}^2$
denoting the minimum portfolio variance in that window.
As shown in Figure~\ref{fig:OOS_comparison2-4_2}, imposing the variance constraint reduces the deviation from Figure~\ref{fig:OOS_comparison2-4}(c), with the cumulative wealth trajectories remaining closer to the true-utility benchmark.

\vspace{-0.1cm}
\subsubsection{Out-of-Sample Performance under MPM and RRUS}
\vspace{-0.1cm}

We compare the RRUS scheme and the MPM in the PRO model with conservatism control, and denote the former PRO problem by (PRO-RRUS).
To evaluate the performance of the portfolio decisions obtained from (PRO-RRUS) and (CR-B), we conduct the out-of-sample tests described in Section~\ref{sec:OoS_utility}.
We consider the utility function in Section~\ref{Sec:Numberial} and set $M=30$, since the dimension of the piecewise linear utility representation in RRUS is $3M+3$, making larger values of $M$ computationally more expensive.

\begin{figure}[!h]
\centering
\begin{minipage}[t]{0.32\linewidth}
\centering
\includegraphics[width=2.0in]{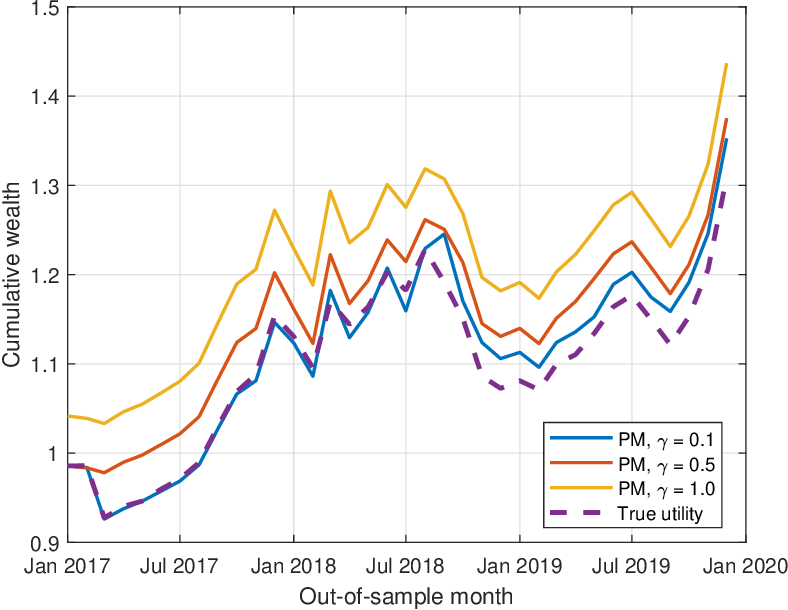}
\text{\footnotesize{\hspace{0.1cm}\textbf{(A)}}}
\end{minipage}
\begin{minipage}[t]{0.32\linewidth}
\centering
\includegraphics[width=2.0in]{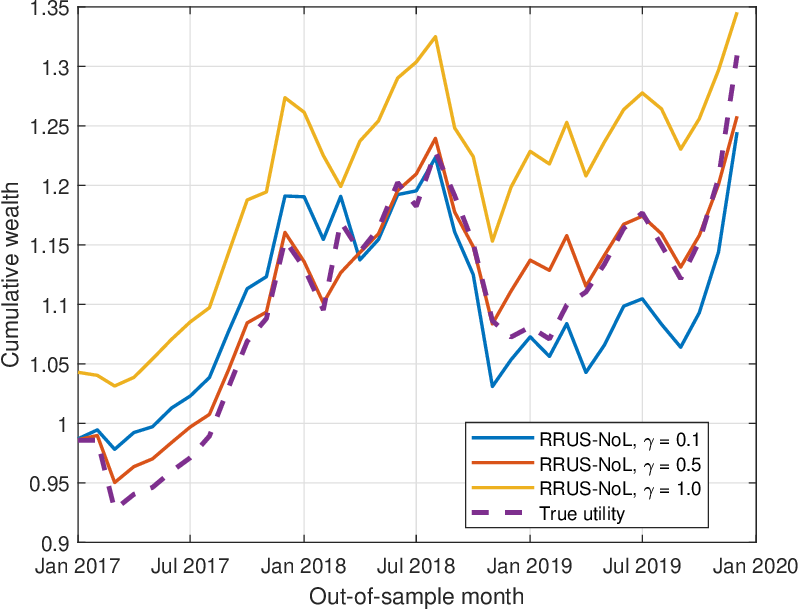}
\text{\footnotesize{\hspace{0.1cm}\textbf{(B1)}}}
\end{minipage}
\begin{minipage}[t]{0.32\linewidth}
\centering
\includegraphics[width=2.0in]{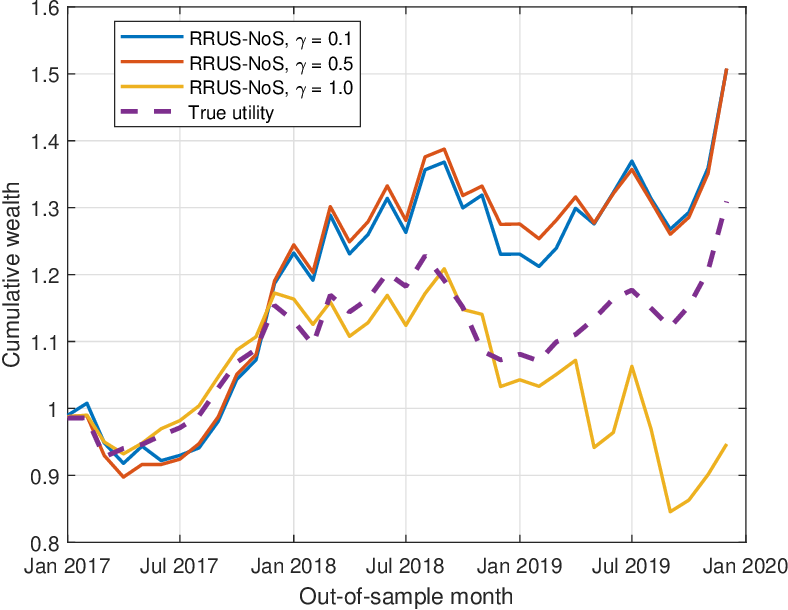}
\text{\footnotesize{\hspace{0.1cm}\textbf{(B2)}}}
\end{minipage}
\caption{\footnotesize{
Out-of-sample performance of the portfolio decisions for $M=30$:
(A) using the proposed MPM;
(B1) using RRUS with $S$-shape information but without information on the Lipschitz constant $L$; and
(B2) using RRUS with the Lipschitz constant $L$ but without $S$-shape information.}}
\label{fig:obj_RRUS1}
\vspace{-0.3cm}
\end{figure}

Figure~\ref{fig:obj_RRUS1}(a) depicts the out-of-sample performance of the portfolio decisions obtained from (CR-B).
Likewise, Figures~\ref{fig:obj_RRUS1}(b)--(c) depict the corresponding out-of-sample performance of (PRO-RRUS) when only $S$-shape information or only information on the Lipschitz modulus is available, respectively.
Comparing the out-of-sample results, we observe that the proposed polyhedral method yields relatively stable portfolio performance across the considered values of $\gamma$.
By contrast, the performance of the RRUS-based portfolios appears to be more sensitive to the degree of robustness. This effect is especially visible when only information on the Lipschitz modulus is available (Figure~\ref{fig:obj_RRUS1}(c)), where the fully robust model ($\gamma=1$) yields lower cumulative out-of-sample wealth, while the moderately robust specifications ($\gamma=0.1$ and $\gamma=0.5$) remain substantially more stable.
These results suggest that MPM yields more stable out-of-sample portfolio performance across different values of $\gamma$ in the considered experiments.

\vspace{-0.2cm}
\section{Extensions}
\label{sec:extension}
\vspace{-0.2cm}

In this section, we 
outline several potential extensions of the MPM,
including its integration with the 
bisection method for preference elicitation, its extension from 
single-attribute to multi-attribute nonlinear utility functions, and its 
application to dynamic pricing problems.

\vspace{-0.1cm}
\subsection{A Hybrid Elicitation Scheme}
\label{sec:hybrid_Bisection}
\vspace{-0.1cm}

Prompted by a reviewer's suggestion to adopt a bisection strategy, we propose a hybrid scheme that integrates the bisection approach with
our MPM to further enhance its performance. This subsection describes the details of the proposed scheme.
Let 
$\underline{u}(x)$ and $\bar{u}(x)$ denote the minimum and maximum 
feasible utility values at $x\in[\underline{x},\bar{x}]$, respectively.
Given the current breakpoint set ${\cal X}$, we identify a 
breakpoint in the set with the largest 
uncertainty in utility,
i.e.,
$x^{\mathrm B}
\in
\arg\max_{x\in{\cal X}}
\left\{
\bar{u}(x)-\underline{u}(x)
\right\}.
$
We then construct a binary-choice query 
that compares
a certain outcome
$Z_2^{\rm B}=x^{\mathrm B}$ with a risky lottery $Z_1^{\rm B}$ that yields 
$\underline{x}$ with probability $1-p^{\rm B}$ and $\bar{x}$ with 
probability $p^{\rm B}$, where
$
p^{\mathrm B}
=
({
\underline{u}(x^{\mathrm B})
+
\bar{u}(x^{\mathrm B})
})/{2}$.
Here we use the superscript ``B" to emphasize ``Bisection".  
Suppose $Z_1^{\rm B}\preceq Z_2^{\rm B}$. 
Under the bisection scheme,
the range of the feasible set of utility functions at point 
$x^{\mathrm B}$ is reduced to
$[p^{\rm B},\bar{u}(x^{\mathrm B})]$; otherwise, it is reduced to
$[\underline{u}(x^{\mathrm B}),p^{\rm B}]$.
Therefore, each query halves the range of feasible utility functions at $x^{\mathrm B}$.

Both the bisection method and the strategy in Section~\ref{sec:query_range}
reduce a wide utility range, but in different ways. The latter uses utility estimate $u_{N_m}^{\rm AC}$ and satisfies the choice-balance condition in Theorem~\ref{prop:EDn=0-2}, 
consistent with complexity analysis in
\cite[Section 3.2]{boyd2003localization}, whereas bisection directly
reduces the utility value range at any specified point by half.
However, when the shape and the Lipschitz modulus
of $u^*$ are unknown, the bisection method 
does not provide an adaptive mechanism for determining where new
breakpoints should be introduced over the domain to
reduce the overall ambiguity 
of $u^*$ efficiently.
In particular, when 
the bisection method is applied to equally spaced points over the domain, 
it may not provide sufficient resolution in regions where the utility function changes sharply.
By contrast,
identifying breakpoints with approximately {\em equally spaced utility values} can improve the efficiency of the bisection method. 
Our MPM can facilitate the adaptive selection of breakpoint locations,
which motivates a hybrid scheme combining the MPM and bisection method.

We consider six strategies for the hybrid scheme applied to a specific problem,
where every five queries consist of 
$k$-bisection queries and $(5-k)$-polyhedral queries, for
$k\in\{0,\ldots,5\}$. 
Note that $k=0$ corresponds to the proposed 
polyhedral method and 
$k=5$ to the bisection method.
$k\in\{1,\ldots,4\}$
reflects 
different degrees of hybridization.
For instance, 
when $k=2$, three consecutive queries are generated 
by the polyhedral method, 
followed by two bisection queries.  
The five-query sequence is repeated throughout the elicitation 
process.
When $k=4$, one polyhedral query is followed by four 
consecutive~bisection~queries,~and~this~sequence~is~repeated.

Next we examine the performance of the six strategies through 
a concave utility function.
For both examples, we set $N_0=4$ with
${\cal X}_0=\{-0.5,0,0.25,0.5\}$.
We consider the concave utility function
$
u^*(x)
=1-\exp(-10x),
$
for $x\in[-0.5,0.5].
$
In the numerical experiment, we set $r_1,r_2,r_3\in[-0.5,0.5]$ and $p\in[0.01,0.99]$, all 
to two decimal places.
\begin{figure}[!ht]
\vspace{-0.3cm}
\begin{minipage}[t]{0.32\linewidth}
\centering
\includegraphics[width=2.0in]{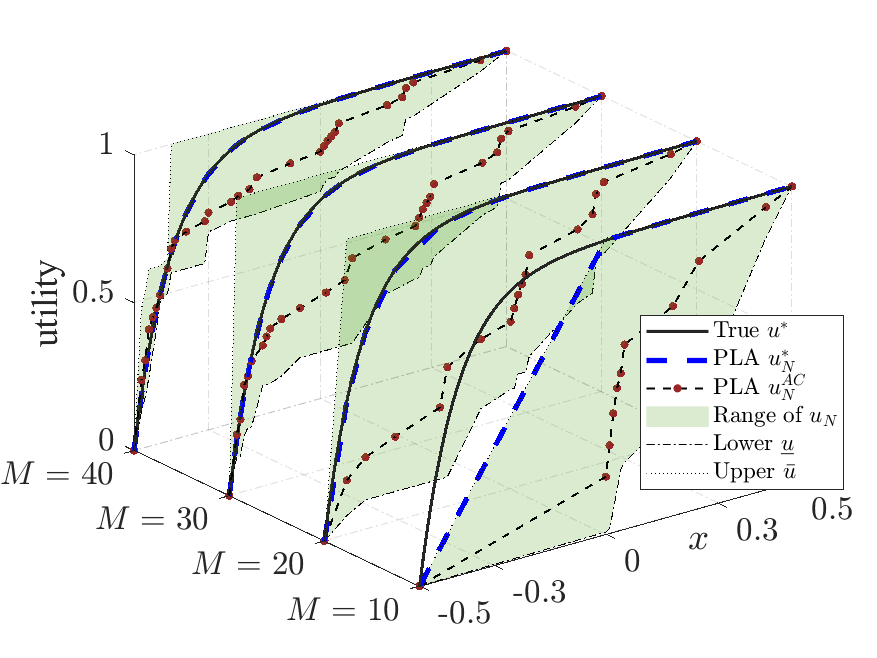}
\text{\footnotesize{\hspace{-0.1cm}\textbf{(a)}}}
\end{minipage}
\hfill
\begin{minipage}[t]{0.32\linewidth}
\centering
\includegraphics[width=2.0in]{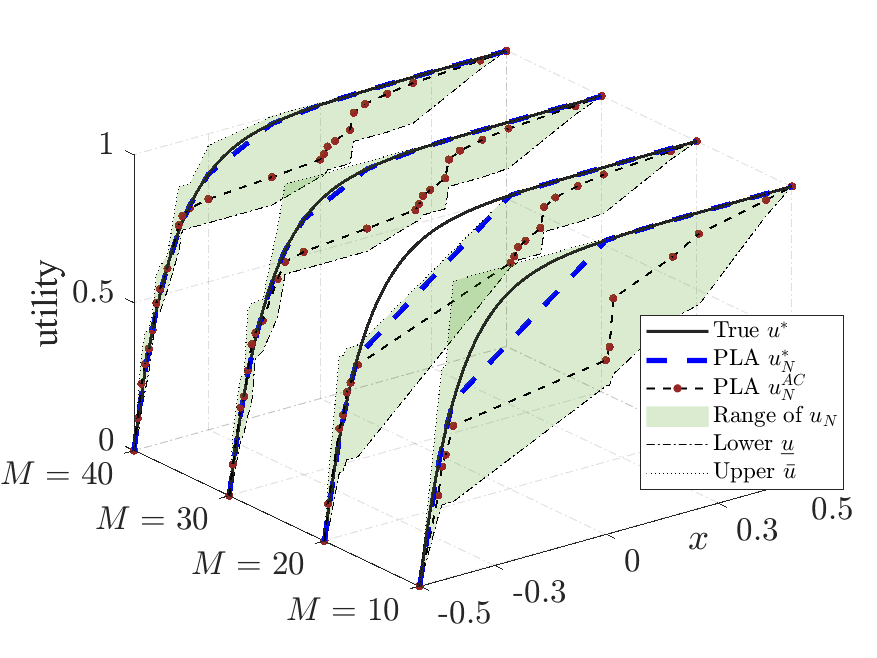}
\text{\footnotesize{\hspace{-0.1cm}\textbf{(b)}}}
\end{minipage}
\begin{minipage}[t]{0.32\linewidth}
\centering
\includegraphics[width=2.0in]{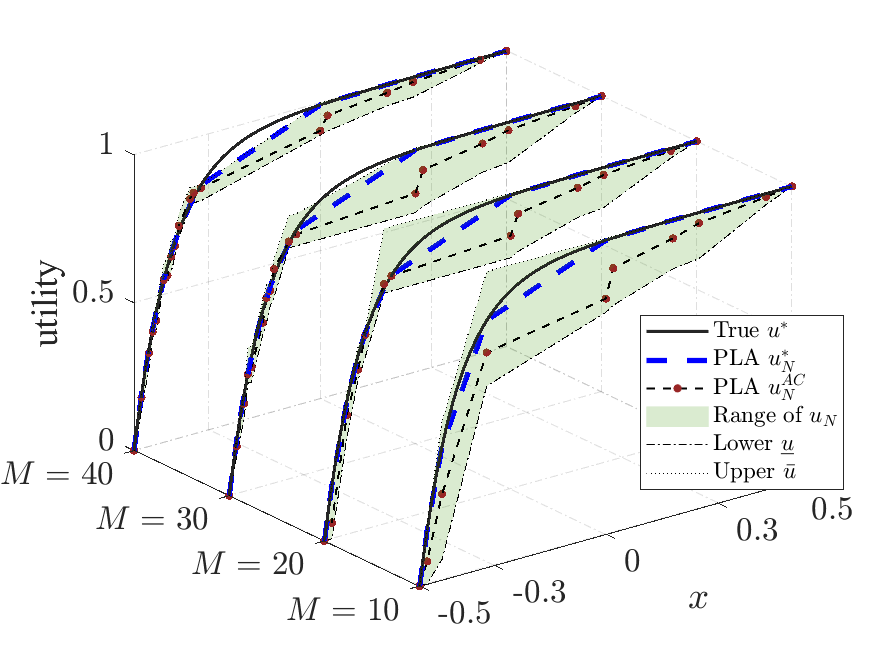}
\text{\footnotesize{\hspace{-0.1cm}\textbf{(c)}}}
\end{minipage}
\hfill
\begin{minipage}[t]{0.32\linewidth}
\centering
\includegraphics[width=2.0in]{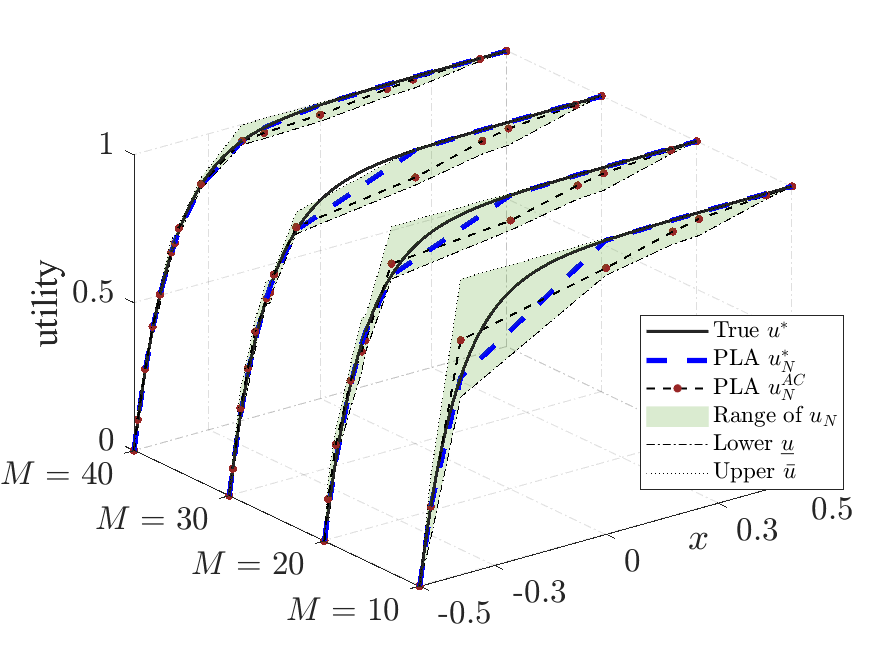}
\text{\footnotesize{\hspace{-0.1cm}\textbf{(d)}}}
\end{minipage}
\begin{minipage}[t]{0.32\linewidth}
\centering
\includegraphics[width=2.0in]{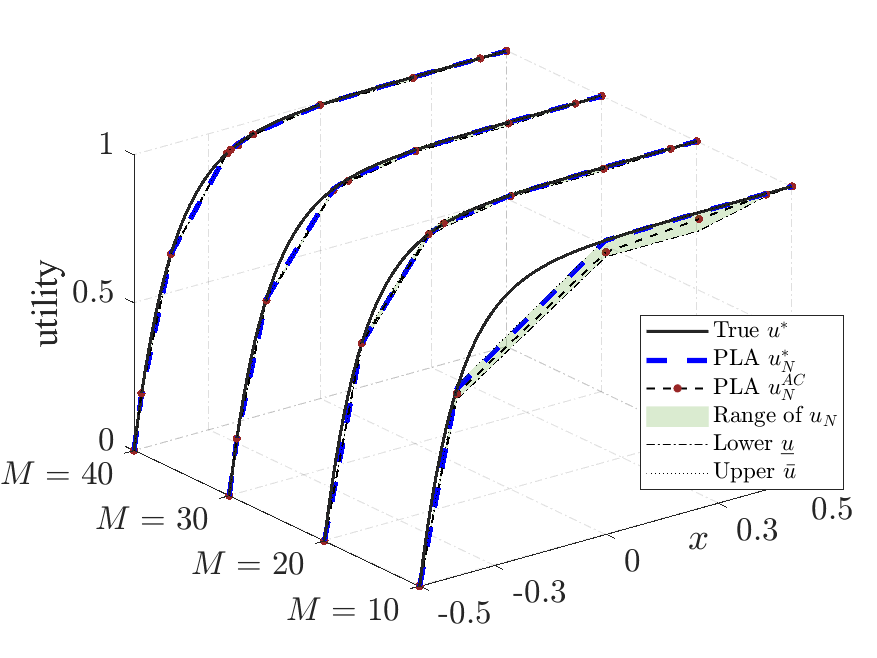}
\text{\footnotesize{\hspace{-0.1cm}\textbf{(e)}}}
\end{minipage}
\hfill
\begin{minipage}[t]{0.32\linewidth}
\centering
\includegraphics[width=2.0in]{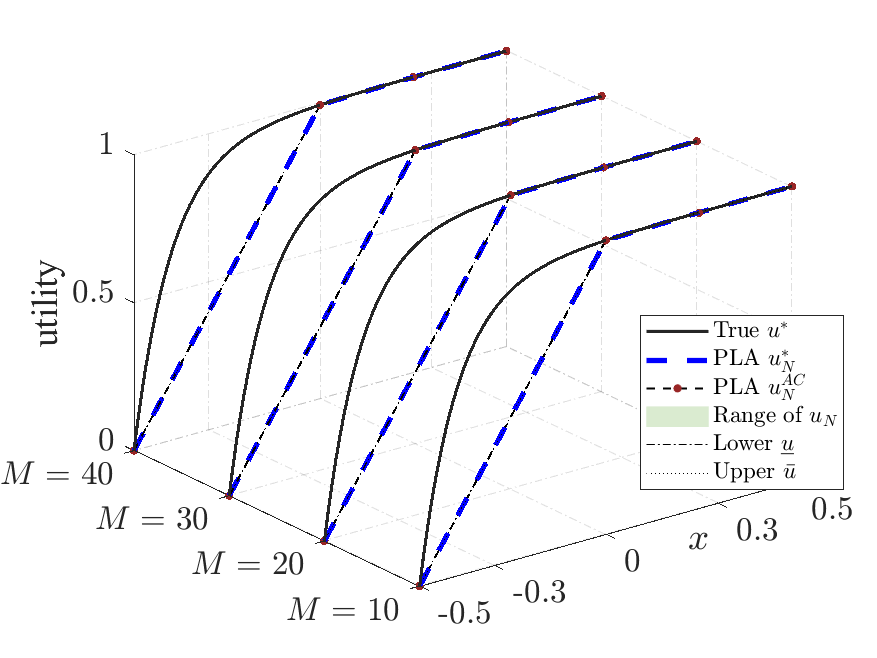}
\text{\footnotesize{\hspace{-0.1cm}\textbf{(f)}}}
\end{minipage}
\caption{\footnotesize{
Concave utility elicitation under different combinations of the polyhedral and bisection methods:
(a) $k=0$: five polyhedral queries;
(b) $k=1$: four polyhedral and one bisection queries;
(c) $k=2$: three polyhedral and two bisection queries;
(d) $k=3$: two polyhedral and three bisection queries;
(e) $k=4$: one polyhedral and four bisection queries;
(f) $k=5$: five bisection queries
in each group of five queries.
}}
\vspace{-0.3cm}
\label{fig:hybrid-elicitation}
\end{figure}

Figure~\ref{fig:hybrid-elicitation}(a) and (f) 
depict the outcomes
of the MPM and the Bisection method, respectively.
As shown in Figure~\ref{fig:hybrid-elicitation}(a), the MPM continuously introduces new breakpoints while gradually reducing the feasible range of utility values. 
By contrast, the bisection method in Figure~\ref{fig:hybrid-elicitation}(f) rapidly reduces the utility uncertainty at the existing breakpoints, but does not introduce new breakpoints.

Figures~\ref{fig:hybrid-elicitation}(b)--(e) depict the transition between these two non-hybrid cases as the number of bisection queries within each group of five queries increases from 1 to 4. 
A larger proportion of bisection queries generally leads to tighter utility bounds at the elicited breakpoints but also introduces
in fewer newly added breakpoints. Hence, the hybrid schemes reveal a trade-off between reducing utility uncertainty at existing breakpoints and introducing new breakpoints. In particular, the hybrid schemes combine the ability of the MPM to introduce breakpoints adaptively with the ability of the bisection method to tighten the utility bounds rapidly at some selected breakpoints.
Similar improvements are observed for S-shaped utility functions.

\vspace{-0.2cm}
\subsection{Multi-attribute Utility Elicitation}
\label{sec:multi}
\vspace{-0.2cm}

We begin by considering a general 
nonlinear parametric utility function
$U(\cdot,{\bm v}):\R^{n_1}\to [0,1]$ which
is linear in the parameter vector
${\bm v}\in \R^{n_2}$.
 We assume that there exists a unique 
$\bm v^*$ such that the DM's true utility function is represented by
$U(\cdot,\bm v^*)$.
Let $\widetilde{\mathcal V}_0$ be an initial polyhedron containing
$\bm v^*$. We iteratively refine it
by generating a
cut at each iteration.
At iteration $m$, we construct two lotteries ${\bm Z}_1^m$ and ${\bm Z}_2^m$,
where ${\bm Z}_1^m$ yields $\bm r_1^m$ with probability $1-p^m$ and
$\bm r_3^m$ with probability $p^m$, whereas ${\bm Z}_2^m$ yields
$\bm r_2^m$ with certainty.
To construct these lotteries, we compute the analytic center
$\bm c^m$ of $\widetilde{\mathcal V}_{m-1}$ and identify the two intersection
points $\bm v_1^m$ and $\bm v_2^m$ between the longest axis
of Sonnevend's ellipsoid 
and the
boundaries of $\widetilde{\mathcal V}_{m-1}$, as described in
Section~\ref{sec:two_problem}.
Analogous to 
\eqref{eq:gene_m_Q},
we determine
 ${\bm r}_1^m, {\bm r}_2^m, {\bm r}_3^m$ and $p^m$
by solving the following optimization problems:
\begin{subequations}
\begin{align*}
& {\bm r}_2^m \in \arg\max_{{\bm r}_2}\;\{U({\bm r}_2, {\bm v}_2^m):U({\bm r}_2,{\bm c}^m)\leq D\},\\
&({\bm r}_1^m,{\bm r}_3^m, p^m) \in \arg\max_{{\bm r}_1,{\bm r}_3,p\in (0,1)}\;\{(1-p)U({\bm r}_1, {\bm v}_1^m)+pU({\bm r}_3, {\bm v}_1^m):(1-p)U({\bm r}_1, {\bm c}^m)+pU({\bm r}_3, {\bm c}^m) \leq D\}.
\end{align*}
\end{subequations}
With ${\bm Z}_1^m$ and ${\bm Z}_2^m$ constructed 
as above, we ask the DM to
indicate a preference between the two lotteries. Suppose 
that the DM prefers ${\bm Z}_2^m$.
We then
update the polyhedron by
$\widetilde{\cal V}_{m} =  \left\{{\bm v}\in \widetilde{\cal V}_{m-1}: (1-p^m)U({\bm r}_1^m, {\bm v})+pU({\bm r}_3^m, {\bm v}) - U({\bm r}_2^m, {\bm v}) \leq 0\right\}$.
In particular, the above scheme can be applied to 
an \ul{additive multi-attribute PLF} 
$U({\bm x},{\bm v})=\sum_{l=1}^m u_l^{N_l}(x_l,{\bm v}_l)=
\sum_{l=1}^m {\bm v}_l^\top {\bm g}^{{\cal X}_{l}}(x_l)$,
where ${\bm v}_l\in \R^{N_l-1}$ and ${\bm g}^{{\cal X}_l}(x_l)\in \R^{N_l-1}$ is defined as in \eqref{eq:g(t)-PRO_m}, for $l\in [m]$.
When the attribute-wise utility functions are not piecewise linear,
they can be approximated by PLFs of the form $U(\bm x,\bm v)$.
A bi-attribute example illustrating the elicitation procedure is provided as follows.

\begin{example}
We consider a bi-attribute additive nonlinear utility function 
$$
U^*(x_1,x_2)=u_1^*(x_1)+u_2^*(x_2),
$$
where 
$
u_1^*(t)=
\left\{
\begin{array}{ll}
2*(e^{8t}-1)/8, &  \inmat{for} \; 
t\in [-0.5,0]\\
(1-e^{-3t})/3,& 
\inmat{for}\;  t\in [0,0.5]
\end{array}
\right.,
$
and
$
u_2^*(t)=
\left\{
\begin{array}{ll} 
-\lambda*(-x)^a, &  \inmat{for} \; 
t\in [-0.5,0]\\
x^b,& 
\inmat{for}\;  t\in [0,0.5]
\end{array}
\right.,
$
where $a=b=0.58$, $\lambda=2.25$.

We normalize $U^*$ over $[\underline{x}_1,\bar{x}_1]\times [\underline{x}_2,\bar{x}_2]$
as
$
\hat{U}^*(x_1,x_2):=\frac{U^*(x_1,x_2)-U^*(\underline{x}_1,\underline{x}_2)}{U^*(\bar{x}_1,\bar{x}_2)-U^*(\underline{x}_1,\underline{x}_2)},
$
so that
$
\hat{U}^*(\underline{x}_1,\underline{x}_2)=0$ 
and 
$\hat{U}^*(\bar{x}_1,\bar{x}_2)=1$.
Let
$\hat{u}^*_1$ and $\hat{u}_2^*$
denote
the corresponding
marginal utility functions, normalized
such that 
$\hat{u}^*_1(\underline{x}_1)=0$ and  $\hat{u}^*_2(\underline{x}_2)=0$. 
Because of
the additive form of $\hat{U}^*$,
the two marginal utility functions can be represented by a univariate function defined over
$[\bar{x}_1,\bar{x}_1+\bar{x}_2-\underline{x}_2]$,
which is the union of 
$[\underline{x}_1,\bar{x}_1]$ and $[\bar{x}_1+\underline{x}_2,\bar{x}_1+\bar{x}_2]$.
Specifically, define
\bgeqn 
\phi(t):=
\left\{
\begin{array}{ll}
\hat{u}_1^*(t) & \inmat{for} \;\; t\in [\underline{x}_1,\bar{x}_1]\\
\hat{u}_1^*(\bar{x}_1) + \hat{u}_2^*(t-\bar{x}_1+\underline{x}_2)
-\hat{u}_2^*(\underline{x}_2)
& \inmat{for} \;\; t\in [\bar{x}_1,\bar{x}_1+\bar{x}_2-\underline{x}_2]\\
\end{array}
\right..
\edeqn 
The piecewise linear approximation of $\phi(t)$ over $[\bar{x}_1,\bar{x}_1+\bar{x}_2-\underline{x}_2]$ is given by 
\bgeqn
\phi_N(t):=
\left\{
\begin{array}{ll}
\hat{u}_1^{N_1}(t) & \inmat{for} \;\; t\in [\underline{x}_1,\bar{x}_1]\\
\hat{u}_1^{N_1}(\bar{x}_1) + \hat{u}_2^{N_2}(t-\bar{x}_1+\underline{x}_2)
-\hat{u}_2^{N_2}(\underline{x}_2)
& \inmat{for} \;\; t\in [\bar{x}_1,\bar{x}_1+\bar{x}_2-\underline{x}_2]\\
\end{array}
\right..
\edeqn 
The piecewise linear function $\phi_N(t)$ can be represented as
\bgeqn
\phi_N(t)=
\left\{
\begin{array}{ll}
{\bm v}_1^\top{\bm g}^{{\cal X}_1}(t) & \inmat{for} \;\; t\in [\underline{x}_1,\bar{x}_1]\\
{\bm v}_1^\top{\bm g}_1^{{\cal X}_1}(\bar{x}_1) + {\bm v}_2^\top{\bm g}^{{\cal X}_2}(t-\bar{x}_1+\underline{x}_2)
& \inmat{for} \;\; t\in [\bar{x}_1,\bar{x}_1+\bar{x}_2-\underline{x}_2]\\
\end{array}
\right..
\edeqn 

Let $\underline{x}_1=\underline{x}_2=-0.5$ and $\bar{x}_1=\bar{x}_2=0.5.$
Figure~\ref{fig:utiltity_MUl} 
depicts the feasible range of function 
$\phi_N(t)$ over the interval $[\bar{x}_1,\bar{x}_1+\bar{x}_2-\underline{x}_2]$
after $M=10,20,30$ iterations of the MPM.

\begin{figure}[H]
\vspace{-0.2cm}
\begin{minipage}[t]{0.32\linewidth}
\centering
\includegraphics[width=1.9in]{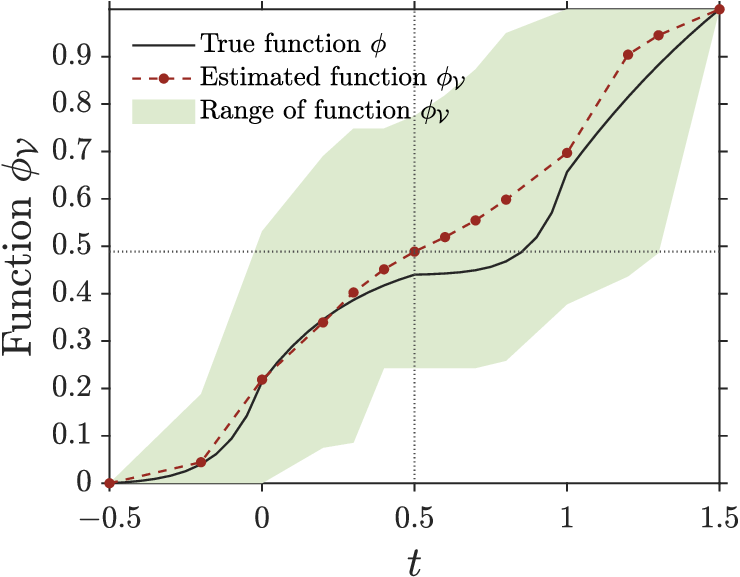}
\text{\footnotesize{\hspace{0.6cm}\textbf{(a)}}}
\end{minipage}
\begin{minipage}[t]{0.32\linewidth}
\centering
\includegraphics[width=1.9in]{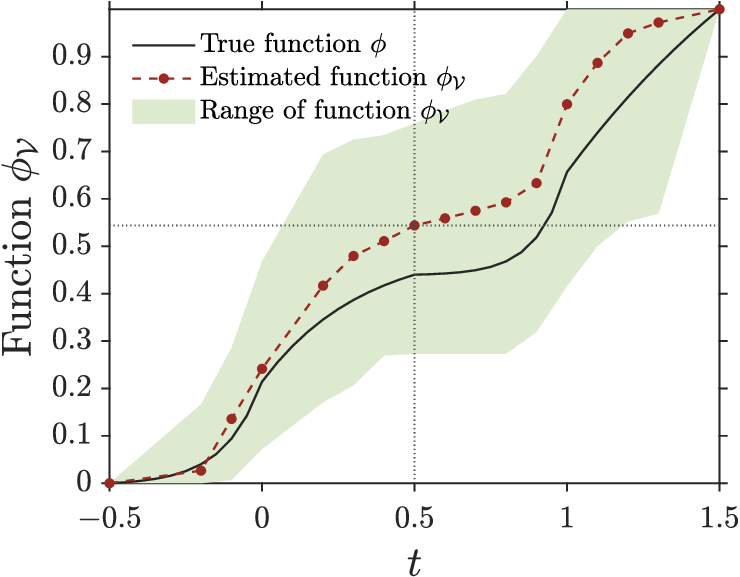}
\text{\footnotesize{\hspace{0.6cm}\textbf{(b)}}}
\end{minipage}
\begin{minipage}[t]{0.32\linewidth}
\centering
\includegraphics[width=1.9in]{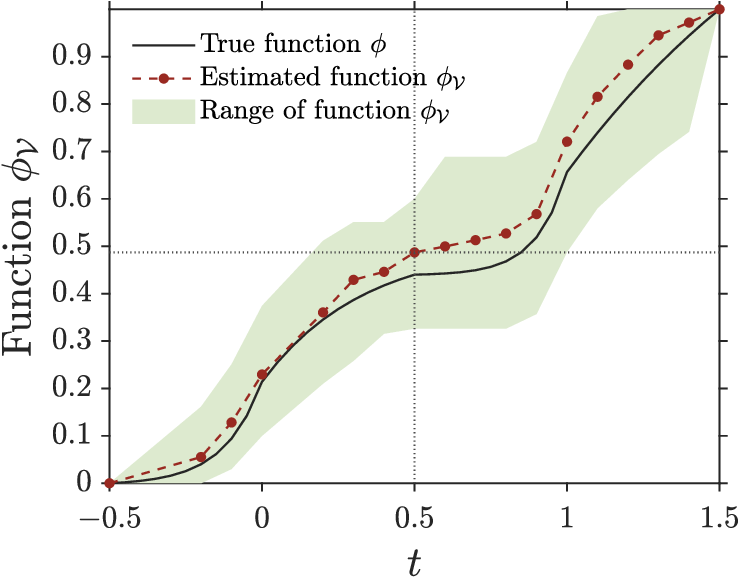}
\text{\footnotesize{\hspace{0.6cm}\textbf{(c)}}}
\end{minipage}
\caption{\footnotesize{
Bi-attribute utility elicitation using the modified polyhedral method:
(a) $M=10$;
(b) $M=20$;
(c) $M=30$.
The breakpoints are rounded to one decimal place.
}}
\vspace{-0.3cm}
\label{fig:utiltity_MUl}
\end{figure}
\end{example}

The scheme can potentially also be extended to the \underline{non-additive nonlinear bi-attribute utility} case. To this end, we can
 use the so-called lambda method of \cite{lee2001polyhedral}
to represent the piecewise linear utility function by a linear combination of the utility values at the breakpoints. 
$ 
U(x,y, {\bm u})=\sum_{i=1}^{N_1}\sum_{j=1}^{N_2} \alpha_{i,j}({x,y}) u_{i,j}= {\bm u}^\top {\bm \alpha}(x,y),
$ 
where 
$
{\bm u}:={\rm vec}\left((u_{ij})_{1\leq i\leq N_1}^{1\leq j\leq N_2}\right)=(u_{1,1},u_{2,1},\cdots,u_{N_1,1},u_{1,2},\\
\cdots,u_{N_1,2},
\cdots,u_{1,N_2},\cdots,u_{N_1,N_2})
$ denotes the vector of utility values at breakpoints $\{(x_i,y_j):i\in [N_1],j\in [N_2]\}$ to be elicited,
and $
{\bm \alpha}({x,y}):={\rm vec}\left((\alpha_{ij}(x,y))_{1\leq i\leq N_1}^{1\leq j\leq N_2}\right)$
is the corresponding weight vector,
see  equations (3.1)-(3.4) in \cite{wu2025bi} for details.

\vspace{-0.1cm}
\subsection{Dynamic Pricing}
\label{sec:dynamic_pricing}
\vspace{-0.1cm}

The basic ideas of the polyhedral method can also be applied to
\underline{dynamic pricing} problems.
\cite{cohen2020feature} consider feature-based dynamic pricing, where the product value is characterized by a multivariate linear function.
Here, 
we extend the setting to the nonlinear case.
Suppose that the market value of product $x_t$ can be represented by an additive utility function of features:
$
v_t=\sum_{j=1}^d u_j(x_{t,j}),
$
where $u_j : [\underline{x}_j, \bar{x}_j] \to [0,1]$ denotes the marginal utility (or market value) of feature $x_{t,j}$, for $j\in [d]$.

Let $u_j^{N_j}$ 
be a piecewise linear utility function 
with breakpoints 
$\mathcal X_j=\{\underline{x}_j=x_1^j<x_2^j<\cdots<
x_{N_j}^j=\bar{x}_j\}$
approximating  $u_j$.
The resulting additive 
piecewise linear utility function
can be written as
$
v_t^N:=\sum_{j=1}^d u_j^{N_j}(x_{t,j}), 
$
where $N=\sum_{j=1}^d N_j$. 
By adopting increment-based PLF representation as in Section \ref{sec:PLA-U}, 
we can write
$
u_j^{N_j}(x_{t,j})
=\bm\alpha_j^\top g^{\mathcal X_j}(x_{t,j}),
$
where ${\bm g}^{\mathcal X_j}(x)=(1,\ldots,({x-x_{l}^j})/({x_{l+1}^j-x_l^j}),0,\ldots,0)$, for $x\in (x_l^j,x_{l+1}^j]$,
and ${\bm \alpha}_j\in \R^{N_j-1}$ is
the vector of increments of $u_j$ over sub-intervals 
$[x_{l}^j, x_{l+1}^j]$, for $l\in [N_j-1]$.
Then the 
function $v_t^N$ can be reformulated:
$
v_t^N={\bm \alpha}^\top {\bm g}(x_t),
$
where 
$
{\bm \alpha}=({\bm \alpha}_1,\ldots,{\bm \alpha}_{d})\in \R^{N_1-1+\ldots+N_d-1}
$
and 
$
{\bm g}({\bm x}_t)=({\bm g}^{\mathcal X_1}(x_{t,1}),\ldots,{\bm g}^{\mathcal X_d}(x_{t,d}))\in \R^{N_1-1+\ldots+N_d-1}.
$

{\bf Pricing Method}:
Using the past price $p^t$ and selling history (product with feature ${\bm x}_t$ is sold out or not),
we propose the following procedure to price a future product ($t'>t$) with feature vector $\bm x_t$.

\begin{itemize}
\item[] {\bf Step 1}. 
Let ${\cal K}_1\subseteq \R^{N-d}$ be  the initial
polyhedron of ${\bm \alpha}$.

\item[]  {\bf Step 2}. {\bf Pricing at period $t$}. Let ${\bm c}^t$ be the analytic center of polyhedron ${\cal K}_{t-1}$. We 
set the price $p^t$ for the product with feature ${\bm x}^t$:  
$p^t=({\bm c}^t)^\top {\bm g}({\bm x}_t)$.

\item[] {\bf Step 3}. {\bf Update of uncertainty set ${\cal K}_{t+1}$}.
If a sale occurs, update the polyhedron as
$
{\cal K}_{t+1}={\cal K}_t\cap \{{\bm \alpha}\in \R^{N-d}: {\bm \alpha}^\top {\bm g}({\bm x}_t) \geq p^t\}.
$
If a sale does not occur, update the polyhedron as
$
{\cal K}_{t+1}={\cal K}_t\cap \{{\bm \alpha}\in \R^{N-d}: {\bm \alpha}^\top {\bm g}({\bm x}_t) \leq p^t\}.
$
\end{itemize}
The process learns customer's preferences,
approximately represented by a piecewise linear utility function, from the dynamic sales data of the products.

\vspace{-0.2cm}
\section{Concluding Remarks}
\label{sec:Conclusion}
\vspace{-0.2cm}

In this paper, we propose a MPM for eliciting a nonlinear univariate utility function. The method constructs a set of piecewise linear utility functions to approximate the unknown true utility function and iteratively shrinks this set based on the DM's responses to adaptively designed pairwise comparison queries.
We then apply the resulting
utility ambiguity set to PRO model and develop
two schemes for controlling the conservatism of the corresponding
maximin problem. Numerical experiments demonstrate the performance
of the proposed elicitation method and the effects of the
conservatism reduction schemes.

There are several promising directions for future research. First, it would be valuable to establish the theoretical complexity of the proposed approach. Although our computational results provide empirical evidence of convergence, and we conjecture that the approach converges asymptotically as the number of elicitation queries tends to infinity, characterizing and quantifying its rate of convergence remains a challenging problem.
Second, response errors may occur during the elicitation process. In such cases, it may be worthwhile to explore a probabilistically MPM.
Third, while the approach can readily be extended to additive multivariate nonlinear utility functions, there is also the potential to extend it to the elicitation of non-additive multivariate utility functions as outlined in Section \ref{sec:multi}. 
Fourth, it would be interesting to adapt the approach to the elicitation of risk measures, such as distortion and spectral risk measures. We leave these directions for future research.

\bibliographystyle{plain}
\bibliography{bibfile_arxiv}

\newpage 
\appendix 
\section*{Appendix}

\section{Mixed-Integer Linear Reformulation of
\texorpdfstring{$\bm{g}(x)$ in \eqref{eq:g(x)-incre}}{g(x)}}
The mapping ${\bm g}(\cdot)$ can instead be represented through a system of linear equalities and inequalities involving auxiliary variables.

\begin{proposition}
[Linear Representation
of ${\bm g}(x)$, \cite{HZXZ22}]
\label{prop:g(x)-rfmlt}
The vector-valued mapping ${\bm g}(\cdot)$ can be 
represented by 
${\bm g}(x)={\bm P}{\bdy}+{\bm Q}{\bm z}$,
where the constant matrices
${\bm P},{\bm Q}\in \R^{(N-1)\times (N-1)}$ are defined as 
\bgeqn
\label{eq:definition-AB}	{\bm P}:= 
\left[
\begin{array}{cccc}
\frac{1}{x_{2} - x_{1}} & & & \\
&\frac{1}{x_{3} - x_{2}} & & \\
& & \ddots & \\
& & & \frac{1}{x_{N} - x_{N-1}}
\end{array}
\right]
\quad \inmat{and}\quad
{\bm Q}:= \left[
\begin{array}{cccc}
\frac{-x_{1}}{x_{2}-x_{1}} & 1 & \dots & 1 \\
& \frac{-x_{2}}{x_{3}-x_{2}} & \dots & 1 \\
&  & \ddots & \vdots \\
&  &  & \frac{-x_{N-1}}{x_{N}-x_{N-1}}
\end{array}
\right].~
\edeqn
For any fixed $x$, the auxiliary variables ${\bdy},{\bm z}$ are uniquely determined by the following mixed-integer linear system:
\bgeq
\left\{({\bm y},{\bm z})\in \R^{N-1}\times\{0,1\}^{N-1}:
\begin{array}{ll}
{\bm e}^\top {\bm z}=1,\;
{\bm e}^\top {\bm y}= x,\;\\
x_{i}z_i-y_i\leq 0,\;
x_{i+1}z_i-y_i\geq 0,\;
z_i\in \{0,1\},\; i\in [N-1]
\end{array}
\right\}.
\edeq
\end{proposition}

\section{Properties of Problems \eqref{eq:gene_m_Q-1-v1-p}-\eqref{eq:gene_m_Q-2-v1-p}}
\label{app_sec:properties}

The main purpose of 
this subsection is to show that
the optimum 
of the programs in 
\eqref{eq:gene_m_Q}
are attained at solutions satisfying the corresponding
equality constraints. 
To this end, we introduce 
notation for the feasible sets, the optimal values, the set of optimal solutions,
and the sets of feasible solutions that satisfy the equality constraints
in problems \eqref{eq:gene_m_Q-1-v1-p}-\eqref{eq:gene_m_Q-2-v1-p}.
\begin{subequations}
\begin{align}
{\cal F}_{{\bm B}}:=&\left\{{r}_2\in [\underline{x},\bar{x}]:({\bm R}{\bm c}^m)^{\top}G_{\bm B}({r}_2) \leq D \right\},\;\;
{\cal F}_{{\bm A}}:=\left\{(r_1,r_3,p)\in F,
({\bm R}{\bm c}^m)^{\top}G_{\bm A}({r}_1,{r}_3,{p})\leq D\right\}, \label{eq:feasible_2}\\
\vt_{\bm B}:=&\max\left\{({\bm R}{\bm v}_1^m)^{\top}G_{\bm B}(r_2): r_2\in {\cal F}_{\bm B}\right\},\;\;
\vt_{\bm A}:=\max\left\{({\bm R}{\bm v}_2^m)^{\top}G_{\bm A}(r_1,r_3,p): (r_1,r_3,p)\in {\cal F}_{\bm A}\right\},\\
{\cal S}_{\bm B}:=&\left\{r_2\in {\cal F}_{\bm B}: \vt_{\bm B}=({\bm R}{\bm v}_1^m)^{\top}\hspace{-0.1cm}G_{\bm B}(r_2)\right\},\;\;
{\cal S}_{\bm A}:=\{(r_1,r_3,p)\in {\cal F}_{\bm A}: \vt_{\bm A}=({\bm R}{\bm v}_2^m)^{\top}\hspace{-0.1cm}G_{\bm A}(r_1,r_3,p)\}, \label{eq:solution_B}\\
{\cal F}_{\bm B}^{=}:=&\left\{{r}_2\in [\underline{x},\bar{x}]:
({\bm R}{\bm c}^m)^{\top}G_{\bm B}({r}_2) =D
\right\},\;\;
 {\cal F}_{\bm A}^{=}:=\left\{(r_1,r_3,p)\in F:
({\bm R}{\bm c}^m)^{\top}G_{\bm A}({r}_1,{r}_3,{p})=D\right\}.
\label{eq:bar_r123-2}
\end{align}
\end{subequations}
Problems \eqref{eq:gene_m_Q-1-v1-p}-\eqref{eq:gene_m_Q-2-v1-p} will be used to generate the parameters $(r_1^m,r_2^m,r_3^m,p^m)$ 
for constructing
the $m$-th pair of lotteries ${\bm A}_m$ and ${\bm B}_m$ defined as in \eqref{eq:AmBm}.
The following proposition characterizes the optimal solutions of these two problems.

\begin{proposition}
[Properties of ${\cal F}_{\bm B}^=$, ${\cal F}_{\bm A}^=$ and $G_{\bm A}(r_1,r_3,p)$]
\label{prop:EDn=0_1}
Let $D\in [0,1]$ be a constant,
and let
${\bm c}^m$ be the analytic center of the polyhedron ${\cal V}_N^{m-1}$ defined as in \eqref{eq:ana_center}.
Then
$|{\cal F}_{\bm B}^{=}|=1$,
${\cal F}_{\bm A}^{=}$
is uncountable,
and
\bgeqn
\label{eq:vt_2}
G_{\bm A}(r_1,r_3,p)
=\Bigg(1,\cdots,{1},\overbrace{(1-{p})
\frac{r_1-x_{j_{r_1}}}{x_{j_{r_1}+1}-x_{j_{r_1}}}
+{p}}^{j_{{r}_1}},{p},\cdots,{p},\overbrace{{p}
\frac{r_3-x_{j_{r_3}}}{x_{j_{r_3}+1}-x_{j_{r_3}}}}^{j_{{r}_3}},0,\cdots,0
\Bigg),
\edeqn
where $|S|$ denotes the cardinality of $S$,
for $i=1,2,3$,
$j_{r_i}$ denotes
the index such that
${r}_i\in (x_{j_{{r}_i}},x_{j_{{r}_i}+1}]$ with $j_{r_1}=1$ when $r_1=\underline{x}$.
\end{proposition}

\begin{proof}{Proof.}
Let 
\bgeqn
\label{eq:UpsilonAB}
\Upsilon_{\bm A}^{=} := \{{\bm d}\in \Upsilon_{\bm A}: 
({\bm R}{\bm c}^m)^\top{\bm d}=D\}
\quad\inmat{and} \quad
\Upsilon_{\bm B}^{=} := \{{\bm d}\in \Upsilon_{\bm B}: 
({\bm R}{\bm c}^m)^\top{\bm d}=D\}.
\edeqn
Then
${\cal F}_{\bm A}^{=} = G_{\bm A}^{-1}(\Upsilon_{\bm A}^{=})$ and ${\cal F}_{\bm B}^{=} = G_{\bm B}^{-1}(\Upsilon_{\bm B}^{=})$,
where $G_{\bm A}^{-1}$ and $G_{\bm B}^{-1}$ are defined as in \eqref{eq:GAB-1}.

\underline{We begin by proving that $|{\cal F}_{\bm B}^{=}|=1$}.
Since $G_{\bm B}^{-1}$ is one-to-one and single-valued, it suffices to show that $|\Upsilon_{\bm B}^=|=1$. 
We consider two cases.
\underline{Case 1. $D=1
$}.
Observe from \eqref{eq:cm} that  
${\bm c}^m>0$, ${\bm e}^{\top}{\bm c}^m<1$, and ${\bm R}{\bm c}^m=({\bm c}^m,1-{\bm e}^\top {\bm c}^m)>0$.
Since ${\bm d}$ lies in
the range of ${\bm g}(\cdot)$,
which is contained in $[0,1]^{N-1}$,
and ${\bm g}(x)$ 
has the form given in
\eqref{eq:g(x)-incre},
the only ${\bm d}$
satisfying $({\bm R}{\bm c}^m)^\top{\bm d}=1$ is
${\bm d}={\bm e}$.
\underline{Case 2. $D\in (0,1)$}.
Since $({\bm R}{\bm c}^m)^\top{\bm e}=1$, then
there exists a unique $i_D\in [N-1]$ such that 
\bgeqn
\label{eq:sumc<=D}
\sum_{j=1}^{i_D-1}c_j^m\leq  D<\sum_{j=1}^{i_D-1}c_j^m+c_{i_D}^m,
\edeqn
which 
implies that
$0\leq \frac{D-\sum_{j=1}^{i_D-1}c_j^m}{c_{i_D}^m}<1$.
Let
$$
{\bm d}_{D} :=\left(1,\cdots,1,\frac{D-\sum_{j=1}^{i_D-1}c_j^m}{c_{i_D}^m},0,\cdots,0\right).
$$
Then
${\bm d}_D\in \Upsilon_{\bm B}$ and $({\bm R}{\bm c}^m)^\top {\bm d}_D=\sum_{j=1}^{i_D-1}c_j^m+\frac{D-\sum_{j=1}^{i_D-1}c_j^m}{c_{i_D}^m} c_{i_D}^m=D$,
and hence
${\bm d}_D\in \Upsilon_{\bm B}^=$.
It remains to
show that ${\bm d}_D$ is the 
unique
solution to $({\bm R}{\bm c}^m)^\top {\bm d}_D=D$.
Suppose, for contradiction,
that there exists  $\tilde{\bm d}_{D}\in \Upsilon_{\bm B}\backslash{\{{\bm 0},{\bm e}\}}$,
with $\tilde{\bm d}_{D}\neq {\bm d}_{D}$, such that
\bgeqn
({\bm R}{\bm c}^m)^\top \tilde{\bm d}_D=\sum_{j=1}^{\tilde{i}_D-1}c_j^m+\tilde{d}_{\tilde{i}_D}{c_{\tilde{i}_D}^m}=D.
\label{eq:d'=D}
\edeqn
Then $\tilde{\bm d}_{D}$ must take the form
$\tilde{\bm d}_{D}=(1,\cdots,1,\tilde{d}_{\tilde{i}_D},0,\cdots,0)$ with some $\tilde{d}_{\tilde{i}_D}\in [0,1)$.
If $\tilde{i}_D\leq i_D-1$,
then 
$$
D\overset{\eqref{eq:d'=D}}{=}\sum_{j=1}^{\tilde{i}_D-1}c_j^m+\tilde{d}_{\tilde{i}_D}{c_{\tilde{i}_D}^m}\overset{\tilde{d}_{\tilde{i}_D}<1}{<}\sum_{j=1}^{\tilde{i}_D}c_j^m \leq \sum_{j=1}^{i_D-1}c_j^m \overset{\eqref{eq:sumc<=D}}{\leq} D,
$$
which yields a contradiction.
Likewise, if $\tilde{i}_D\geq i_D+1$,
then 
$$
D\overset{\eqref{eq:d'=D}}{=}\sum_{j=1}^{\tilde{i}_D-1}c_j^m+\tilde{d}_{\tilde{i}_D}{c_{\tilde{i}_D}^m}\overset{\tilde{d}_{\tilde{i}_D}\geq 0}{\geq } \sum_{j=1}^{\tilde{i}_D-1}c_j^m\geq \sum_{j=1}^{i_D}c_j^m\overset{\eqref{eq:sumc<=D}}{>} D,
$$
again yielding a contradiction.
Therefore,
$\tilde{i}_D=i_D$. 
It then follows from \eqref{eq:d'=D} that 
$D=\sum_{j=1}^{i_D-1}c_j^m+\tilde{d}_{i_D}{c_{i_D}^m}$ and thus $\tilde{d}_{i_D}=\frac{D-\sum_{j=1}^{i_D-1}c_j^m}{c_{i_D}^m}$, 
implying
${\bm d}_{D}= \tilde{\bm d}_D$, a contradiction!

\underline{Next, we show that
${\cal F}_{\bm A}^{=}$ is uncountable}.
For any $D\in(0,1)$ and any $D_-\in(0,D)$, there exists a unique
${\bm d}_{D_-}=(1,\ldots,1,d_{D_-,i_-},0,\ldots,0)$
in the range of ${\bm g}(\cdot)$ such that
\[
(R{\bm c}^m)^\top {\bm d}_{D_-}=D_-.
\]
Similarly, for any $D_+\in(D,1]$, there exists a unique
${\bm d}_{D_+}=(1,\ldots,1,d_{D_+,i_+},0,\ldots,0)$
in the range of ${\bm g}(\cdot)$ such that
\[
(R{\bm c}^m)^\top {\bm d}_{D_+}=D_+.
\]
Let $p:=\frac{D-D_-}{D_+-D_-}$ and ${\bm d}_{D}:=(1-p){\bm d}_{D_-}+p{\bm d}_{D_+}\in \Upsilon_{\bm A}$.
Let ${\cal D}$ denote
the set of all such ${\bm d}_D$ 
obtained by varying $D_-\in (0,D)$ and $D_+\in (D,1]$.
Then
$({\bm R}{\bm c}^m)^\top {\bm d}_{D}=(1-p)D_-+p D_+=D$ for
every
${\bm d}_D\in {\cal D}$.
Thus,
${\cal D}\subset \Upsilon_{\bm A}^=$.
Since ${\cal D}$ is uncountable and $G^{-1}_{\bm A}$ is set-valued,
it follows that
${\cal F}_{\bm A}^{=}=G_{\bm A}^{-1}(\Upsilon_{\bm A}^=)\supset G_{\bm A}^{-1}({\cal D})$ is uncountable.

\underline{Finally, we establish \eqref{eq:vt_2}.}
By the definitions of $G_{\bm A}(\cdot,\cdot,\cdot)$ in \eqref{eq:def_G} and ${\bm g}(\cdot)$ in \eqref{eq:g(x)-incre},
for any $(r_1,r_3,p)\in {\cal F}_{\bm A}$,
we have
\bgeqn
G_{\bm A}(r_1,r_3,p)
&=&(1-p)\Big(1,\cdots,{1},\overbrace{
\frac{r_1-x_{j_{r_1}}}{x_{j_{r_1}+1}-x_{j_{r_1}}}}^{j_{{r}_1}},0,\cdots,0
\Big)
+p\Big(1,\cdots,{1},\overbrace{
\frac{r_3-x_{j_{r_3}}}{x_{j_{r_3}+1}-x_{j_{r_3}}}}^{j_{{r}_3}},0,\cdots,0
\Big) \nonumber \\
&=&\Big(1-p,\cdots,{1-p},\overbrace{
(1-p)\frac{r_1-x_{j_{r_1}}}{x_{j_{r_1}+1}-x_{j_{r_1}}}}^{j_{{r}_1}},0,\cdots,0
\Big)
+\Big(p,\cdots,{p},\overbrace{
p\frac{r_3-x_{j_{r_3}}}{x_{j_{r_3}+1}-x_{j_{r_3}}}}^{j_{{r}_3}},0,\cdots,0
\Big) \nonumber \\
&=&\Big(1,\cdots,{1},\overbrace{(1-{p})
\frac{r_1-x_{j_{r_1}}}{x_{j_{r_1}+1}-x_{j_{r_1}}}
+{p}}^{j_{{r}_1}},{p},\cdots,{p},\overbrace{{p}
\frac{r_3-x_{j_{r_3}}}{x_{j_{r_3}+1}-x_{j_{r_3}}}}^{j_{{r}_3}},0,\cdots,0
\Big).
\label{eq:G_A}
\edeqn
The proof is complete.
\hfill 
$\Box$
\end{proof}

Proposition~\ref{prop:EDn=0_1} characterizes the cardinality of ${\cal F}_{\bm B}^{=}$,
${\cal F}_{\bm A}^{=}$,
and provides an explicit representation of the vector-valued function $G_{\bm A}(\cdot,\cdot,\cdot)$,
which plays an important role in
analyzing the properties of problem \eqref{eq:gene_m_Q}.
We next show that, for any choice of $D$,
the optimal values of problems \eqref{eq:gene_m_Q-1-v1-p}-\eqref{eq:gene_m_Q-2-v1-p}
are attained at feasible solutions 
that satisfy the corresponding equality constraints.

\begin{proposition}
[Properties of the Optimal Values of Problems \eqref{eq:gene_m_Q-1-v1-p}-\eqref{eq:gene_m_Q-2-v1-p}]
\label{prop:EDn=0_2-1st}
Let 
${\bm v}_1^m$
and ${\bm v}_2^m$
be defined as in problems \eqref{eq:gene_m_Q-1-v1-p}-\eqref{eq:gene_m_Q-2-v1-p},
respectively, and let
${\bm c}^m$ be the analytic center of polyhedron ${\cal V}_N^{m-1}$ defined as in \eqref{eq:ana_center}.
Then
\begin{subequations}
\begin{align}
\max_{r_2\in {\cal F}_{\bm B} }\;({\bm R}{\bm v}_1^m)^{\top}G_{\bm B}(r_2)
=& \max_{r_2\in {\cal F}_{\bm B}^{=} }\;({\bm R}{\bm v}_1^m)^{\top}G_{\bm B}(r_2),
\label{eq:F_B=F_B=}
\\
\max_{(r_1,r_3,p) \in {\cal F}_{\bm A} } ({\bm R}{\bm v}_2^m)^{\top}G_{\bm A}(r_1,r_3,p)
=&\max_{(r_1,r_3,p) \in {\cal F}_{\bm A}^{=} } ({\bm R}{\bm v}_2^m)^{\top}G_{\bm A}(r_1,r_3,p),
\label{eq:F_A=F_A=}
\end{align}
\end{subequations}
where ${\cal F}_{\bm A}$, ${\cal F}_{\bm B}$, ${\cal F}_{\bm A}^=$, and ${\cal F}_{\bm B}^=$ are defined as in \eqref{eq:feasible_2} and \eqref{eq:bar_r123-2}.
\end{proposition}

\proof{Proof of Proposition~\ref{prop:EDn=0_2-1st}.}
\underline{\textbf{First}, we prove \eqref{eq:F_B=F_B=}}.
Since ${\cal F}_{\bm B}\supset {\cal F}_{\bm B}^=$,  
$$
\inmat{the left-hand-side (lhs) of \eqref{eq:F_B=F_B=}} 
\geq
\inmat{the right-hand-side (rhs) of \eqref{eq:F_B=F_B=}}.
$$
Suppose, for contradiction,
that the equality does not hold. Then 
$$
\inmat{the lhs of \eqref{eq:F_B=F_B=}}  
> \inmat{the rhs of \eqref{eq:F_B=F_B=}},
$$
implying that
${\cal S}_{\bm B}\subset ({\cal F}_{\bm B}\backslash{{\cal F}_{\bm B}^=})$. By the definition of $\Upsilon_{\bm B}$ 
in \eqref{eq:Upsilon_B},
the two sides of \eqref{eq:F_B=F_B=} can be written as
\begin{subequations} 
\begin{align}
\inmat{lhs of \eqref{eq:F_B=F_B=} }
&= \max_{{\bm d}\in \Upsilon_{\bm B} }\{({\bm R}{\bm v}_1^m)^{\top}{\bm d}:({\bm R}{\bm c}^m)^{\top}{\bm d}\leq D\},
\label{eq:F_B=F_B=-lhs}
\\
\inmat{rhs of \eqref{eq:F_B=F_B=} }
&= \max_{{\bm d}\in \Upsilon_{\bm B}^= }\;({\bm R}{\bm v}_1^m)^{\top}{\bm d}.
\label{eq:F_B=F_B=-rhs}
\end{align}
\end{subequations}
\underline{Case 1: $D=1$}.
It follows
from the proof of Proposition~\ref{prop:EDn=0_1} that
$\Upsilon_B^{=}=\{{\bm e}\}$.
Since 
$
({\bm R}{\bm v}_1^m)^{\top}{\bm d} \leq 
({\bm R}{\bm v}_1^m)^{\top}{\bm e},
$ 
for every
${\bm d}\in \Upsilon_{\bm B}$ satisfying  
$({\bm R}{\bm c}^m)^\top {\bm d}\leq 1$,
${\bm e}$ is also an optimal solution of the rhs of \eqref{eq:F_B=F_B=-lhs}.
Hence, the optimal values in \eqref{eq:F_B=F_B=-lhs} and \eqref{eq:F_B=F_B=-rhs} coincide, contradicting our assumption.

\noindent\underline{Case 2: $D\in (0,1)$}.
Define
$\Upsilon_{\bm B}^{<}:=\{{\bm d}\in \Upsilon_{\bm B}:({\bm R}{\bm c}^m)^\top {\bm d}<D\}$.
Then $\{{\bm d}\in \Upsilon_{\bm B}: ({\bm R}{\bm c}^m)^\top {\bm d} \leq D\}= \Upsilon_{\bm B}^{<}\bigcup \Upsilon_{\bm B}^{=}$.
Let ${\bm d}^B$ be an optimal solution to the rhs of \eqref{eq:F_B=F_B=-lhs}. 
Under our contradiction hypothesis,
${\bm d}^B\in \Upsilon_{\bm B}^{<}$.
Hence,
there exists $i^B\in [N-1]$ such that ${\bm d}^{B}=(1,\cdots,1,d_{i^B}^B,0,\cdots,0)$,
where
$0\leq d_{i^B}^B<1$ and $
({\bm R}{\bm c}^m)^{\top}{\bm d}^B<D
$.
Then 
\bgeqn
\label{eq:d_iB}
d_{i^B}^B<\frac{D-\sum_{j=1}^{i^B-1}c_j^m}{c_{i^B}^m}.
\edeqn
To derive a contradiction, we construct a feasible solution
$\tilde{\bm d}^B \in \Upsilon_{\bm B}^=$
whose
objective value
is at least as large as that attained by the optimal solution
${\bm d}^B\in \Upsilon_{\bm B}^{<}$. 
Let 
\bgeqn
\label{eq:tilded^B}
\tilde{d}_{i^B+{l}}^B:=\min\left\{\max\left\{\frac{D-\sum_{j=1}^{i^B-1+{l}}c_j^m}{c_{i^B+{l}}^m},0\right\},1\right\},
\inmat{ for }{l}\in \{0\}\cup [N-1-i^B],
\edeqn
and 
$\tilde{\bm d}^B:=(1,\cdots,1,\tilde{d}_{i^B}^B,\tilde{d}_{i^B+1}^B,\cdots,\tilde{d}_{N-1}^B)$.
By construction,
$$
\tilde{d}_{i^B+{l}}^B 
=
\left\{\begin{array}{lll}
\min\left\{\frac{D-\sum_{j=1}^{i^B-1+{l}}c_j^m}{c_{i^B+{l}}^m},1\right\}, &&  \sum_{j=1}^{i^B-1+{l}}c_j^m < D, \\
0,  && \inmat{otherwise}.
\end{array}
\right.
$$
Then  $\tilde{\bm d}^B\in \Upsilon_{\bm B}$ and 
$({\bm R}{\bm c}^m)^{\top}\tilde{\bm d}^B=D$,
and hence
$\tilde{\bm d}^B\in \Upsilon_{\bm B}^=$.
Let $\bar{l}^B\in [N-2-i^B]$ 
denote the largest index
such that $\tilde{d}^B_{{i}^B+\bar{l}^B}>0$, i.e.,
\bgeqn
\label{eq:lable_lB}
\bar{l}^B:=\max\left\{l:  \tilde{d}^B_{{i}^B+l}>0, l \in [N-2-i^B]\right\}.
\edeqn
\underline{If $v_{1,i^B}^m>0$},
then the objective value  of 
the rhs of \eqref{eq:F_B=F_B=-lhs} attained by
$\tilde{\bm d}^B$
satisfies
\bgeqn
\label{eq:v_i*>0}
({\bm R}{\bm v}_1^m)^\top \tilde{\bm d}^B&=&\sum_{j=1}^{i^B-1}v_{1,j}^m+\tilde{d}_{i^B}^Bv_{1,i^B}^m+\tilde{d}_{i^B+1}^Bv_{1,i^B+1}^m+\cdots+ \tilde{d}_{N-1}^Bv_{1,N-1}^m \nonumber \\
&
\geq&\sum_{j=1}^{i^B-1}v_{1,j}^m+\min\left\{\frac{D-\sum_{j=1}^{i^B-1}c_j^m}{c_{i^B}^m},1\right\}v_{1,i^B}^m\\
&\overset{\eqref{eq:d_iB}}{>}&\sum_{j=1}^{i^B-1}v_{1,j}^m+d_{i^B}^Bv_{1,i^B}^m=\vt_{\bm B},\nonumber 
\edeqn
which contradicts 
the fact that $\vt_{\bm B}$ is the optimal value of the rhs of \eqref{eq:F_B=F_B=-lhs}. 

\noindent\underline{If $v_{1,i^B}^m=0$ and $v_{1,i^B+l^B}^m>0$} for some $l^B\in [\bar{l}^B]$,
then
\bgeq
({\bm R}{\bm v}_1^m)^\top \tilde{\bm d}^B
\geq \sum_{j=1}^{i^B-1}v_{1,j}^m+
\tilde{d}^B_{i^B+l^B} v_{1,i^B+l^B}^m \overset{\eqref{eq:lable_lB}}{>} \sum_{j=1}^{i^B-1}v_{1,j}^m+0
=\sum_{j=1}^{i^B-1}v_{1,j}^m+d_{i^B}^B\cdot 0=\vt_{\bm B},
\edeq
again contradicting
the fact that $\vt_{\bm B}$ is the optimal value of the rhs of \eqref{eq:F_B=F_B=-lhs}.

\noindent\underline{If $v_{1,i^B+l}^m=0$ for all $l\in \{0\}\cup [\bar{l}^B]$},
then 
$$
({\bm R}{\bm v}_1^m)^\top \tilde{\bm d}^B
=\sum_{j=1}^{i^B-1}v_{1,j}^m+0=\sum_{j=1}^{i^B-1}v_{1,j}^m+d_{i^B}^B\cdot 0=\vt_{\bm B},
$$
which implies that
$\tilde{\bm d}^B\in \Upsilon_{\bm B}^=$ 
attains
the optimal value $\vt_{\bm B}$.
This contradicts our
assumption that any optimal 
solutions lie in
$\Upsilon_{\bm B}^{<}$.

\underline{\textbf{Next}, we prove \eqref{eq:F_A=F_A=}}. The proof is
like
that of \eqref{eq:F_B=F_B=},
but we provide the details for completeness.
Since ${\cal F}_{\bm A}\supset {\cal F}_{\bm A}^=$, 
$$
\inmat{the lhs of \eqref{eq:F_A=F_A=}} 
\geq
\inmat{the rhs of \eqref{eq:F_A=F_A=}}.
$$
Suppose, for contradiction, that the equality in \eqref{eq:F_A=F_A=} does not hold.
Then 
$$
\inmat{the lhs of \eqref{eq:F_A=F_A=}}
>
\inmat{the rhs of \eqref{eq:F_A=F_A=}},
$$
which implies that
${\cal S}_{\bm A}\subset ({\cal F}_{\bm A}\backslash{{\cal F}_{\bm A}^=})$. By the definition of $\Upsilon_{\bm A}$ in \eqref{eq:Upsilon_A},
the lhs and rhs of \eqref{eq:F_A=F_A=} can be written as
\begin{subequations}
\begin{align}
\inmat{lhs of \eqref{eq:F_A=F_A=} }
&= \max_{{\bm d}\in \Upsilon_{\bm A} }\{({\bm R}{\bm v}_2^m)^{\top}{\bm d}:({\bm R}{\bm c}^m)^{\top}{\bm d}\leq D\},
\label{eq:F_A=F_A=-lhs}
\\
\inmat{rhs of \eqref{eq:F_A=F_A=} }
&= \max_{{\bm d}\in \Upsilon_{\bm A}^= }\;({\bm R}{\bm v}_2^m)^{\top}{\bm d}.
\label{eq:F_A=F_A=-rhs}
\end{align}
\end{subequations}
\underline{Case 1: $D=1$}.
In this case, 
$\Upsilon_{\bm A}^{=}=\{{\bm e}\}$.
As shown earlier,
${\bm e}$ is an optimal solution of the rhs of \eqref{eq:F_A=F_A=-lhs}. 
Therefore, the lhs and rhs of \eqref{eq:F_A=F_A=-lhs} are equal, contradicting our assumption.

\noindent\underline{Case 2: $D\in (0,1)$}.
Let $\Upsilon_{\bm A}^{<}:=\{{\bm d}\in \Upsilon_{\bm A}:({\bm R}{\bm c}^m)^\top {\bm d}<D\}$.
Then $ \{{\bm d}\in \Upsilon_{\bm A}: ({\bm R}{\bm c}^m)^\top {\bm d}\leq D\}= \Upsilon_{\bm A}^{<}\bigcup \Upsilon_{\bm A}^{=}$.
Under our contradiction hypothesis,
${\bm d}^A\in \Upsilon_{\bm A}^{<}$.

By the definition of $\Upsilon_{\bm A}$ in \eqref{eq:Upsilon_A},
there exist $p\in (0,1)$, and ${\bm d}_1$, ${\bm d}_3$ in the range of ${\bm g}(\cdot)$ of the form
$$
{\bm d}_1=(1,\cdots,1,\overbrace{d_{1,i_1}}^{i_1},0,\cdots,0)
\quad\inmat{and}\quad 
{\bm d}_3=(1,\cdots,1,\overbrace{d_{3,i_3}}^{i_3},0,\cdots,0)
$$ 
where
$0\leq d_{1,i_1}<1$, $0\leq d_{3,i_3}<1$ and $i_1\leq i_3$ such that
\bgeqn
\label{eq:dA}
{\bm d}^A
=(1-p){\bm d}_1+p{\bm d}_3=\Big(1,\cdots,{1},\overbrace{(1-p)d_{1,i_1}+p}^{i_1},{p},\cdots,{p},\overbrace{{p}
{d}_{3,i_3}}^{i_3},0,\cdots,0
\Big).
\edeqn
Since
${\bm d}^A\in \Upsilon_{\bm A}^{<}$, we have
\bgeqn
\label{eq:T123}
({\bm R}{\bm c}^m)^\top {\bm d}^A=\underbrace{\sum_{j=1}^{i_1-1}c_j^m}_{T_1}+\underbrace{(1-p)d_{1,i_1}c_{i_1}^m}_{T_2}+\underbrace{p\sum_{j=i_1}^{i_3-1} c_j^m}_{T_3} +{p}
{d}_{3,i_3} c_{i_3}^m <D,
\edeqn
which yields that 
${d}_{3,i_3} < \frac{D-(T_1+T_2+T_3)}{pc_{i_3}^m}$.
We then define
$$
\tilde{d}_{3,i_3+l}:=\min\left\{\max\left\{\frac{D-(T_1+T_2+T_3)-p\sum_{j=i_3}^{i_3-1+l}c_j^m}{pc_{i_3+l}^m},0\right\},1\right\}, \inmat{ for }l\in \{0\}\cup [N-1-i_3].
$$
By convention, we set
$\sum_{j=i}^{i-1} c_j^m=0$.
By construction,
\bgeqn
\label{eq:defi_tilded3}
\tilde{d}_{3,i_3+l} 
=
\left\{\begin{array}{lll}
\min\left\{\frac{D-(T_1+T_2+T_3)-p\sum_{j=i_3}^{i_3-1+l}c_j^m}{pc_{i_3+l}^m},1\right\},
&& (T_1+T_2+T_3)+p\sum_{j=i_3}^{i_3+l-1}c_j^m< D, \\
0,  
&& \inmat{otherwise}.
\end{array}
\right.
\edeqn
Define
 \bgeqn
 \label{eq:label_lA}
 \bar{l}^A:=\max\left\{l: \tilde{d}_{3,{i}_3+{l}}>0,l\in [N-2-i_3]\right\}.
  \edeqn
We first show that
\bgeqn
\label{eq:tilded_3d_3}
\tilde{d}_{3,j}>{d}_{3,j}, \;\; j\in [i_3+\bar{l}^A].
\edeqn
For $j=i_3$, we have
$$
\frac{D-(T_1+T_2+T_3)-p\sum_{j=i_3}^{i_3-1}c_j^m}{pc_{i_3}^m}= \frac{D-(T_1+T_2+T_3)}{pc_{i_3}^m}\overset{\eqref{eq:T123}}{>}
{d}_{3,i_3} \quad \inmat{and}\quad 1>{d}_{3,i_3}.
$$
Hence,
$\tilde{d}_{3,i_3}=\min\left\{\frac{D-(T_1+T_2+T_3)-p\sum_{j=i_3}^{i_3-1}c_j^m}{pc_{i_3}^m},1\right\}>d_{3,i_3}$.
For $l\in [\bar{l}^A]$,
the definition of $ \bar{l}^A$ in \eqref{eq:label_lA}
implies that
$\tilde{d}_{3,i_3+{l}}>0=d_{3,i_3+l}$, for $l\in [\bar{l}^A]$. 
Let 
$\tilde{\bm d}_3:=(1,\cdots,1,\tilde{d}_{3,i_3},\tilde{d}_{3,i_3+1},\cdots,\tilde{d}_{3,N-1})$.
We next show that $\tilde{\bm d}_3$ lies in the range of ${\bm g}(\cdot)$.
First, suppose that
$\tilde{d}_{3,N-1}>0$.
Then
$$
(T_1+T_2+T_3)+p\sum_{j=i_3}^{N-2}c_j^m<D,
$$
which implies that
$$
\frac{D-(T_1+T_2+T_3)-p\sum_{j=i_3}^{N-3}c_j^m}{pc_{N-2}^m}>1,
$$
and hence 
$$
\tilde{d}_{3,N-3} =
\min\left\{\frac{D-(T_1+T_2+T_3)-p\sum_{j=i_3}^{N-3}c_j^m}{pc_{N-2}^m},1\right\}=1.
$$
Applying the same argument recursively yields
$\tilde{d}_{3,j}=1$ for $j=i_3,\cdots,N-3$.
Next, to construct
$\tilde{\bm d}_1$,
define
$$
T_4:=p\sum_{j=i_3}^{N-1}\tilde{d}_{3,j}c_j^m \leq p\sum_{j=i_3}^{N-1} 1 \cdot c_j^m.
$$
Then 
$
({\bm R}{\bm c}^m)^\top ((1-p){\bm d}_1+p \tilde{\bm d}_3)= T_1+T_2+T_3+T_4.
$
For $l\in \{0\}\cup [N-1-i_1]$,
 define
\bgeq
\tilde{d}_{1,i_1+l}
:=\left\{\begin{array}{lll}
\min\left\{\frac{D-(T_1+T_3+T_4)-(1-p)\sum_{j=i_1}^{i_1-1+l}c_j^m}{(1-p)c_{i_1+l}^m},1\right\}, && 
 (T_1+T_3+T_4)+(1-p)\sum_{j=i_1}^{i_1-1+l}c_j^m< D, \\
0, &&  \inmat{otherwise}.
\end{array}
\right.
\edeq
Note that
\bgeq
T_1+(1-p)\sum_{j=i_1}^{N-1} 1\cdot c_{j}^m+T_3+{p}
\sum_{j=i_3}^{N-1}  1\cdot c_{j}^m 
=\sum_{j=1}^{N-1} c_j^m=1 \quad \inmat{and} \quad 0<D<1.
\edeq 
If $T_1+T_2+T_3+T_4=D$,
then 
$$
T_1+T_2+T_3+(1-p)\sum_{j=i_1}^{i_1-1}c_j^m =T_1+T_2+T_3<D
$$
and hence
$$
\tilde{d}_{1,i_1}=\min\left\{\frac{D-(T_1+T_3+T_4)-(1-p)\sum_{j=i_1}^{i_1-1}c_j^m}{(1-p)c_{i_1}^m},1\right\}=
\min\left\{\frac{T_2}{(1-p)c_{i_1}^m},1\right\}\overset{\eqref{eq:T123}}{=}d_{1,i_1}.
$$
If $T_1+T_2+T_3+T_4<D$,
then $D-(T_1+T_2+T_3)-p\sum_{j=i_3}^{N-2}c_j^m> pc_{N-1}^m$.
It then follows from \eqref{eq:tilded_3d_3} that
$$
\tilde{d}_{3,N-1}=\min\left\{\frac{D-(T_1+T_2+T_3)-p\sum_{j=i_3}^{N-2}c_j^m}{pc_{N-1}^m},1\right\}=1.
$$
Consequently, $\tilde{ d}_{3,j}=1$ for $j=i_3,\cdots,N-1$, 
and therefore
$$
T_4=p\sum_{j=i_3}^{N-1} 1 \cdot c_j^m,
$$
and 
$$
\tilde{d}_{1,i_1}= \min\left\{\frac{D-(T_1+T_3+T_4)-(1-p)\sum_{j=i_1}^{i_1-1}c_j^m}{(1-p)c_{i_1}^m},1\right\}\overset{D>\sum_{j=1}^4T_j}{>}
\min\left\{\frac{T_2-(1-p)\sum_{j=i_1}^{i_1-1}c_j^m}{(1-p)c_{i_1}^m},1\right\}\overset{\eqref{eq:T123}}{=}d_{1,i_1}.
$$
Define 
\bgeqn
\label{eq:label_lA2}
\tilde{l}^A:=\max\left\{l:  \tilde{d}_{1,{i}_1+{l}}>0, l\in \{0\}\cup [N-2-i_1]\right\}.
\edeqn
By an argument analogous to that used for $\tilde{d}_{3}$,
we have
\bgeqn
\label{eq:tilded_1d_1}
\tilde{d}_{1,j}> {d}_{1,j}, \;\; j\in [i_1+\tilde{l}^A]
\edeqn 
and 
$\tilde{\bm d}_1:=(1,\cdots,1,\tilde{d}_{1,i_1},\tilde{d}_{1,i_1+1},\cdots,\tilde{d}_{1,N-1})$ 
lies 
in 
the range of ${\bm g}(\cdot)$.
We now construct
\bgeqn
\label{eq:tilded^A}
\underline{\tilde{\bm d}^A
:=(1-{p})\tilde{\bm d}_1+{p}\tilde{\bm d}_3}
&=&\Big(1,\cdots,{1},(1-p)\tilde{d}_{1,i_1}+p,(1-p)\tilde{d}_{1,i_1+1}+{p},\cdots,(1-p)\tilde{d}_{1,i_3-1}+{p}, \nonumber\\
&& ~~~(1-p)\tilde{d}_{1,i_3}+p\tilde{d}_{3,i_3},\cdots,(1-p)\tilde{d}_{1,N-1}+p\tilde{d}_{3,N-1}\Big)\in \Upsilon_{\bm A}.
\edeqn
By the definitions of $\tilde{\bm d}_1$ and $\tilde{\bm d}_3$,
\bgeq
({\bm R}{\bm c}^m)^{\top}\tilde{\bm d}^A&=&\sum_{j=1}^{i_1-1}c_j^m+(1-p)\sum_{j=i_1}^{N-1}\tilde{d}_{1,j}c_{j}^m+{p}\sum_{j=i_1}^{i_3-1} c_j^m +{p}
\sum_{j=i_3}^{N-1}\tilde{d}_{3,j} c_{j}^m\\
&=& T_1+(1-p)\sum_{j=i_1}^{N-1}\tilde{d}_{1,j}c_{j}^m+T_3+T_4\\
&=&D,
\edeq
and hence $\tilde{\bm d}^A\in \Upsilon_{\bm A}^=$. 
To compare the objective values attained by $\tilde{\bm d}^A$ and ${\bm d}^A$, we consider 
two cases: $T_1+T_2+T_3+T_4=D$ and 
$T_1+T_2+T_3+T_4<D$.
\begin{itemize}
\item \underline{Case 1. $T_1+T_2+T_3+T_4=D$}.
We consider the following three subcases.

\underline{If $v_{2,i_3}^m>0$}, then
the objective function of the rhs of \eqref{eq:F_A=F_A=-lhs} attained by
$\tilde{\bm d}^A$ satisfies
\bgeq
({\bm R}{\bm v}_2^m)^{\top}\tilde{\bm d}^A&=&\sum_{j=1}^{i_1-1}v_{2,j}^m+(1-p)\tilde{d}_{1,i_1}v_{2,i_1}^m+{p}\sum_{j=i_1}^{i_3-1} v_{2,j}^m +
{p}\tilde{d}_{3,i_3} v_{2,i_3}^m\\
&&
+(1-p)\sum_{j=i_1+1}^{N-1}\tilde{d}_{1,j}v_{2,j}^m
+p\sum_{j=i_3+1}^{N-1}\tilde{d}_{3,j}v_{2,j}^m\\
&\geq & \sum_{j=1}^{i_1-1}v_{2,j}^m
+(1-p)\tilde{d}_{1,i_1}v_{2,i_1}^m
+{p}\sum_{j=i_1}^{i_3-1} v_{2,j}^m
+{p}\tilde{d}_{3,i_3} v_{2,i_3}^m\\ 
&\overset{\eqref{eq:tilded_3d_3},\tilde{\bm d}_1={\bm d}_1}{>}&\sum_{j=1}^{i_1-1}v_{2,j}^m
+(1-p){d}_{1,i_1}v_{2,i_1}^m
+{p}\sum_{j=i_1}^{i_3-1} v_{2,j}^m 
+{p}{d}_{3,i_3} v_{2,i_3}^m\\
 &\overset{\eqref{eq:dA}}{=} & \vt_{\bm A},
\edeq
which contradicts the fact that $\vt_{\bm A}$ is the optimal value of the rhs of \eqref{eq:F_A=F_A=-lhs}. 

\noindent
\underline{If $v_{2,i_3}^m=0$ and $v_{2,i_3+\bar{l}^0}^m>0$ for some $\bar{l}^0\in [\bar{l}^A]$},
then
\bgeq
({\bm R}{\bm v}_2^m)^{\top}\tilde{\bm d}^A&=&\sum_{j=1}^{i_1-1}v_{2,j}^m
+(1-p)\tilde{d}_{1,i_1}v_{2,i_1}^m
+{p}\sum_{j=i_1}^{i_3-1} v_{2,j}^m 
+{p}\tilde{d}_{3,i_3} v_{2,i_3}^m\\
&& +(1-p)\sum_{j=i_1+1}^{N-1}\tilde{d}_{1,j}v_{2,j}^m
+p\sum_{j=i_3+1}^{N-1}\tilde{d}_{3,j}v_{2,j}^m\\
&\geq & \sum_{j=1}^{i_1-1}v_{2,j}^m
+(1-p)\tilde{d}_{1,i_1}v_{2,i_1}^m
+{p}\sum_{j=i_1}^{i_3-1} v_{2,j}^m 
+{p}\tilde{d}_{3,i_3} v_{2,i_3}^m
+{p}\tilde{d}_{3,i_3+\bar{l}^0} v_{2,i_3+\bar{l}^0}^m\\ 
&\overset{\eqref{eq:tilded_3d_3},\tilde{\bm d}_1={\bm d}_1}{>}&\sum_{j=1}^{i_1-1}v_{2,j}^m+(1-p){d}_{1,i_1}v_{2,i_1}^m+{p}\sum_{j=i_1}^{i_3-1} v_{2,j}^m +{p}{d}_{3,i_3} \cdot 0+ 0
\\
&\overset{\eqref{eq:dA}}{=}&\vt_{\bm A},
\edeq
again yielding a contradiction.

\noindent\underline{If $v_{2,i_3+j}^m=0$ for all $j\in \{0\}\cup [\bar{l}^{A}]$},
then 
\bgeq
({\bm R}{\bm v}_2^m)^{\top}\tilde{\bm d}^A&=&\sum_{j=1}^{i_1-1}v_{2,j}^m+(1-p)\tilde{d}_{1,i_1}v_{2,i_1}^m+{p}\sum_{j=i_1}^{i_3-1} v_{2,j}^m +
{p}\tilde{d}_{3,i_3} v_{2,i_3}^m\\
&& +(1-p)\sum_{j=i_1+1}^{N-1}\tilde{d}_{1,j}v_{2,j}^m
+p\sum_{j=i_3+1}^{N-1}\tilde{d}_{3,j}v_{2,j}^m\\
&\overset{\tilde{\bm d}_1={\bm d}_1}{=}&\sum_{j=1}^{i_1-1}v_{2,j}^m+(1-p){d}_{1,i_1}v_{2,i_1}^m+{p}\sum_{j=i_1}^{i_3-1} v_{2,j}^m +0\\
&\overset{\eqref{eq:dA}}{=}&\vt_{\bm A},
\edeq
which implies that 
$\tilde{\bm d}^A\in \Upsilon_{\bm A}^=$ 
attains
the optimal value $\vt_{\bm A}$.
This contradicts our assumption that all optimal solutions lie in $\Upsilon_{\bm A}^<$.

\item \underline{Case 2. $T_1+T_2+T_3+T_4<D$}.
We consider the following five subcases.

\underline{If $v_{2,i_1}^m>0$}, then
the objective function of the rhs of \eqref{eq:F_A=F_A=-lhs} attained by $\tilde{\bm d}^A$ satisfies
\bgeq
({\bm R}{\bm v}_2^m)^{\top}\tilde{\bm d}^A&=&\sum_{j=1}^{i_1-1}v_{2,j}^m+(1-p)\tilde{d}_{1,i_1}v_{2,i_1}^m+{p}\sum_{j=i_1}^{i_3-1} v_{2,j}^m +
{p}\tilde{d}_{3,i_3} v_{2,i_3}^m\\
&&
+(1-p)\sum_{j=i_1+1}^{N-1}\tilde{d}_{1,j}v_{2,j}^m
+p\sum_{j=i_3+1}^{N-1}\tilde{d}_{3,j}v_{2,j}^m\\
&\geq & \sum_{j=1}^{i_1-1}v_{2,j}^m
+(1-p)\tilde{d}_{1,i_1}v_{2,i_1}^m
+{p}\sum_{j=i_1}^{i_3-1} v_{2,j}^m
+{p}\tilde{d}_{3,i_3} v_{2,i_3}^m\\ 
&\overset{\eqref{eq:tilded_3d_3},\eqref{eq:tilded_1d_1}}{>}&\sum_{j=1}^{i_1-1}v_{2,j}^m
+(1-p){d}_{1,i_1}v_{2,i_1}^m
+{p}\sum_{j=i_1}^{i_3-1} v_{2,j}^m 
+{p}{d}_{3,i_3} v_{2,i_3}^m\\
 &\overset{\eqref{eq:dA}}{=} & \vt_{\bm A},
\edeq
which contradicts the fact that $\vt_{\bm A}$ is the optimal value of the rhs of \eqref{eq:F_A=F_A=-lhs}. 

\underline{If $v_{2,i_1}^m=0$ and $v_{2,i_1+\tilde{l}^0}>0$ for some $\tilde{l}^0\in [\tilde{l}^A]$}, then
\bgeq
({\bm R}{\bm v}_2^m)^{\top}\tilde{\bm d}^A&=&\sum_{j=1}^{i_1-1}v_{2,j}^m+(1-p)\tilde{d}_{1,i_1}v_{2,i_1}^m+{p}\sum_{j=i_1}^{i_3-1} v_{2,j}^m +
{p}\tilde{d}_{3,i_3} v_{2,i_3}^m\\
&&
+(1-p)\sum_{j=i_1+1}^{N-1}\tilde{d}_{1,j}v_{2,j}^m
+p\sum_{j=i_3+1}^{N-1}\tilde{d}_{3,j}v_{2,j}^m\\
&\geq & \sum_{j=1}^{i_1-1}v_{2,j}^m
+(1-p)\tilde{d}_{1,i_1}v_{2,i_1}^m
+{p}\sum_{j=i_1}^{i_3-1} v_{2,j}^m
+{p}\tilde{d}_{3,i_3} v_{2,i_3}^m+(1-p)\tilde{d}_{1,i_1+\tilde{l}^0}v_{2,i_1+\tilde{l}^0}^m\\ 
&\overset{\eqref{eq:tilded_3d_3},\eqref{eq:tilded_1d_1}}{>}&\sum_{j=1}^{i_1-1}v_{2,j}^m
+(1-p){d}_{1,i_1}v_{2,i_1}^m
+{p}\sum_{j=i_1}^{i_3-1} v_{2,j}^m 
+{p}{d}_{3,i_3} v_{2,i_3}^m+0\\
 &\overset{\eqref{eq:dA}}{=} & \vt_{\bm A},
\edeq
again yielding a contradiction.

\underline{If $v_{2,i_1+j}^m=0$ for all $j\in\{0\}\cup [\tilde{l}^A]$ and $v_{2,i_3}^m>0$},
then 
\bgeq
({\bm R}{\bm v}_2^m)^{\top}\tilde{\bm d}^A&=&\sum_{j=1}^{i_1-1}v_{2,j}^m+(1-p)\tilde{d}_{1,i_1}v_{2,i_1}^m+{p}\sum_{j=i_1}^{i_3-1} v_{2,j}^m +
{p}\tilde{d}_{3,i_3} v_{2,i_3}^m\\
&& +(1-p)\sum_{j=i_1+1}^{N-1}\tilde{d}_{1,j}v_{2,j}^m
+p\sum_{j=i_3+1}^{N-1}\tilde{d}_{3,j}v_{2,j}^m\\
&\geq & \sum_{j=1}^{i_1-1}v_{2,j}^m+(1-p)\tilde{d}_{1,i_1}v_{2,i_1}^m+{p}\sum_{j=i_1}^{i_3-1} v_{2,j}^m +{p}\tilde{d}_{3,i_3} v_{2,i_3}^m\\ 
&\overset{\eqref{eq:tilded_3d_3},\eqref{eq:tilded_1d_1}}{>}&\sum_{j=1}^{i_1-1}v_{2,j}^m+(1-p){d}_{1,i_1}\cdot 0+{p}\sum_{j=i_1}^{i_3-1} v_{2,j}^m +{p}{d}_{3,i_3} v_{2,i_3}^m\\
\\
&\overset{\eqref{eq:dA}}{=}&\vt_{\bm A},
\edeq
again yielding a contradiction.

\underline{If $v_{2,i_1+j}^m=0$ for all $j\in  [\tilde{l}^A]$, $v_{2,i_3}=0$, and  $v_{2,i_3+l^A}^m>0$ for some $l^A\in [\bar{l}^{A}]$},
then
\bgeq
({\bm R}{\bm v}_2^m)^{\top}\tilde{\bm d}^A&=&\sum_{j=1}^{i_1-1}v_{2,j}^m+(1-p)\tilde{d}_{1,i_1}v_{2,i_1}^m+{p}\sum_{j=i_1}^{i_3-1} v_{2,j}^m
+
{p}\tilde{d}_{3,i_3} v_{2,i_3}^m\\
&& +(1-p)\sum_{j=i_1+1}^{N-1}\tilde{d}_{1,j}v_{2,j}^m
+p\sum_{j=i_3+1}^{N-1}\tilde{d}_{3,j}v_{2,j}^m\\
&\geq & \sum_{j=1}^{i_1-1}v_{2,j}^m+(1-p)\tilde{d}_{1,i_1}v_{2,i_1}^m+{p}\sum_{j=i_1}^{i_3-1} v_{2,j}^m +{p}\tilde{d}_{3,i_3} v_{2,i_3}^m
+p \tilde{d}_{3,i_3+l^A}v_{2,i_3+l^A}^m\\ 
&\overset{\eqref{eq:tilded_3d_3},\eqref{eq:tilded_1d_1}}{>}&\sum_{j=1}^{i_1-1}v_{2,j}^m+(1-p){d}_{1,i_1}\cdot 0+{p}\sum_{j=i_1}^{i_3-1} v_{2,j}^m +{p}{d}_{3,i_3}\cdot 0\\
&\overset{\eqref{eq:dA}}{=}&\vt_{\bm A},
\edeq
again yielding a contradiction.

\noindent\underline{If $v_{2,i_1+j}^m=0$ for all $j\in  [\tilde{l}^A]$, and  $v_{2,i_3+j}^m=0$ for all $j\in [\bar{l}^{A}]$},
then 
\bgeq
({\bm R}{\bm v}_2^m)^{\top}\tilde{\bm d}^A&=&\sum_{j=1}^{i_1-1}v_{2,j}^m+(1-p)\tilde{d}_{1,i_1}v_{2,i_1}^m+{p}\sum_{j=i_1}^{i_3-1} v_{2,j}^m +
{p}\tilde{d}_{3,i_3} v_{2,i_3}^m\\
&& +(1-p)\sum_{j=i_1+1}^{N-1}\tilde{d}_{1,j}v_{2,j}^m
+p\sum_{j=i_3+1}^{N-1}\tilde{d}_{3,j}v_{2,j}^m\\
&=&\sum_{j=1}^{i_1-1}v_{2,j}^m+p\sum_{j=i_1}^{i_3-1} v_{2,j}^m+0\\
&=&\sum_{j=1}^{i_1-1}v_{2,j}^m+(1-p){d}_{1,i_1}\cdot 0+{p}\sum_{j=i_1}^{i_3-1} v_{2,j}^m +{p}{d}_{3,i_3} \cdot 0\\
&\overset{\eqref{eq:dA}}{=}&\vt_{\bm A},
\edeq
which implies that
$\tilde{\bm d}^A\in \Upsilon_{\bm A}^=$ 
attains
the optimal value $\vt_{\bm A}$.
This contradicts our assumption that all optimal solutions lie in
$\Upsilon_{\bm A}^{<}$.
\end{itemize}
The proof is complete. 
\hfill $\Box$

Proposition~\ref{prop:EDn=0_2-1st}
shows that the optimum of problems \eqref{eq:gene_m_Q-1-v1-p}-\eqref{eq:gene_m_Q-2-v1-p}
are attained in 
${\cal F}_{\bm A}^{=}$ and ${\cal F}_{\bm B}^{=}$,
respectively.
However, this does not rule out
the case that the optimum 
of problem \eqref{eq:gene_m_Q-2-v1-p}
is also attained in  ${\cal F}_{\bm B} \backslash {\cal F}_{\bm B}^{=}$.
The next proposition
provides conditions on
${\bm v}_1^m$ and ${\bm v}_2^m$,
under which all
optimal solutions lie in 
${\cal F}_{\bm A}^{=}$ and $ {\cal F}_{\bm B}^{=} $,
respectively.

\begin{proposition}
[Properties of the Optimal Solutions of Problems \eqref{eq:gene_m_Q-1-v1-p}-\eqref{eq:gene_m_Q-2-v1-p}]
\label{prop:EDn=0_2}
Assume the setting and conditions of Proposition~\ref{prop:EDn=0_2-1st}.
If ${\bm v}_1^m$ and ${\bm v}_2^m$ satisfy
\begin{subequations}
\label{eq:condition_thm}
\begin{align}
&\inmat{C1:}\;\;{v}_{1,j_{r_2^m}}^m>0\;\; \inmat{for}\;\;r_2^m \notin {\cal X},\;\;\inmat{or}
\;\;{v}_{1,j_{r_2^m}+1}^m>0\;\;\inmat{for}\;\; r_2^m \in {\cal X},\\
& \inmat{C2:}\;\;
{v}_{2,j_{{r}_\tau^m}}^m> 0\;\;\inmat{and}\;r_{\tau}^m \notin {\cal X} \;\;\inmat{for}\;\;
\tau=1,3,\label{eq:condition_thm-b}
\end{align}
\end{subequations}
where $r_2^m\in {\cal S}_{\bm B}$,
$(r_1^m,r_3^m,p)\in {\cal S}_{\bm A}$,
and $j_{{r}_i^m}$ 
denote
the index such that
${r}_i^m\in (x_{j_{{r}_i^m}},x_{j_{{r}_i^m}+1}]$
for $i=1,2,3$, 
then 
\bgeqn
\label{eq:rel_solution}
{\cal S}_{\bm B}\subset {\cal F}_{\bm B}^{=}
\quad\inmat{and}\quad
{\cal S}_{\bm A}\subset {\cal F}_{\bm A}^{=}.
\edeqn
\end{proposition}

\proof{Proof of Proposition~\ref{prop:EDn=0_2}.}
\underline{First,
we show  ${\cal S}_{\bm B}\subset {\cal F}_{\bm B}^{=}$ in \eqref{eq:rel_solution}}.
Suppose, for contradiction， 
${\cal S}_{\bm B}\cap ({\cal F}_{\bm B}\backslash{{\cal F}_{\bm B}^{=}})\neq \emptyset$.
Then
there exists an optimal solution to 
the rhs problem of \eqref{eq:F_B=F_B=-lhs},
denoted by
${\bm d}^B$,
lying in 
$\Upsilon_{\bm B}^{<}$.
Let $r_2^m$ be 
the corresponding optimal solution to the lhs of \eqref{eq:F_B=F_B=-lhs}
so hat
${\bm d}^B={\bm g}(r_2^m)$.
Since
${\bm R}{\bm c}^m>0$,
the structure  
of 
${\bm g}(\cdot)$ implies that
there exists  $i^B\in [N-1]$ such that 
\bgeqn
\label{eq:d^B_unique}
{\bm d}^B=(1,\cdots,1,d_{i^B},0,\cdots,0)=\left(1,\cdots,1,\frac{r_2^m-x_{i^B}}{x_{i^B+1}-x_{i^B}},0,\cdots,0\right)\in \Upsilon_{\bm B}^{<},
\edeqn
where
$i^B=j_{r_2^m}$.
Then
$({\bm R}{\bm c}^m)^\top {\bm d}^B<D$,
which implies that
$\frac{r_2^m-x_{i^B}}{x_{i^B+1}-x_{i^B}}<\frac{D-\sum_{j=1}^{i^B-1}c_j^m}{c_{i^B}^m}$.
Let $\tilde{\bm d}^B$ be defined as in \eqref{eq:tilded^B}.
Under the condition C1,
if
$v_{1,j_{r_2^m}}^m>0$ and $r_2^m \notin {\cal X}$,
then $\frac{r_2^m-x_{j_{r_2^m}}}{x_{j_{r_2^m}+1}-x_{j_{r_2^m}}}<1$.
Hence,
the 
objective value of the rhs of \eqref{eq:F_B=F_B=-lhs} 
evaluated at $\tilde{\bm d}^B$
satisfies  
\bgeq
({\bm R}{\bm v}_1^m)^\top \tilde{\bm d}^B&=&\sum_{j=1}^{j_{r_2^m}-1}v_{1,j}^m+\tilde{d}^B_{j_{r_2^m}}v_{1,j_{r_2^m}}^m+\tilde{d}^B_{j_{r_2^m}+1}v_{1,j_{r_2^m}+1}^m+\cdots+ \tilde{d}^B_{N-1}v_{1,N-1}^m \nonumber \\
&
\geq&\sum_{j=1}^{j_{r_2^m}-1}v_{1,j}^m+\min\left\{\frac{D-\sum_{j=1}^{j_{r_2^m}-1}c_j^m}{c_{j_{r_2^m}}^m},1\right\}v_{1,j_{r_2^m}}^m\\
&>&\sum_{j=1}^{j_{r_2^m}-1}v_{1,j}^m+\frac{r_2^m-x_{j_{r_2^m}}}{x_{j_{r_2^m}+1}-x_{j_{r_2^m}}}v_{1,j_{r_2^m}}^m\overset{\eqref{eq:d^B_unique}}{=}\vt_{\bm B},\nonumber 
\edeq
which contradicts the assumption that $\vt_{\bm B}$ is the optimal value of the rhs of \eqref{eq:F_B=F_B=-lhs}.
Under the condition C1,
if 
${v}_{1,j_{r_2^m}+1}^m>0\,\inmat{for}\;r_2^m \in {\cal X}$,
then $\frac{r_2^m-x_{j_{r_2^m}}}{x_{j_{r_2^m}+1}-x_{j_{r_2^m}}}=1$, 
$\frac{D-\sum_{j=1}^{j_{r_2^m}}c_j^m}{c_{j_{r_2^m}+1}^m}>0$.
Hence,
the objective value of the rhs of \eqref{eq:F_B=F_B=-lhs}
attained at $\tilde{\bm d}^B$ satisfies  
\bgeq
({\bm R}{\bm v}_1^m)^\top \tilde{\bm d}^B&=&\sum_{j=1}^{j_{r_2^m}}v_{1,j}^m+\tilde{d}^B_{j_{r_2^m}+1}v_{1,j_{r_2^m}+1}^m+\cdots+ \tilde{d}^B_{N-1}v_{1,N-1}^m \nonumber \\
&
\geq&\sum_{j=1}^{j_{r_2^m}}v_{1,j}^m+
\min\left\{\frac{D-\sum_{j=1}^{j_{r_2^m}}c_j^m}{c_{j_{r_2^m}+1}^m},1\right\}v_{1,j_{r_2^m}+1}^m\\
&>&\sum_{j=1}^{j_{r_2^m}}v_{1,j}^m\overset{\eqref{eq:d^B_unique}}{=}\vt_{\bm B},\nonumber 
\edeq
which again contradicts the assumption that $\vt_{\bm B}$ is the optimal value of the rhs of \eqref{eq:F_B=F_B=-lhs}.
In either case, we obtain a contradiction. Therefore,
${\cal S}_{\bm B}\cap ({\cal F}_{\bm B}\backslash{{\cal F}_{\bm B}^{=}})=\emptyset$,
which establishes
${\cal S}_{\bm B}\subset {\cal F}_{\bm B}^{=}$.

\underline{To complete the proof,
we show  ${\cal S}_{\bm A}\subset {\cal F}_{\bm A}^{=}$ in \eqref{eq:rel_solution}}.
Suppose, for contradiction, that
${\cal S}_{\bm A} \cap ({\cal F}_{\bm A}\backslash{{\cal F}_{\bm A}^{=}})\neq \emptyset$.
Then there exists an optimal solution to rhs of \eqref{eq:F_A=F_A=-lhs}, 
denoted by
${\bm d}^A$,
that lies in
$\Upsilon_{\bm A}^<$.
Let 
$(r_1^m,r_3^m,p)$ be
a corresponding
optimal solution of lhs of of \eqref{eq:F_A=F_A=-lhs},
with
$r_1^m,r_3^m\in [\underline{x},\bar{x}]$, 
such that ${\bm d}_1={\bm g}(r_1^m)=\left(1,\cdots,1,\frac{r_1^m-x_{i_1}}{x_{i_1+1}-x_{i_1}},0,\cdots,0\right)$ and ${\bm d}_3={\bm g}(r_3^m)=\left(1,\cdots,1,\frac{r_3^m-x_{i_3}}{x_{i_3+1}-x_{i_3}},0,\cdots,0\right)$,
with  $0\leq \frac{r_1^m-x_{i_1}}{x_{i_1+1}-x_{i_1}}<1$,
$0\leq \frac{r_3^m-x_{i_3}}{x_{i_3+1}-x_{i_3}}<1$,
and
\bgeq
{\bm d}^A&=&(1-p){\bm g}(r_1^m)+p{\bm g}(r_3^m)\\
&=& \left(1,\cdots,1,(1-p)\frac{r_1^m-x_{i_1}}{x_{i_1+1}-x_{i_1}}+p,p,\cdots,p,p\frac{r_3^m-x_{i_3}}{x_{i_3+1}-x_{i_3}},0,\cdots,0\right),
\edeq
where $j_{r_1^m}=i_1$
and
$j_{r_3^m}=i_3$.
Since ${\bm d}^A\in \Upsilon_{\bm A}^{<}$,
we have
$({\bm R}{\bm c}^m)^\top {\bm d}^{A}<D$,
which implies that 
$$
(1-p)\frac{r_1^m-x_{i_1}}{x_{i_1+1}-x_{i_1}}c_{i_1}^m+p\frac{r_3^m-x_{i_3}}{x_{i_3+1}-x_{i_3}}c_{i_3}^m<D-\left(\sum_{j=1}^{i_1-1}c_j^m
+p\sum_{j=i_1}^{i_3-1} c_j^m \right).
$$
Choose sufficiently small
$\epsilon_1>0$ and $\epsilon_2>0$
such that 
$0<\left(\frac{r_1^m-x_{i_1}}{x_{i_1+1}-x_{i_1}}+\epsilon_1\right)<1$,
$0<\left(\frac{r_3^m-x_{i_3}}{x_{i_3+1}-x_{i_3}}+\epsilon_2\right)<1$,
and
$$
(1-p)\left(\frac{r_1^m-x_{i_1}}{x_{i_1+1}-x_{i_1}}c_{i_1}^m+\epsilon_1\right)+p\left(\frac{r_3^m-x_{i_3}}{x_{i_3+1}-x_{i_3}}c_{i_3}^m+\epsilon_2\right)<D-\left(\sum_{j=1}^{i_1-1}c_j^m
+p\sum_{j=i_1}^{i_3-1} c_j^m \right).
$$
Let
$$
\hat{\bm d}_1(\epsilon_1):= (1,\cdots,1,\frac{r_1^m-x_{i_1}}{x_{i_1+1}-x_{i_1}}+\epsilon_1,0,\cdots,0)
\quad \inmat{and} \quad  
\hat{\bm d}_3(\epsilon_2):= (1,\cdots,1,\frac{r_3^m-x_{i_3}}{x_{i_3+1}-x_{i_3}}+\epsilon_2,0,\cdots,0).
$$
Then, define
\bgeq
\hat{\bm d}^A&:=&(1-{p})\hat{\bm d}_1(\epsilon_1)+{p}\hat{\bm d}_3(\epsilon_2)\\
&=&\Big(1,\cdots,{1},\overbrace{(1-{p})\left(\frac{r_1^m-x_{i_1}}{x_{i_1+1}-x_{i_1}}+\epsilon_1\right)+{p}}^{i_1},{p},\cdots,{p},\overbrace{p\left(\frac{r_3^m-x_{i_3}}{x_{i_3+1}-x_{i_3}}+\epsilon_2\right)}^{i_3},0,\cdots,0 
\Big).
\edeq
By construction,
$\hat{\bm d}^A \in \Upsilon_{\bm A}$ and 
$({\bm R}{\bm c}^m)^\top \hat{\bm d}^{A}< D$.
Under the condition C2 in \eqref{eq:condition_thm-b},
the
objective value of rhs of \eqref{eq:F_A=F_A=-lhs},
evaluated at
$\hat{\bm d}^A$ satisfies
\bgeq
({\bm R}{\bm v}_2^m)^{\top}\hat{\bm d}^A&=&\sum_{j=1}^{i_1-1}v_{2,j}^m+(1-{p})\left(\frac{r_1^m-x_{i_1}}{x_{i_1+1}-x_{i_1}}+\epsilon_1 \right)v_{2,i_1}^m
+{p}\sum_{j=i_1}^{i_3-1} v_{2,j}^m +{p} \left( \frac{r_3^m-x_{i_3}}{x_{i_3+1}-x_{i_3}}c_{i_3}^m+\epsilon_2\right)v_{2,i_3}^m\\
 &>&\sum_{j=1}^{i_1-1}v_{2,j}^m+(1-{p})\frac{r_1^m-x_{i_1}}{x_{i_1+1}-x_{i_1}}v_{2,i_1}^m
+{p}\sum_{j=i_1}^{i_3-1} v_{2,j}^m +{p} \frac{r_3^m-x_{i_3}}{x_{i_3+1}-x_{i_3}}v_{2,i_3}^m\\
&=&({\bm R}{\bm v}_2^m)^{\top}{\bm d}^A=\vt_{\bm A},
\edeq
which contradicts the assumption that $\vt_{\bm A}$ is the optimal value of rhs of \eqref{eq:F_A=F_A=-lhs}.
\hfill
\Box

Condition C1 ensures that the objective 
function of program \eqref{eq:gene_m_Q-1-v1-p} is
strictly increasing in the
neighborhood of  $r_2^m$ and thereby  
$r_2^m$ must satisfy the equality constraint.
Likewise Condition C2 ensures the objective function of program \eqref{eq:gene_m_Q-2-v1-p} to be strictly increasing in the neighborhoods of  $r_1^m$ and $r_3^m$ and  
hence
$(r_1^m,r_3^m,p)$ must satisfy 
the equality constraint.

Based on the discussion above, 
we are ready to establish the relationship between the optimal solutions of problems
\eqref{eq:gene_m_Q-1-v1-p}-\eqref{eq:gene_m_Q-2-v1-p}.

\begin{proposition}
\label{pro:new-prop}
Let ${\cal F}_{{\bm A}}^=$ and ${\cal F}_{{\bm B}}^=$
be defined as in \eqref{eq:bar_r123-2}, and 
${\cal S}_{\bm A}$ and 
${\cal S}_{\bm B}$ 
as in \eqref{eq:solution_B}.
For any $({r}_1^m,{r}_3^m,{p}^m)\in {\cal S}_{\bm A}\cap {\cal F}_{\bm A}^=$
and ${r}_2^m\in {\cal S}_{\bm B}\cap{\cal F}_{{\bm B}}^=$,
we have 
$r_1^m\leq r_2^m\leq r_3^m$.
\end{proposition}

\proof{Proof of Proposition~\ref{pro:new-prop}.}
Suppose, for contradiction, that
$r_2^m>r_3^m$.

Since $({\bm R}{\bm c}^m)^{\top}G_{\bm A}({r}_1,{r}_3,{p})= 
({\bm c}^m,1-{\bm e}^\top {\bm c}^m)^\top 
((1-{p}){\bm g}({r}_1)+{p}{\bm g}({r}_3))
$ is strictly increasing w.r.t.~each variate,
we have 
$$
({\bm R}{\bm c}^m)^{\top}G_{\bm A}({r}_1^m,{r}_3^m,{p})< ({\bm R}{\bm c}^m)^{\top}G_{\bm A}({r}_2^m,{r}_2^m,{p})=({\bm R}{\bm c}^m)^{\top}G_{\bm B}({r}_2^m)=D,
$$
which implies that $(r_1^m,r_3^m,p^m)\not\in {\cal F}_{\bm A}^{=}$, a contradiction.
Likewise,
if 
$r_1^m> r_2^m$,
then
$$
({\bm R}{\bm c}^m)^{\top}G_{\bm A}({r}_1^m,{r}_3^m,{p}^m)> ({\bm R}{\bm c}^m)^{\top}G_{\bm A}({r}_2^m,{r}_2^m,{p}^m)=({\bm R}{\bm c}^m)^{\top}G_{\bm B}({r}_2^m)=D,
$$
which also contradicts the assumption that $(r_1^m,r_3^m,p^m)\in {\cal F}_{\bm A}^{=}$.

Therefore,
$r_1^m \leq r_2^m \leq r_3^m$, 
which completes the proof.
\hfill $\Box$

\section{Fixed-Breakpoint Variant}
\label{sec:variant_fixed}

We fix the set of breakpoints ${\cal X}_0$ with $N=|{\cal X}_0|$ and develop the following variant of Algorithm~\ref{alg:RUS}.
First, we control the generation of $r_2^m$ by proper choice of $D_s$, e.g.,
\begin{eqnarray}
\label{eq:D_s_set}
D_s\in {\cal D}^m:=
\left\{
c_1^m,\,
\sum_{j=1}^{2}c_j^m,\,
\ldots,\,
{\bm e}^{\top}{\bm c}^m
\right\}.
\end{eqnarray}
The corresponding solution $r_2^{m,s}$ to
$
({\bm R}{\bm c}^m)^{\top}
G_{\bm A}(r_2^{m,s})=D_s
$
lies in ${\cal X}_{m-1}$, 
since each value of $D_s$ corresponds to the cumulative utility increment up to that breakpoint.

Second,
for each admissible pair $(r_1^m,r_3^m)\in{\cal X}_{m-1}\times{\cal X}_{m-1}$, we solve
$
({\bm R}{\bm c}^m)^\top
G_{\bm A}(r_1^m,r_3^m,p^m)=D
$
for $p^m$ and retain combinations $(r_1^m,r_3^m,p^m)$
if $p^m\in[\underline p,\bar p]$.
We select the combination from all feasible combinations
whose normal vector of the induced cutting plane is most closely aligned with the longest axis of the Sonnevend inner ellipsoid of ${\cal V}_N^{m-1}$.
Algorithm~\ref{alg:RUS-fixed} summarizes this procedure.

\begin{algorithm}
\caption{Variant of the Modified Polyhedral Method (with Fixed Breakpoints)}
\label{alg:RUS-fixed}
\begin{algorithmic}[1]
\footnotesize
\STATE \textbf{Input:}
Initial breakpoint set ${\cal X}_0$, number of queries $M$, and $m=1$.
Set $N=|{\cal X}_0|$.
\REPEAT
\STATE
Compute the analytic center ${\bm c}^m$ of ${\cal V}_{N}^{m-1}$
and the intersections ${\bm v}_1^m,{\bm v}_2^m$ of the longest
axis with the boundary of ${\cal V}_{N}^{m-1}$.
\STATE
Solve problem \eqref{eq:gene_m_Q-1-v1-p} with a sequence of parameters
$D_s\in{\cal D}^m$ to obtain the candidates $r_2^{m,s}$.
For each $D_s\in{\cal D}^m$, construct
$(r_1^{m,s},r_3^{m,s},p^{m,s})$ as follows:
consider all feasible $(r_1,r_3,p)$ satisfying
\[
({\bm R}{\bm c}^m)^\top
G_{\bm A}(r_1,r_3,p)=D_s,
\qquad
p\in[\underline p,\bar p],
\]
and select $(r_1^{m,s},r_3^{m,s},p^{m,s})$ whose normal vector of the induced cutting plane is most
closely aligned with the longest-axis direction of the Sonnevend inner
ellipsoid of ${\cal V}_N^{m-1}$.
Ask the DM to select between ${\bm A}_m=(r_1^m,r_3^m,p^m)$ and ${\bm B}_m=r_2^m$.
\STATE
\textbf{Update the Ambiguity Set.}
Set
\begin{eqnarray}
\label{sec:V_fixed}
{\cal V}_{N}^{m}
:=
{\cal V}_{N}^{m-1}
\bigcap
\left\{
{\bm v}\in\mathbb R^{N-2}:
h_m({\bm R}{\bm v})^\top
\left(
G_{\bm A}^{{\cal X}_0}(r_1^m,r_3^m,p^m)
-
G_{\bm B}^{{\cal X}_0}(r_2^m)
\right)
\leq 0
\right\},
\end{eqnarray}
where $h_m=1$ if ${\bm A}_m\preceq{\bm B}_m$ and $h_m=-1$
otherwise.
\STATE
Update $m:=m+1$.
\UNTIL $m>M$.
\STATE \textbf{Output:}
Polyhedron ${\cal V}_{N}^{M}$.
\end{algorithmic}
\end{algorithm}
We next provide an illustration of 
Algorithm~\ref{alg:RUS-fixed}. 
For each query $m\in[M]$,
let ${\bm d}^m$ be defined as in \eqref{eq:postchoice-1}.
We denote the linear inequality induced by the $m$-th query in
${\cal V}_{N}^{m}$ (see \eqref{sec:V_fixed}) by
$({\bm a}^{\,m})^\top{\bm v}\leq b_m$,
where
${\bm a}^{\,m}
:=
h_m
\frac{
{\bm d}^{m}_{[N-2]}-d_{N-1}^{m}{\bm e}
}{
\left\|
{\bm d}^{m}_{[N-2]}-d_{N-1}^{m}{\bm e}
\right\|
}
$
and
$b_m
:=
h_m
\frac{
d_{N-1}^{m}
}{
\left\|
{\bm d}^{m}_{[N-2]}-d_{N-1}^{m}{\bm e}
\right\|
}$.
Accordingly, the polyhedron ${\cal V}_N^M$ can be represented using normalized inequalities:
\begin{eqnarray}
{\cal V}_N^M
=
\left\{
{\bm v}\in\mathbb{R}^{N-2}:
-{\bm v}\leq 0,\; \; \frac{{\bm e}^\top {\bm v}}{\sqrt{N-2}}\leq \frac{1}{\sqrt{N-2}},\;\;
({\bm a}^{\,m})^\top{\bm v}\leq b_m,\; m\in[M]
\right\}.
\end{eqnarray}
Let ${\bm \delta}^m$ denote the longest-axis direction of the Sonnevend
inner ellipsoid of ${\cal V}_N^{m-1}$. 
We also define the alignment between the normal vector of the cutting plane and the longest-axis direction as
$
\eta_m:=
\frac{|({\bm a}^m)^\top{\bm \delta}^m|}
{\|{\bm \delta}^m\|_2}.
$
Hence, $\eta_m=|\cos\theta_m|\in[0,1]$, where $\theta_m$ is the angle
between the two directions, and a larger $\eta_m$ indicates a stronger
alignment with the longest-axis direction.
Table~\ref{tab:selected_alignment} reports the distribution of the selected alignments over the 50 elicitation queries. 
Notably, 42\% of the selected cutting hyperplanes attain an alignment of at least 0.8, including 22\% with an alignment of at~least~0.9.

\begin{table}[!h]
\vspace{-0.1cm}
\centering
\caption{Distribution of the selected alignment} \label{tab:selected_alignment} 
{\small \setlength{\tabcolsep}{3pt} 
\begin{tabular}{c|cccccccccc} \hline $\eta_m$ & $[0,0.1)$ & $[0.1,0.2)$ & $[0.2,0.3)$ & $[0.3,0.4)$ & $[0.4,0.5)$ & $[0.5,0.6)$ & $[0.6,0.7)$ & $[0.7,0.8)$ & $[0.8,0.9)$ & $[0.9,1]$ \\ \hline Number & 0 & 4 & 3 & 2 & 6 & 5 & 4 & 5 & 10 & 11 \\ Percentage & 0\% & 8\% & 6\% & 4\% & 12\% & 10\% & 8\% & 10\% & 20\% & 22\% \\ 
\hline
\end{tabular}}
 \vspace{-0.2cm}
\end{table}

\begin{example}[Elicitation Effectiveness: Continuation of Example~\ref{exa:concave_breakpoints}]We examine the 
effectiveness of the elicitation procedure.
Figure~\ref{fig:Fig_accuracy_50} reports the evolution of the associated utility function 
ambiguity set from three perspectives.
Figure~\ref{fig:Fig_accuracy_50}(a) depicts the range of utility functions after
different numbers of queries.
Next, we investigate the shrinkage of 
${\cal V}_N^M$.
Let $r_{\rm AC}^M:=\min\{\min_{m\in[M]}\{b^m-({\bm a}^m)^\top{\bm c}^{m+1}\},
\min_{i\in[N-2]}c_i^{M+1},
(1-{\bm e}^\top{\bm c}^{M+1})/\sqrt{N-2}\}$
denote the radius of the largest
Euclidean ball centered at ${\bm c}^{M+1}$ and contained in ${\cal V}_N^M$.
Figure~\ref{fig:Fig_accuracy_50}(b) depicts the minimum slack $r_{\rm AC}^m$
against the increase in the number of queries.
It exhibits a clear overall decreasing trend as the elicitation proceeds.
Notice that $r_{\rm AC}^m$ need not decrease monotonically even though
${\cal V}_N^{m+1}\subseteq{\cal V}_N^m$.
This is because the analytic center is recomputed after each query and
may relocate within the updated ambiguity set.

To distinguish this local analytic-center measure from the overall size
of the ambiguity set, Figure~\ref{fig:Fig_accuracy_50}(c) reports three additional
measures. 
As shown in Figure~\ref{fig:Fig_accuracy_50}(c), both $\dd_{\rm R1}$ and $\dd_{\rm R2}$
decrease monotonically as additional preference constraints are
introduced, consistent with the nesting
${\cal V}_N^{m+1}\subseteq{\cal V}_N^m$. Meanwhile,
$\dd_{\rm AC}$ also exhibits a substantial overall reduction, indicating
that the utility function represented by the analytic center 
approaches
increasingly closer
to the true utility function.
Together,
Figures~\ref{fig:Fig_accuracy_50}(a)--(c) show that the elicitation procedure not
only reduces the ambiguity set globally but also improves the accuracy
of the representative utility function, while the nonmonotonic local
behavior of the minimum slack is explained by the relocation of the
analytic center.
From Figure~\ref{fig:Fig_accuracy_50}(b)-(c),
we can easily identify the number of queries required for a specified precision.

\begin{figure}[!htbp]
 \vspace{-0.5cm}
\centering
\begin{minipage}[t]{0.32\linewidth}
    \centering
    \includegraphics[width=2.1in]{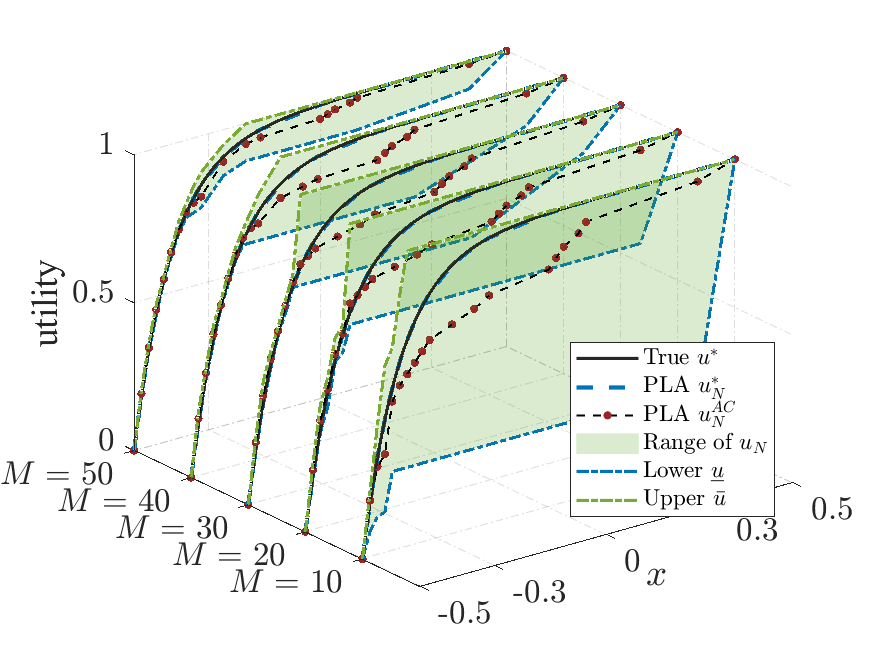}
    \text{\footnotesize{\textbf{(a)}}}
\end{minipage}
\hfill
\begin{minipage}[t]{0.32\linewidth}
    \centering
    \includegraphics[width=2.1in]{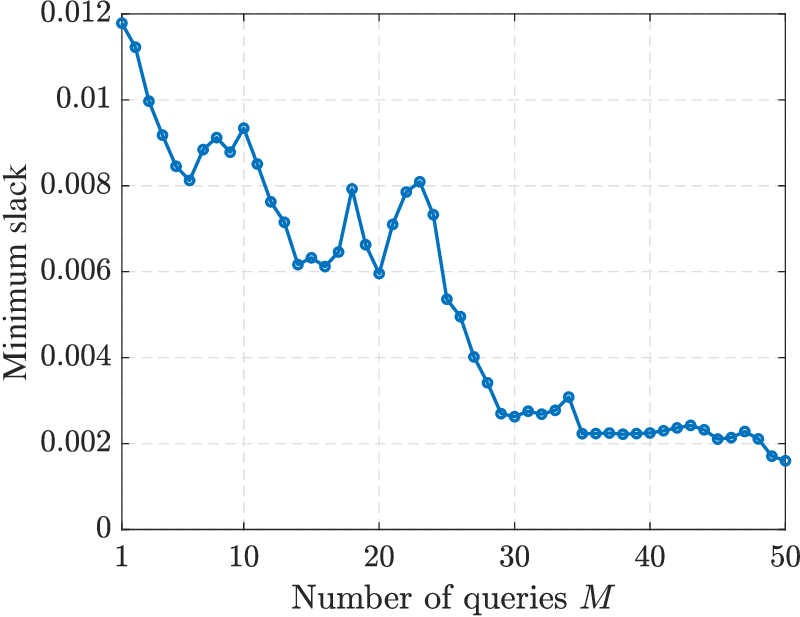}
    \text{\footnotesize{\hspace{0.6cm}\textbf{(b)}}}
\end{minipage}
\hfill
\begin{minipage}[t]{0.32\linewidth}
    \centering
    \includegraphics[width=2.2in]{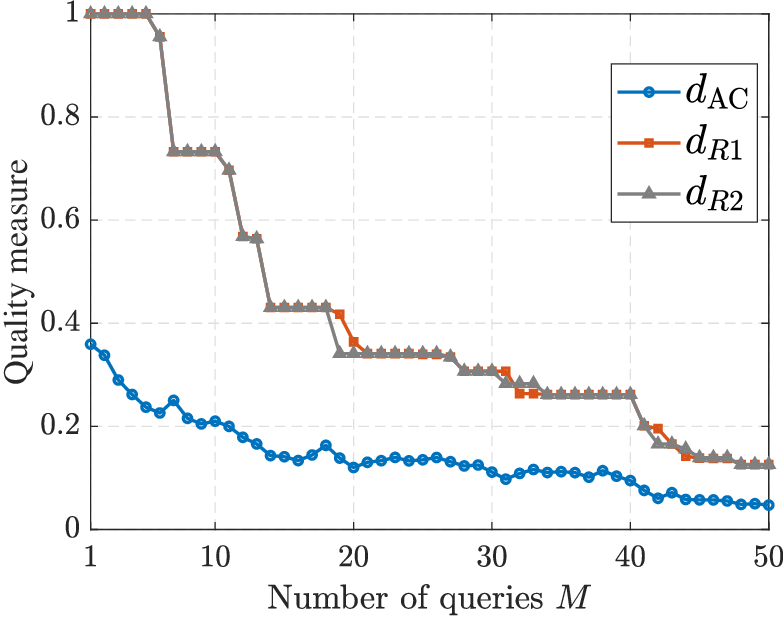}
    \text{\footnotesize{\hspace{0.6cm}\textbf{(c)}}} 
\end{minipage}
\caption{(a) Range of utility functions.
(b) Minimum slack value of ${\cal V}_{N_M}^M$.
(c) Three metrics characterizing the size of 
${\cal V}_{N_M}^M$.
}
\label{fig:Fig_accuracy_50}
\vspace{-0.5cm}
\end{figure}
\end{example}

The strengthened strategy in line~4 of Algorithm~\ref{alg:RUS-fixed}
selects $(r_1^m,r_3^m)$ from different combinations such that the normal
vector of the resulting cutting hyperplane is as parallel as possible to
the longest axis of the Sonnevend ellipsoid. Moreover, by
Theorem~\ref{prop:EDn=0-2}, each cutting hyperplane passes through the
analytic center of the current polyhedron. These two properties allow
us to adapt the analysis in
\cite[Section~3.2]{boyd2003localization} to investigate the shrinkage of
${\cal V}_N^M$ through the change in its minimum slack and to derive a
complexity result for the fixed-breakpoint variant of the MPM.

\section{Proof of Proposition~\ref{prop:reformulate_1}}
\proof{Proof of Proposition~\ref{prop:reformulate_1}.}
By the definitions of $\underline{\bm v}$ and $\bar{\bm v}$ in \eqref{eq:bounds_vo-0},
we have 
${v}_i\in [{c}_i^{M}-\underline{v}_i,{c}_i^{M}+\bar{v}_i]$ for $i\in [ N-1]$,
where
$c_{ N-1}^{M}=1-{\bm e}^\top {\bm c}^{M}$
and 
${\bm c}^{M}\in \R^{ N-2}$. 
Problem \eqref{eq:minmin_bounds-a0}
can be written as
\begin{subequations}
\label{eq:min_expand}
 \begin{align}
\displaystyle  \min_{{\bm v}_{[N-1]}} \quad & \left({\bm v}_{[N-1]}-{\bm c}_{[N-1]}^{M}\right)^\top \left(\frac{1}{K}\sum_{k=1}^K{\bm g} ({\bm z}^\top {\bm \xi}^k)\right)+({\bm c}_{[N-1]}^{M})^\top\left(\frac{1}{K}\sum_{k=1}^K{\bm g} ({\bm z}^\top {\bm \xi}^k)\right) \\
{\rm s.t.}\quad~ &  {v}_i-{c}_i^{M}\in [-\underline{v}_i,\bar{v}_i], ~i\in [N-1],\\
& {\bm e}^\top \left({\bm v}_{[N-1]}-{\bm c}^{M}_{[N-1]}\right)+{\bm e}^\top {\bm c}_{[N-1]}^{M}=1,\\
& {v}_i-{c}^{M}_i\geq -{c}^{M}_i,i\in [N-1],\\
& h_m \left({\bm v}_{[N-1]}-{\bm c}_{[N-1]}^{M}\right)^{\top}
 \left(G_{\bm A} (r_{1}^{m},r_3^{m},p^{m})-G_{\bm B} (r_2^{m})\right) \nonumber \\
 &
 +h_m ({\bm c}_{[N-1]}^{M})^{\top}
 \left(G_{\bm A} (r_{1}^{m},r_3^{m},p^{m})-G_{\bm B} (r_2^{m})\right)\leq 0, \;\; m\in [M].
\end{align}
\end{subequations}
By replacing
${\bm v}_{[N-1]}-{\bm c}_{[N-1]}^{M}$
with a new variable 
$\widetilde{\bm v}_{[N-1]}$,
we can 
reformulate 
problem \eqref{eq:min_expand}  
as 
\begin{subequations}
\begin{align}
\displaystyle \min_{\widetilde{\bm v}_{[N-1]}} \quad &
\widetilde{\bm v}_{[N-1]}^\top \left(\frac{1}{K}\sum_{k=1}^K{\bm g} ({\bm z}^\top{\bm \xi}^k)\right)+
(
{\bm c}_{[N-1]}^{M})^\top\left(\frac{1}{K}\sum_{k=1}^K{\bm g} ({\bm z}^\top {\bm \xi}^k)\right) \\
{\rm s.t.} \quad~ &  \widetilde{v}_i\in [-\underline{v}_i,\bar{v}_i], ~ i\in [N-1],\label{eq:bound}\\
& {\bm e}^\top \widetilde{\bm v}_{[N-1]}=0,\\
& -\widetilde{v}_i-{c}^{M}_i\leq 0, ~i\in [N-1],\\
& h_m \widetilde{\bm v}_{[N-1]}^{\top}
 \left(G_{\bm A} (r_{1}^{m},r_3^{m},p^{m})-G_{\bm B} (r_2^{m})\right) \nonumber \\
 &
 +h_m ({\bm c}_{[N-1]}^{M})^{\top}
 \left(G_{\bm A} (r_{1}^{m},r_3^{m},p^{m})-G_{\bm B} (r_2^{m})\right)\leq 0, \;\; m\in [M].
\end{align}
\end{subequations}
Introducing
dual variables $\underline{\bm \lambda}\in \R^{N-1}$, $\bar{\bm \lambda}\in \R^{N-1}$,
 $\beta\in \R$, ${\bm \eta}\in \R^{N-1}$ and ${\bm \theta}\in \R^M$,
the Lagrangian is
\bgeq
&& L(\widetilde{v}_1,\cdots,\widetilde{v}_{N-1}, \underline{\lambda}_1,\cdots,\underline{\lambda}_{N-1},\bar{\lambda}_1,\cdots,\bar{\lambda}_{N-1},\beta,\eta_1,\cdots,\eta_{N-1},\theta_1,\cdots,\theta_{M})\\
&=&\widetilde{\bm v}_{[N-1]}^\top\left(\frac{1}{K}\sum_{k=1}^K{\bm g} ({\bm z}^\top {\bm \xi}^k)\right)
+({\bm c}_{[N-1]}^{M})^\top\left(\frac{1}{K}\sum_{k=1}^K{\bm g} ({\bm z}^\top {\bm \xi}^k)\right)
+\underline{\bm \lambda}^\top (-\widetilde{\bm v}_{[N-1]})
+\underline{\bm \lambda}^\top(-\underline{\bm v})
\\
&& +\bar{\bm \lambda}^\top \widetilde{\bm v}_{[N-1]}-\bar{\bm \lambda}^\top\bar{\bm v}+\beta{\bm e}^\top \widetilde{\bm v}_{[N-1]}
+{\bm \eta}^\top(-\widetilde{\bm v}_{[N-1]}-{\bm c}_{[N-1]}^{M})\\
&& +\sum_{m=1}^M {\theta}_m \Big(h_m \widetilde{\bm v}_{[N-1]}^{\top}
 \left(G_{\bm A} (r_{1}^{m},r_3^{m},p^{m})-G_{\bm B} (r_2^{m})\right)
 + h_m ({\bm c}_{[N-1]}^{M})^{\top}
 \left(G_{\bm A} (r_{1}^{m},r_3^{m},p^{m})-G_{\bm B} (r_2^{m})\right)\Big)\\
&=& ({\bm c}_{[N-1]}^{M})^\top\left(\frac{1}{K}\sum_{k=1}^K{\bm g} ({\bm z}^\top {\bm \xi}^k)\right)
-\underline{\bm \lambda}^\top\underline{\bm v}
-\bar{\bm \lambda}^\top\bar{\bm v} 
-{\bm \eta}^\top {\bm c}_{[N-1]}^{M}\\
&& +\sum_{m=1}^M {\theta}_m h_m  ({\bm c}_{[N-1]}^{M})^{\top}
 \left(G_{\bm A} (r_{1}^{m},r_3^{m},p^{m})-G_{\bm B} (r_2^{m})\right)\\
&&
+ \widetilde{\bm v}_{[N-1]}^\top\left(\frac{1}{K}\sum_{k=1}^K{\bm g} ({\bm z}^\top {\bm \xi}^k)
-\underline{\bm \lambda}
+\bar{\bm \lambda} +\beta{\bm e}-{\bm \eta}+\sum_{m=1}^M {\theta}_m h_m \left(G_{\bm A} (r_{1}^{m},r_3^{m},p^{m})-G_{\bm B} (r_2^{m})\right)
\right).
\edeq
The resulting 
dual problem is 
\bgeq
& \displaystyle \max_{\underline{\bm \lambda}\geq 0,\bar{\bm \lambda}\geq 0,\beta,{\bm \eta}\geq 0,{\bm \theta}\geq 0} & ~({\bm c}_{[N-1]}^{M})^\top\left(\frac{1}{K}\sum_{k=1}^K{\bm g} ({\bm z}^\top {\bm \xi}^k)\right)
-\underline{\bm \lambda}^\top\underline{\bm v}
-\bar{\bm \lambda}^\top\bar{\bm v}
-{\bm \eta}^\top {\bm c}_{[N-1]}^{M}\\
&&
\quad +\sum_{m=1}^M {\theta}_m h_m ({\bm c}_{[N-1]}^{M})^{\top}
 \left(G_{\bm A} (r_{1}^{m},r_3^{m},p^{m})-G_{\bm B} (r_2^{m})\right) \\
&{\rm s.t.}& 
\frac{1}{K}\sum_{k=1}^K{\bm g} ({\bm z}^\top {\bm \xi}^k)
-\underline{\bm \lambda}
+\bar{\bm \lambda} +\beta{\bm e}-{\bm \eta}+\sum_{m=1}^M {\theta}_m h_m \left(G_{\bm A} (r_{1}^{m},r_3^{m},p^{m})-G_{\bm B} (r_2^{m})\right)=0.
\edeq
The proof is complete.
\hfill
\Box

\section{Proof of Theorem~\ref{thm:reformulate_n}}
\proof{Proof of Theorem~\ref{thm:reformulate_n}.}
We first show that the Lagrange dual of \eqref{eq:robust-pi0} is \eqref{eq:reformulation_conser}.
By replacing
${\bm v}_{[N-1]}-{\bm c}_{[N-1]}^{M}$
with a new variable  
$\widetilde{\bm v}_{[N-1]}$,
we can reformulate
problem \eqref{eq:robust-pi0}  
as 
\begin{subequations}
\begin{align}
\displaystyle    \min_{\widetilde{\bm v}_{[N-1]},{\bm \Pi}} \quad &
\widetilde{\bm v}_{[N-1]}^\top \left(\frac{1}{K}\sum_{k=1}^K{\bm g} ({\bm z}^\top{\bm \xi}^k)\right)+
(
{\bm c}_{[N-1]}^{M})^\top\left(\frac{1}{K}\sum_{k=1}^K{\bm g} ({\bm z}^\top {\bm \xi}^k)\right) \\
{\rm s.t.} \quad ~~&  \widetilde{v}_i\in [-({\bm \Pi}_i \sigma)\underline{v}_i,({\bm \Pi}_i \sigma)\bar{v}_i], ~ i\in [N-1],\label{eq:pi-bound}\\
& {\bm e}^\top \widetilde{\bm v}_{[N-1]}=0,\\
& -\widetilde{v}_i-{c}^{M}_i\leq 0, \;\;i\in [N-1],\\
& h_m \widetilde{\bm v}_{[N-1]}^{\top}
 \left(G_{\bm A} (r_{1}^{m},r_3^{m},p^{m})-G_{\bm B} (r_2^{m})\right) \nonumber \\
 &
 +h_m ({\bm c}_{[N-1]}^{M})^{\top}
 \left(G_{\bm A} (r_{1}^{m},r_3^{m},p^{m})-G_{\bm B} (r_2^{m})\right)\leq 0,\;\; m\in [M],\\
  &0\leq {\bm \Pi}_{i,j}\leq 1,\;\;\forall i,j\in [N-1],\\
 &{\bm \Pi}{\bm e}={\bm e},
 {\bm \Pi}^\top{\bm e}={\bm e}.
\end{align}
\end{subequations}
Introducing
dual variables $\underline{\bm \lambda}\in \R^{N-1}$, $\bar{\bm \lambda}\in \R^{N-1}$,
 $\beta\in \R$, ${\bm \eta}\in \R^{N-1}$, ${\bm \theta}\in \R^M$, $\Theta\in \R^{(N-1)\times (N-1)}$,
 $\underline{\bm \tau} \in \R^{N-1}$ and $\bar{\bm \tau}\in \R^{N-1}$,
the Lagrangian is
\bgeq
&& L(\widetilde{v}_{[N-1]}, \Pi,\underline{\bm \lambda},\bar{\bm \lambda},\beta,{\bm\eta}, {\bm \theta}, {\bm \Theta}, \underline{\bm \tau}, \bar{\bm \tau})\\
&=&\widetilde{\bm v}_{[N-1]}^\top\left(\frac{1}{K}\sum_{k=1}^K{\bm g} ({\bm z}^\top {\bm \xi}^k)\right)
+({\bm c}_{[N-1]}^{M})^\top\left(\frac{1}{K}\sum_{k=1}^K{\bm g} ({\bm z}^\top {\bm \xi}^k)\right)
+\sum_{i=1}^{N-1} \underline{\lambda}_i(-({\bm \Pi}_i \sigma) \underline{v}_i- \tilde{v}_i)
\\
&& +\sum_{i=1}^{N-1} \bar{\lambda}_i(\tilde{v}_i-({\bm \Pi}_i \sigma) \bar{v}_i)+\beta{\bm e}^\top \widetilde{\bm v}_{[N-1]}
+{\bm \eta}^\top(-\widetilde{\bm v}_{[N-1]}-{\bm c}_{[N-1]}^{M})\\
&& +\sum_{m=1}^M {\theta}_m^\top\Big(h_m \widetilde{\bm v}_{[N-1]}^{\top}
 \left(G_{\bm A} (r_{1}^{m},r_3^{m},p^{m})-G_{\bm B} (r_2^{m})\right)
 + h_m ({\bm c}_{[N-1]}^{M})^{\top}
 \left(G_{\bm A} (r_{1}^{m},r_3^{m},p^{m})-G_{\bm B} (r_2^{m})\right)\Big)\\
 &&+\sum_{i,j} {\bm \Theta}_{ij} (\Pi_{ij}-1)+\underline{\bm \tau}^\top({\bm \Pi}{\bm e}-{\bm e})
 +\bar{\bm \tau}^\top({\bm \Pi}^\top{\bm e}-{\bm e})\\
&=& ({\bm c}_{[N-1]}^{M})^\top\left(\frac{1}{K}\sum_{k=1}^K{\bm g} ({\bm z}^\top {\bm \xi}^k)\right)
-{\bm \eta}^\top {\bm c}_{[N-1]}^{M} +\sum_{m=1}^M {\theta}_m h_m  ({\bm c}_{[N-1]}^{M})^{\top}
 \left(G_{\bm A} (r_{1}^{m},r_3^{m},p^{m})-G_{\bm B} (r_2^{m})\right)\\
 &&-\sum_{i,j} {\bm \Theta}_{ij}-\underline{\bm \tau}^\top {\bm e}-\bar{\bm \tau}^\top {\bm e}
 +\langle (-\underline{\bm \lambda} \circ \underline{\bm v}-\bar{\bm \lambda} \circ \bar{\bm v})\sigma^\top+\underline{\bm \tau}
{\bm e}^\top+{\bm e} \bar{\bm \tau}^\top+{\bm \Theta},{\bm \Pi}\rangle\\
&&+ \widetilde{\bm v}_{[N-1]}^\top\left(\frac{1}{K}\sum_{k=1}^K{\bm g} ({\bm z}^\top {\bm \xi}^k)
-\underline{\bm \lambda}
+\bar{\bm \lambda} +\beta{\bm e}-{\bm \eta}+\sum_{m=1}^M {\theta}_m h_m \left(G_{\bm A} (r_{1}^{m},r_3^{m},p^{m})-G_{\bm B} (r_2^{m})\right)
\right).
\edeq
Then the dual problem is given by \eqref{eq:reformulation_conser}.
To complete the proof, we can derive the Lagrange dual of the inner minimization problem in \eqref{eq:defi_new_control-a},
which also yields
\eqref{eq:reformulation_conser}. We omit the details.

\hfill 
\Box

\section{Sensitivity Analysis of Problem 
\texorpdfstring{$\max_{\bm z\in Z}\rho_{\gamma}^{\mathrm{B}}(\bm z)$}{}}

In this section, we investigate the changes in the optimal value of problem 
$\max_{\bm z\in Z}\rho_{\gamma}^{\mathrm{B}}(\bm z)$, whose MILP reformulation 
is given in \eqref{eq:reformulation_conser}. Specifically, the vector 
${\bm \sigma}$ in \eqref{eq:sigma0} is generated by the function 
$\varsigma(t):=\gamma t^{\gamma-1}$ and is replaced by 
${\bm \sigma}\circ {\bm e}_s$, where 
${\bm e}_s=(\overbrace{1,\cdots,1}^{s},0,\cdots,0)\in \R^{N-1}$.
Accordingly, ${\bm \sigma}$ is replaced by 
${\bm \sigma}\circ {\bm e}_s$ in reformulation 
\eqref{eq:reformulation_conser}, and the number of nonzero components 
of the resulting weight vector is controlled by 
$s\in\{1,\ldots,N-1\}$. This construction combines replacing the binary 
weights $y_i\in\{0,1\}$ with continuous weights $\sigma_i\in(0,1)$ 
and controlling the number of nonzero components through 
$s\in\{1,\ldots,N-1\}$ instead of $\Gamma\leq N-1$.

\begin{figure}[!ht]
\begin{minipage}[t]{0.32\linewidth}
\centering
\includegraphics[width=2.2in]{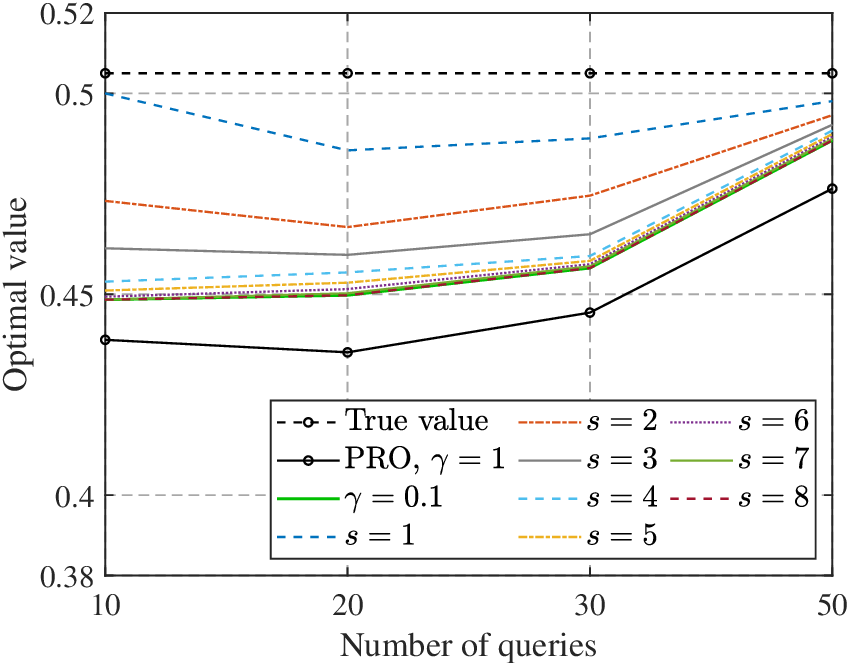}
\text{\footnotesize{\hspace{0.9cm}\textbf{(a)}}}
\end{minipage}
 \hfill
\begin{minipage}[t]{0.32\linewidth}
\centering
\includegraphics[width=2.2in]{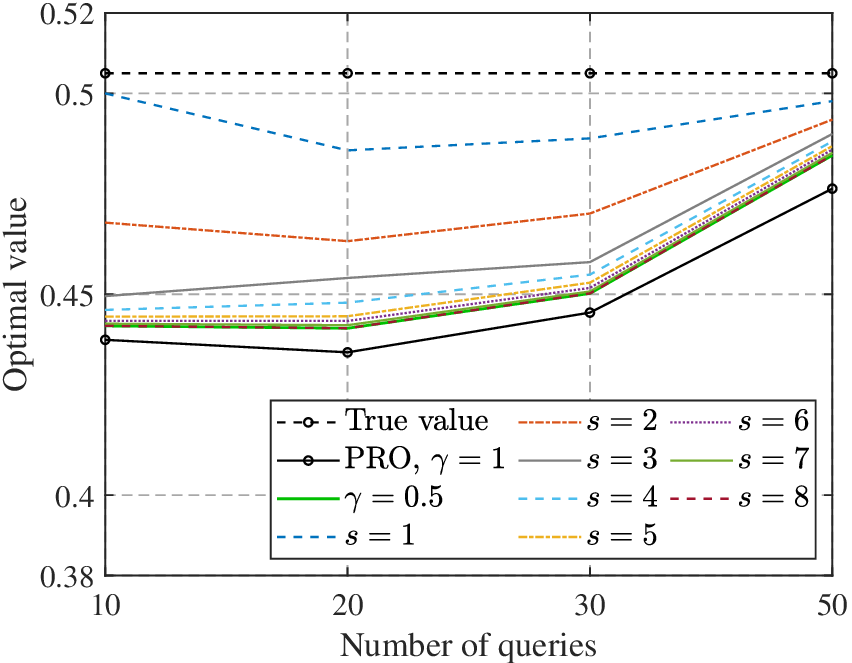}
\text{\footnotesize{\hspace{0.9cm}\textbf{(b)}}}
\end{minipage}
\hfill
\begin{minipage}[t]{0.32\linewidth}
\centering
\includegraphics[width=2.2in]{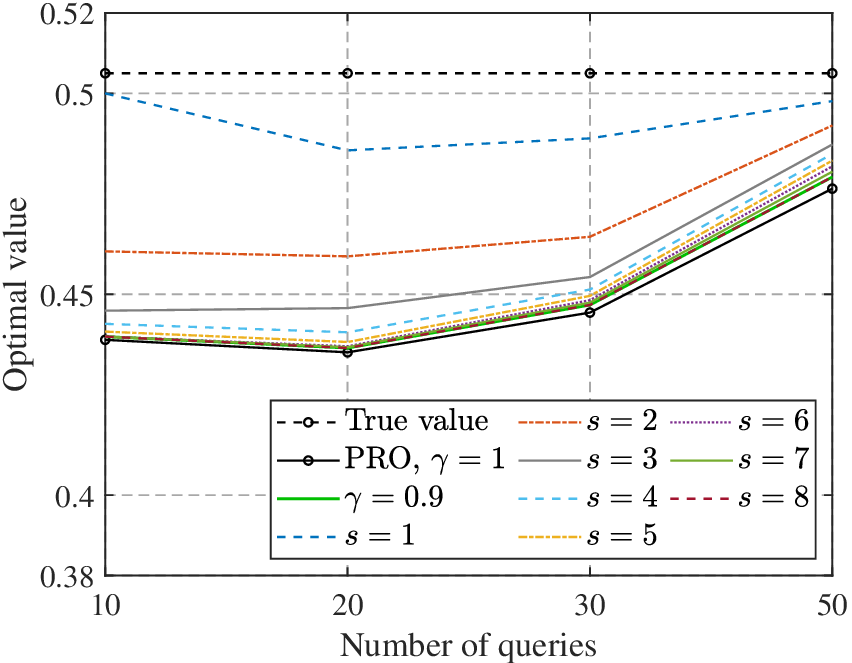}
\text{\footnotesize{\hspace{0.9cm}\textbf{(c)}}}
\end{minipage}
\caption{\footnotesize{
Optimal values of problem~\eqref{eq:PRO-PLA-N-dis-Gamma} under different values of $\gamma$ and number of 
queries $M$: 
(a)--(c) correspond to $\gamma=0.1,0.5,0.9$, respectively.
In each panel, the curves report the optimal values for different values of $s$
as $M$ varies, together with the true-utility benchmark and the result for 
$\gamma=1$.
}}
\label{fig:optimal values_compare}
\vspace{-0.2cm}
\end{figure}

Figure~\ref{fig:optimal values_compare}(a) depicts the optimal values
for different $s\in\{1,\ldots,N_{10}-1\}$ with $\gamma=0.1$ and
$M\in\{10,20,30,50\}$.
Figures~\ref{fig:optimal values_compare}(b)--(c) present the
corresponding results for $\gamma=0.5$ and $0.9$, respectively.
We observe that, for each fixed $\gamma$, the differences among the
optimal values become smaller as $s$ increases.
Moreover, as $\gamma$ increases from $0.1$ to $0.5$ and $0.9$,
the resulting optimal values gradually approach the benchmark
corresponding to preference robust optimization without conservatism
control, i.e., $\gamma=1$ and ${\bm \sigma}={\bm e}$, under which
problem $\max_{\bm z\in Z}\rho_{\gamma}^{\mathrm{B}}(\bm z)$ reduces
to the robust problem \eqref{eq:PRO-PLA-N-dis-2-mm}.

\begin{figure}[!h]
\begin{minipage}[t]{0.32\linewidth}
\centering
\includegraphics[width=2.1in]{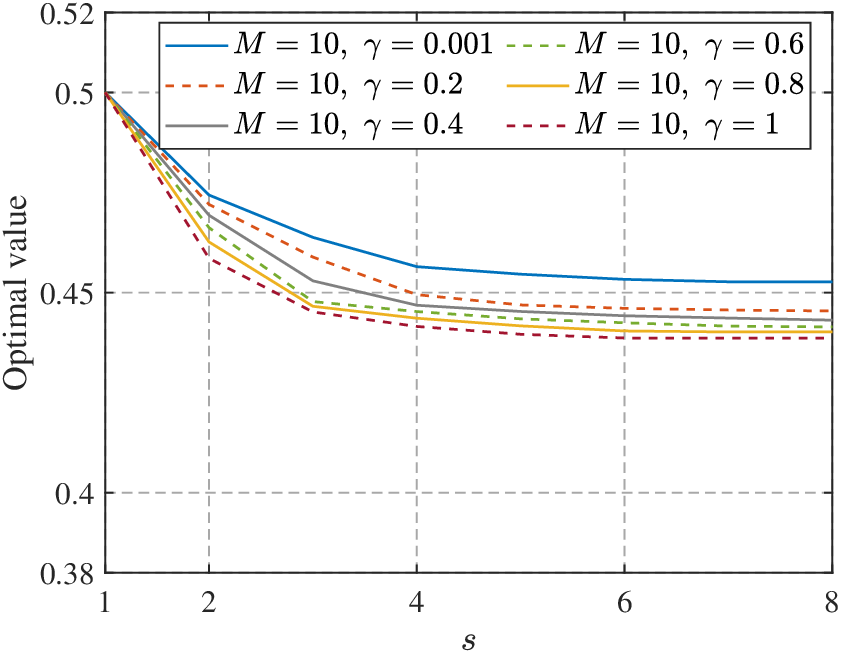}
\text{\footnotesize{\hspace{0.6cm}\textbf{(a)}}}
\end{minipage}
\hfill
\begin{minipage}[t]{0.32\linewidth}
\centering
\includegraphics[width=2.1in]{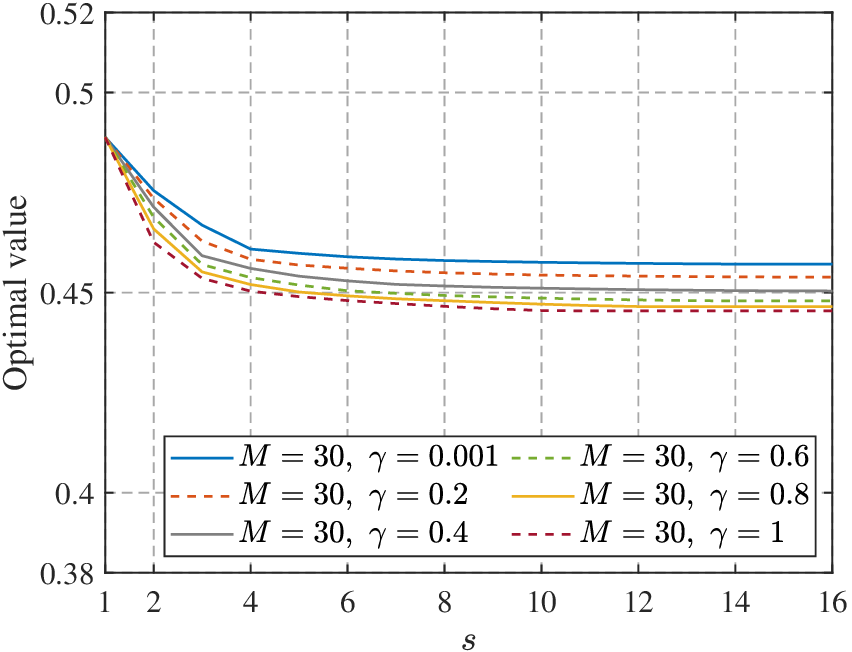}
\text{\footnotesize{\hspace{0.6cm}\textbf{(b)}}}
\end{minipage}
\hfill
\begin{minipage}[t]{0.32\linewidth}
\centering
\includegraphics[width=2.1in]{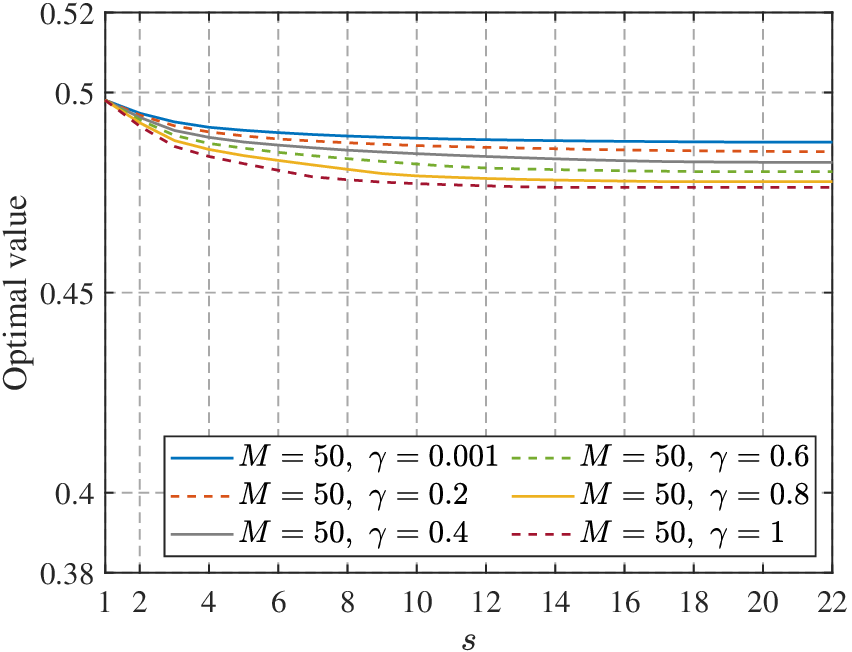}
\text{\footnotesize{\hspace{0.6cm}\textbf{(c)}}}
\end{minipage}
\caption{\footnotesize{
Optimal values of problem~\eqref{eq:PRO-PLA-N-dis-Gamma} as $s$ varies under different values 
of $M$ and $\gamma$: 
(a) $M=10$, (b) $M=30$, and (c) $M=50$.
Each panel reports the optimal values for 
$\gamma\in\{0.001,0.2,0.4,0.6,0.8,1\}$.
}}
\label{fig:optimal values_compare2}
\end{figure}

Further,
Figure~\ref{fig:optimal values_compare2}(a) shows the changes in the
optimal values with respect to $s$ for
$\gamma\in\{0.001,0.2,0.4,0.6,0.8,1\}$ when $M=10$.
Likewise, Figures~\ref{fig:optimal values_compare2}(b) and (c)
report the corresponding results for $M=30$ and $M=50$, respectively.
We observe that, for all considered values of $\gamma$, the optimal value increases as the number of queries $M$ increases from $20$ to $40$, while it decreases as the number $s$ of nonzero components of ${\bm \sigma}$ increases from $1$ to $N_{M}-1$. In addition, the differences among the optimal values associated with different values of $\gamma$ and $s$ become smaller as $M$ increases, indicating that the effects of these parameters become less pronounced as more preference information is elicited.

\section{Modified Polyhedral Method with Response Error}
\label{sec:inconsistenty}

Throughout this paper, we assume that the DM's responses are fully consistent with his or her underlying preferences. Although response errors and inconsistencies are beyond the scope of the present study, the proposed polyhedral method can be adapted to accommodate such situations with appropriate modifications.  To this end, we draw on the ideas of \cite{BeO13}, who propose a self-correcting mechanism for handling inconsistent responses and elicitation errors in multi-attribute linear utility assessment.
We can incorporate this
technique into the nonlinear utility case.
To illustrate
the technique,
suppose that we have already asked the DM the comparison question and received the response  
${\bm A}_m\succeq {\bm B}_m$,
and construct the inequality $\bbe[u_N^*({\bm A}_m)] -\bbe[u_N^*({\bm B}_m)]\geq 0$.
To account for the possibility of inconsistencies or errors in response,
we can introduce a binary variable $\phi_m$ for question $m$
and consider 
\bgeq
&& ({\bm R}{\bm v})^\top (G_{{\bm A}_m}(r_1^m,r_3^m,p^m)-G_{{\bm B}_m}(r_2^m)) +(N-1+\epsilon)\phi_m \geq \epsilon, \\
&& ({\bm R}{\bm v})^\top (G_{{\bm A}_m}(r_1^m,r_3^m,p^m)-G_{{\bm B}_m}(r_2^m)) +(N-1-\epsilon)\phi_m \leq N-1.
\edeq
In both constraints, if $\phi_i=0$,
 the first constraint enforces the correct response and the second constraint is redundant.
If $\phi_i=1$,
then it either introduces an inconsistency or makes a response error,
where the first constraint becomes redundant and the second constraint forces the inequality to flip.
If we only expect the user to respond incorrectly to 
a small fraction of the questions, we can use the constraint $\sum_{l=1}^m\phi_l\leq r m$,
where $r\in (0,1)$ is a parameter that indicates the maximum fraction of responses that we allow to be incorrect. We envisage that problem \eqref{eq:ana_center} for the analytic center and problems \eqref{eq:gene_m_Q} for optimal strategy will 
still be mixed-integer linear programs.

\end{document}